\documentclass[12pt]{article}
\usepackage{amsmath,amssymb,amsfonts,latexsym,amsthm}
\usepackage{amsmath,amssymb,amsxtra,epsf,amscd,graphics,color,array,ulem,etex}
\usepackage[all]{xy}
\usepackage{amsmath,amssymb,amsxtra,epsf,amscd,graphics,color,array,ulem,etex,pstricks,tikz,pgf,texgraph}
\usepackage[all]{xy}
\usepackage{mathrsfs}
\usepackage[dvips]{geometry}
\usepackage{array}
\usepackage{epsf}
\usepackage{epsfig}
\usepackage{verbatim,ifthen, fancyvrb,xkeyval,xcolor,epic,eepic,rotating}

\usepackage{mathrsfs}
\usepackage[dvips]{geometry}
\usepackage{array}
\usepackage{epsf}
\usepackage{epsfig}
\newtheorem*{lemma}{Lemma}
\newtheorem*{prop}{Proposition}
\newtheorem*{thm}{Theorem}
\newtheorem*{cor}{Corollary}

\newcommand{\ad}{\operatorname{ad}}
\newcommand{\nc}{\newcommand}
\nc{\Spec}{\operatorname{Spec}}
\nc{\Index}{\operatorname{Index}}
\nc{\GKdim}{\operatorname{GKdim}}
\def\q{\mathfrak q}
\def\h{\mathfrak h}
\def\g{\mathfrak g}
\def\a{\mathfrak a}
\def\p{\mathfrak p}
\def\s{\mathfrak s}
\def\l{\mathfrak l}

\title{ADAPTED PAIRS IN TYPE A AND REGULAR NILPOTENT ELEMENTS 
\footnote{Work supported in part by Israel Science Foundation Grant No.  710724.
The main content of this paper was the subject of a four hour lecture course given by F. Fauquant-Millet 
at the  Israel Science Foundation Research Workshop On Orbits Primitive Ideals And Quantum Groups,
held on 3-7 March 2013 in the Weizmann Institute of Science.
The taped lectures may be viewed on the site http://www.wisdom.weizmann.ac.il/$\sim$joseph/orbits/.}}

\author{Florence Fauquant-Millet 
\footnote{email : florence.millet@univ-st-etienne.fr}
and Anthony Joseph
\footnote{email : anthony.joseph@weizmann.ac.il}}

\date{March 2015}

\begin{document}

\maketitle


\begin{abstract}

Let $\mathfrak g$ be a simple Lie algebra over an algebraically closed field $\bf k$ of characteristic zero and $\bf G$ its adjoint group. A biparabolic subalgebra $\mathfrak q$ of $\mathfrak g$ is the intersection of two parabolic subalgebras whose sum is $\mathfrak g$.
The algebra $Sy(\mathfrak q)$ of semi-invariants on $\mathfrak q^*$ of a proper biparabolic subalgebra $\mathfrak q$ of $\mathfrak g$ is polynomial in most cases, in particular when $\mathfrak g$ is simple of type $A$ or $C$.
On the other hand $\mathfrak q$ admits a canonical truncation $\mathfrak q_{\Lambda}$ such that $Sy(\mathfrak q)=Sy(\mathfrak q_{\Lambda})=Y(\mathfrak q_{\Lambda})$ where 
$Y(\mathfrak q_{\Lambda})$ denotes the algebra of invariant functions on $(\mathfrak q_{\Lambda})^*$.
An adapted pair for $\mathfrak q_{\Lambda}$ is a pair $(h,\,\eta)\in \mathfrak q_{\Lambda}\times(\mathfrak q_{\Lambda})^*$ such that $\eta$ is regular in $(\mathfrak q_{\Lambda})^*$ and $(ad\,h)(\eta)=-\eta$.
In [A. Joseph, Slices for biparabolic coadjoint actions in type A. J. Algebra 319 (2008), no. 12, 5060-5100] adapted pairs for every truncated biparabolic subalgebra $\mathfrak q_{\Lambda}$ of a simple Lie algebra $\mathfrak g$ of type $A$ were constructed and then
provide Weierstrass sections  for $Y(\mathfrak q_{\Lambda})$ in $(\mathfrak q_{\Lambda})^*$.
These Weierstrass sections are linear subvarieties $\eta+V$ of $(\mathfrak q_{\Lambda})^*$ such that the restriction map induces an algebra isomorphism of $Y(\mathfrak q_{\Lambda})$ onto the algebra of regular functions on $\eta+V$.
The main result of the present work is to show that for each of the adapted pairs $(h,\,\eta)$ constructed in [A. Joseph, Slices for biparabolic coadjoint actions in type A. J. Algebra 319 (2008), no. 12, 5060-5100] one can express $\eta$ (not quite uniquely) as the image of a \textit{regular nilpotent} element $y$ of $\mathfrak g^*$ under the restriction map $\mathfrak g^* \rightarrow \mathfrak q^*$. This is a significant extension of [A. Joseph and F. Fauquant-Millet, Slices for biparabolics of index one. Transformation Groups 16 (2011), no. 4, 1081-1113], which obtains this result in the rather special case of a truncated biparabolic of index one.  Observe that $y$ must be a $\textbf{G}$ translate of the standard regular nilpotent element defined in terms of the already chosen set $\pi$ of simple roots.  Consequently one may attach to $y$ a unique element of the Weyl group $W$ of $\mathfrak g$.  Ultimately one can then hope to be able to describe adapted pairs (in general, that is not only for $\mathfrak g$ of type $A$) through the Weyl group.

\end{abstract}


\



\section{Introduction.}\label{00}

Throughout this paper the base field $\bf{k}$ is assumed to be algebraically closed of characteristic zero.

\subsection{Basic Goals.}\label{00.0}

Let $\mathfrak g$ be a simple Lie algebra and $\textbf{G}$ the corresponding adjoint group.  A fundamental result of Chevalley asserts that the algebra $Y(\mathfrak g)$ of $\textbf{G}$ invariant regular functions on $\mathfrak g^*$ is polynomial. A key component of the proof was the use of the Weyl group $W$.  A further fundamental result of Kostant gave what we call (\ref {00.1}) a Weierstrass section $y+V$ for $Y(\mathfrak g)$ in $\g^*$.  Identifying $\mathfrak g^*$ with $\mathfrak g$ through the Killing form, the Kostant construction was made possible through the existence of a ``principal s-triple" $(x,\,h,\,y)$ taking $V=\mathfrak g^x$ in the above and comparing the $\ad h$ eigenvalues on $\mathfrak g^x$ with the degrees of homogeneous generators of $Y(\mathfrak g)$.

The Weyl group and the principal s-triple are very special to simple Lie algebras.  Thus it would seem almost inconceivable that similar results could hold for other families of Lie algebras.  Yet we showed that Chevalley's theorem extends to most (truncated - see \ref {00.3}) parabolic \cite {FJ0} (and indeed biparabolic \cite {J0.1}) subalgebras of $\mathfrak g$.  This provided large families of Lie algebras admitting invariant algebras which are polynomial.  In this it should be stressed that the invariant generators are even for $\mathfrak {sl}_n(\bf k)$ given by complicated expressions of high degree which can only be explicitly described in some very special cases.

One meaning of a Weierstrass section is that gives a ``linearization" of the invariant generators. Thus its existence can provide their better understanding.  At first these were constructed by brute force, literally parabolic by parabolic; but eventually a general though extremely complicated construction was found \cite {J1} for type $A$.  A disappointment was that in type $C$ a Weierstrass section may fail to exist.

In \cite {J2} we found that in the case of the (truncated) Borel $\mathfrak b_\Lambda$ a slight modification of the invariant algebra admits a Weierstrass section for all $\mathfrak g$ simple.  This fact, by no means trivial, uses some special features of the invariant generators.  Its full significance has yet to be understood.   It involves a very special choice of a principal nilpotent element $y$ whose image in $\mathfrak b_\Lambda^*$ eventually gives rise to a Weierstrass section.  This brings the Weyl group into the picture, as $y$ is determined by a very special element $w \in W$ relative to the standard ``Jordan" form of $y$ fixed by our choice of Borel. Roughly speaking $w$ was found to be ``the square root" of the longest element (see \cite [Thm. 3.2]{J2}). Moreover we found \cite [Thm. 9.4]{J2} that it was possible to thereby obtain a Weierstrass section for $\mathfrak b_\Lambda \subset \mathfrak g$ under the condition that a simple Lie algebra of type $C_2$ is not the penultimate element of the Kostant cascade for $\mathfrak g$, explicitly that the latter is not of type $C, B_{2n},F_4$.

The above circumstance raised the question as to whether our previous construction of a Weierstrass section for a truncated biparabolic in type $A$ could be expressed in terms of taking an image of a regular nilpotent element of the ambient Lie algebra.  This question was settled \cite {FJ1} for the very simplest case of a truncated biparabolic of index $1$.  This is just the case when the biparabolic is the parabolic subalgebra whose Levi factor has two blocks of relative coprime size.  We call it the coprime two-block case. The analysis was purely combinatorial and rather complicated, particularly involving meanders. Notably these have received considerable attention \cite {DG}, \cite {M}, \cite {MY} recently in studying related problems. A meander, as the name suggests, is globally indescribable.  For example it is practically impossible to decide if a meander passing through one point will eventually return to pass through one of its nearest neighbours.  This feature is a serious stumbling block to proving any result needing such a description.  To be specific, recall that a biparabolic subalgebra in $\s\l_n(\bf k)$ is Frobenius, that is of index $0$, exactly when \cite {DG} (see also \cite [Appendix]{J0.5}) the meander defined by the Kostant cascades associated to the choice of biparabolic is a single ``edge" (in the language of Section \ref {1}) passing through every integer point on the line $[1,n]$.  One example is provided by the coprime two-block case.  However a similar general description and even a ``closed" formula for the number of Frobenius biparabolics in $\s\l_n(\bf k)$ seems quite unattainable.

In the present work we settle the general question of whether all the Weierstrass sections constructed in \cite {J1} for truncated biparabolics in type $A$ can be obtained through a regular nilpotent element. This is far more complex than the coprime two-block case not only because arbitrary choices of blocks are involved, but also because the Kostant cascade has to be altered to obtain the required meanders.

As in our previous work \cite{FJ1} a fundamental notion is that of a turning point of a meander.  Roughly speaking we must eliminate the internal turning points and join the subsequent ``straightened" meanders.  A huge simplification in what could otherwise be a problem of insurmountable difficulty is that we need only the local description of a meander as it passes through a ``component" of the double partition (see \ref {0.4}) associated to the biparabolic in question.  A further huge simplification comes from the subdivision of turning points into sources and sinks.

Unlike the case when the truncated biparabolic $\mathfrak q_\Lambda$  is $\mathfrak g$ itself, or of index $1$, a Weierstrass section does not necessarily meet every regular co-adjoint orbit.  This leads to a number of open problems which can be expressed in terms of the nullcone $\mathscr N(\q_{\Lambda})$ (see \ref{00.2})
of $\q_{\Lambda}$.  For example is $\mathscr N(\mathfrak q_\Lambda)$ equidimensional?  If so, does every regular orbit of $\mathscr N(\mathfrak q_\Lambda)$ admit a ``minimal representative" (see \ref {00.4}) and more precisely can the minimal representative be chosen in the form described in 
\cite{J1}?  In other words does the construction in 
\cite{J1} recover all regular co-adjoint orbits in $\mathscr N(\mathfrak q_\Lambda)$?  These questions have a positive answer if $\mathscr N(\mathfrak q_\Lambda)$ is irreducible; but irreducibility generally fails.

Again the nullcone $\mathscr N(\mathfrak q_\Lambda)$ may admit components with no regular elements \cite[11.3]{J1}.
Outside type $A$ the situation is even worse, $\mathscr N(\mathfrak q_\Lambda)$ may admit no regular element.

A general question which can be posed for all $\mathfrak g$ semisimple is whether the image of $\mathscr N(\mathfrak g)_{reg}$ under the restriction map contains $\mathscr N(\mathfrak q_\Lambda)_{reg}$ (see \ref{00.2} for a definition)?  Even if this were true the image contains far too much, so a further question is to find a systematic way to cut down the image so that just $\mathscr N(\mathfrak q_\Lambda)_{reg}$ is obtained.

 Indeed in the coprime two-block case we find (\ref {10.3}) that \textit{every} regular element is an image of a regular \textit{nilpotent} element of $\mathfrak g^*$.  This does not seem to have a very easy proof (for example that it is a consequence of the parabolic being Frobenius).   Again there does not seem to be an obvious reason why this should be true in general.

\subsection{Weierstrass Sections.}\label{00.1}

Let $X$ be a finite dimensional vector space. Identify the algebra of regular functions ${\bf k}[X]$ on $X$ with the symmetric algebra $S(X^*)$ of the dual $X^*$ of $X$. A subvariety of $X$ of the form $y+V$ with $y \in X$ and $V$ a subspace of $X$ is called a linear subvariety.  Let $A$ be a subalgebra of $S(X^*)$.  A Weierstrass section for $A$ in $X$ is a linear subvariety $y+V$ of $X$ such that the restriction of $A$ to $y+V$ induces an algebra isomorphism of $A$ onto ${\bf k}[y+V]$.  One may remark that the latter identifies with $S(V^*)$ and so is a polynomial algebra and quotient of $S(X^*)$.

The above terminology was introduced by the Russian school (see \cite {P}, \cite {PV} for example), it being noted that Weierstrass canonical form for elliptic curves  can be expressed through the existence of a Weierstrass section \cite [2.2.2.1] {P}, \cite [Sect. 2]{J3}.

Let $\mathfrak a$ be a finite dimensional algebraic Lie algebra over $\textbf {k}$. Let $Y(\mathfrak a)$ denote the algebra of invariant polynomial functions on $\mathfrak a^*$ which we identify as a subalgebra of $S(\mathfrak a)$.   In principle one could use the existence of a Weierstrass section for $Y(\mathfrak a)$ in $\mathfrak a^*$ to show that $Y(\mathfrak a)$ is polynomial, though such a procedure is seldom practical.  Instead of this  one attempts to construct a Weierstrass section for $Y(\mathfrak a)$ given that $Y(\mathfrak a)$ is polynomial.

\subsection{Adapted Pairs.}\label{00.2}

Continue to assume that the Lie algebra $\a$ is defined as in \ref{00.1}.
One calls an element  $f \in \mathfrak a^*$ regular if the dimension  of the stabilizer $\mathfrak a^f:=\{a \in \mathfrak a\mid (\ad a)(f)=0\}$ has minimum possible dimension called the index of $\mathfrak a$, denoted by $\ell(\a)$. Here $\ad$ designates  the co-adjoint action of $\a$ on $\a^*$.  The $\bf k$-vector space $\a^f$ is a subalgebra of $\a$. Let $\mathfrak a^*_{reg}$ denote the set of regular elements of $\mathfrak a^*$.  It is an open dense subset.  Let $\bf{A}$ be the adjoint group of $\mathfrak a$. Then the regular elements of $\mathfrak a^*$ are just those generating $\bf{A}$ orbits of maximal dimension.

An adapted pair for $\mathfrak a$ is a pair $(h,\,\eta)\in\mathfrak a\times\mathfrak a^*_{reg}$ such that 
$(\ad h)(\eta)=-\eta$.  Since $\mathfrak a$ is algebraic we can assume that $h$ is ad-semisimple and \textit{this we shall always do}.
Observe that $\bf{A}$ acts by simultaneous conjugation on the set of adapted pairs,
giving then an equivalence relation on this set.

Let $\mathscr N(\a)$ be the zero locus in $\a^*$ of the ideal of $S(\a)$ generated by the space $Y(\a)_+$ of homogeneous vectors of $Y(\a)$ with positive degree. The variety $\mathscr N(\a)$ may be called the nullcone of $\a$. We set $\mathscr N(\a)_{reg}=\mathscr N(\a)\cap\a^*_{reg}$.

Let $(h,\,\eta)$ be an adapted pair for $\a$ and $f\in Y(\a)_+$ homogeneous of degree $d$. Since $f$ is $\ad\,h$ invariant one has, for all $c\in{\bf k}$, $f(\eta)=f(exp(c\ad\,h)(\eta))=f(exp(-c)(\eta))$ since $(\ad\,h)(\eta)=-\eta$. Thus $f(\eta)=(exp(-c))^df(\eta)$
and $\eta$ must belong to $\mathscr N(\mathfrak a)$.

 A semi-invariant of $S(\mathfrak a)$ is defined to be the generator of a non-trivial one dimensional representation of $\mathfrak a$ under adjoint action.  Denote the space generated by all semi-invariants in $S(\mathfrak a)$  by $Sy(\mathfrak a)$.  
 
 Assume that $Sy(\mathfrak a)=Y(\mathfrak a)$ and that $Y(\mathfrak a)$ is polynomial and 
 let $(h,\,\eta)$ be an adapted pair for $\a$ and
$V$ be an $\ad h$ stable complement of $(\ad \mathfrak a)(\eta)$ in $\a^*$.
  Then by \cite [Cor. 2.3] {JS}, $\eta +V$ is a Weierstrass section for $Y(\mathfrak a)$ in $\a^*$ (shortly a Weierstrass section for $\a$).  Moreover in this case every element of  $\eta+V$ is regular (see \cite [7.8(ii)]{J2}, for example).

Adapted pairs need not exist.  It is not at all easy to find adapted pairs if they do exist, nor is it to classify their $\bf{A}$ orbits.  If $Sy(\mathfrak a)=Y(\mathfrak a)$ and is polynomial, a formal description of such pairs is given by \cite[Prop. 2.1.2]{FJ3}.


\subsection{Truncation.}\label{00.3}

Define $\mathfrak a$ as in \ref {00.1}.  Let $\mathfrak h$ be a Cartan subalgebra for $\mathfrak a$.   In general $S(\mathfrak a)$ may admit proper semi-invariants.  Let $\mathfrak a_\Lambda$ denote the common kernel of the corresponding one dimensional representations in $\mathfrak a$.   It is an algebraic Lie algebra and an ideal in $\mathfrak a$, called the canonical truncation of $\mathfrak a$.  Moreover $\mathfrak h_\Lambda:= \mathfrak h \cap \mathfrak a_\Lambda$ is a Cartan subalgebra for $\mathfrak a_\Lambda$. One has $Sy(\mathfrak a)=Y(\mathfrak a_\Lambda)=Sy(\a_{\Lambda})$, a result referred to as a generalization of a lemma of Borho in \cite{RV}.

Thus in the following we shall assume that $Sy(\mathfrak a)=Y(\mathfrak a)$, equivalently that $S(\mathfrak a)$ admits no proper semi-invariants.  This further implies by a well-known result of Chevalley-Rosenlicht that the growth rate dimension or Gelfand-Kirillov dimension $\GKdim Y(\mathfrak a)$ of $Y(\mathfrak a)$ equals $\ell(\mathfrak a)$. (It turns out \cite [Prop. 4.1]{OM} that $\mathfrak a$ being algebraic is not needed here !)

Let $\mathfrak g$ be a semisimple Lie algebra.  A biparabolic subalgebra $\mathfrak q$ of $\mathfrak g$ is defined to be the intersection of two parabolic subalgebras whose sum is $\mathfrak g$. Identify $\mathfrak g^*$  with $\mathfrak g$ through the Killing form.
Call $y \in \mathfrak g^*$ nilpotent if it is ad-nilpotent as an element of $\mathfrak g$.

Let $\mathfrak q_\Lambda$ be a truncated biparabolic subalgebra of a semisimple Lie algebra and let $\textbf{Q}_\Lambda$ be its adjoint group.  It  was shown in \cite{FJ0} and \cite{J0.1} that in most cases $Y(\mathfrak q_\Lambda)$ is a polynomial $\textbf{k}$-algebra.  In particular this is true for any truncated biparabolic subalgebra of a simple Lie algebra of type $A$ or $C$.

\subsection{Presentation of Adapted Pairs.}\label{00.4}

In \cite {J1}, it was shown in type $A$ that $\mathfrak q_\Lambda$ admits at least one and possibly several orbits of adapted pairs.  It should be stressed that in the latter the second element $\eta$ of an adapted pair was presented in a very special form, namely as a sum of root vectors running over a subset $S$ of roots with the property that the restriction of $S$ to $\mathfrak h_\Lambda$ is a basis for $\mathfrak h^*_\Lambda$. Of course few elements in the $\textbf{Q}_\Lambda$ orbit of $\eta$ will have such a presentation, which we call minimal. 

For an algebraic Lie algebra $\mathfrak a$ (and even for a truncated biparabolic) we do not
 know if every $\textbf{A}$ orbit of regular elements in $\mathscr N(\mathfrak a)$ admits such a minimal presentation.

Assume now that $\a=\q_{\Lambda}$ and $Y(\q_{\Lambda})$ is polynomial. Then (\cite[Rem. 2.2.2 (1)]{FJ3}) the first element $h\in\h_{\Lambda}$ of an adapted pair for $\q_{\Lambda}$  is uniquely determined by the second one, namely $\eta$.
Moreover when the dimension of the $-1$ eigenspace of $\ad h_{\mid_{\mathfrak q^*_\Lambda}}$ has dimension $\dim \mathfrak h_\Lambda$ then every second element $\eta$ of an adapted pair for $\q_{\Lambda}$ with $h\in\h_{\Lambda}$ as its first element is of minimal presentation (\cite[2.3.3]{FJ3}).

\subsection{Introducing the Weyl Group.}\label{00.5}

As discussed in \ref {00.0} it is interesting to bring the Weyl group into the description of adapted pairs by attempting to express the second element $\eta$ of an adapted pair of a truncated biparabolic subalgebra $\q_{\Lambda}$ of a simple Lie algebra $\g$ of type $A$ as the image of a $\it{regular}$ nilpotent element $y$ of $\mathfrak g^*$, under the restriction map $\rho:\mathfrak g^* \rightarrow \mathfrak  q^*_\Lambda$.

It seemed at first that this problem must be 
of almost insurmountable difficulty, since adapted pairs form an enormous and rather varied family.   Nevertheless in this paper we resolve it for all the adapted pairs described in 
\cite{J1}. The combinatorics which extends over many pages is rather intricate, yet has a structure which on careful study can be seen to be intrinsic if not essentially canonical.  In the next section we give an overview of the proof.


\section{Outline of Proof.}

\subsection{Link Patterns and Meanders.}\label{0.0}

Let $n$ be a positive integer and $[1,\,n]$ the set $\{1,\,2,\,\ldots,\,n\}$ viewed as $n$ points on a horizontal line $L$. A link pattern is a set of arcs joining the points of $L$, 
the end-points of each arc being images under an involution $\sigma$.  In fact we shall consider a pair of involutions $\sigma^+,\,\sigma^-$ giving rise to a pair of link patterns drawn (for aesthetic reasons) above and below $L$.  A meander is defined to be an orbit of the group $<\sigma^+,\,\sigma^->$.  This notion plays a fundamental role in our analysis.

\subsection{The subset $S$.}\label{0.1}

A biparabolic subalgebra $\mathfrak q:=\mathfrak p^- \cap \mathfrak p^+$ of 
$\s\l_n(\bf k)$ in standard form is specified by a double ordered partition $\mathscr  J^\pm$ of $[1,n]$.  A directed arc above (resp. below) $L$ joining distinct points of a fixed component of 
$\mathscr J^-$ (resp. $\mathscr J^+$) of cardinal $>1$ corresponds to a root vector in $ \mathfrak p^-$ (resp. $\mathfrak p^+$).

Through the Killing form on $\mathfrak g$ we may identify $\mathfrak q^*$ with the opposed subalgebra $\q^-$ of $\mathfrak q$ in $\mathfrak g$.  Then the starting point of our construction is to write the second element $\eta$ of an adapted pair of $\q_{\Lambda}$ in the form
$$\eta =\sum_{\alpha \in S}x_\alpha,$$
 for some subset $S$ of the roots of $\mathfrak q^-$.  Notice that this means that we can consider $\eta$ as an element of $\mathfrak g^*$.

 Recall that the Kostant cascade for  
$\s\l_n(\bf k)$ is the set of positive roots defined by the arcs of the Dynkin diagram involution. Similarly there are involutions $\kappa^\pm$ defined by the Dynkin diagram involutions associated to the components of $\mathscr  J^\pm$.  The corresponding arcs form a set which is a union of Kostant cascades.  This does not quite give an appropriate choice for $S$, but in \cite [Sect. 5]{J1} it was shown how to define $S$ by modifying this choice. In Sections \ref {1} - \ref {3}, we describe this construction in a slightly more general form using what we call anti-Toeplitz involutions.  In more detail we select parts of the Dynkin involutions $\kappa^\pm$ according to the classification of generators of $Y(\mathfrak q_\Lambda)$ and replace these parts by anti-Toeplitz involutions.

A key fact is that $S$ has in particular the following two properties.

The first property of $S$ is that it is a disjoint union of two subsets $S^+$ and $S^-$.  The elements of $S^+$ and $S^-$ are represented by directed arcs above and below $L$ and correspond to root vectors in 
$(\p^+)^*$ and $(\p^-)^*$ respectively.
 \textit{In this no two arcs in $S^+$ (resp. $S^-$) meet the same point of $L$}.  This latter condition means that we may view $S^+$ (resp. $S^-$) as being exactly given by an involution $\sigma^+$ (resp. $\sigma^-$) of $[1,\,n]$ which in particular interchanges the two end-points of each arc.  This property was shared by the union of the Kostant cascades defined by $\kappa^\pm$.

A second property of $S$ is that its restriction to $\mathfrak h_\Lambda$ is a basis of $\mathfrak h^*_\Lambda$, Proposition \ref {3.2}. This property was \textit{not} shared by the union of the Kostant cascades.  It implies in particular that the elements of $S$ are linearly independent.  This has the consequence that every $<\sigma^+ ,\,\sigma^->$ orbit is in fact a $<\sigma^+ \sigma^->$ orbit and as such is called an edge $E$.  We may consider that $S$ is described by the arcs of a union of edges. 

We remark that these two properties are not sufficient to ensure that $\eta$ is a regular element of $\mathfrak q^*_\Lambda$.  However this is so if the anti-Toeplitz involutions are of anti-Coxeter or anti-Jordan type \cite [Sect. 8]{J1}.  In particular we do not know if our present construction gives rise to any new (equivalence classes of) adapted pairs.

\subsection{Turning Points.}\label{0.2}

First we remark that every root occurring in $S^+$ (resp. $S^-$) is a positive (resp. negative) root.  Consequently the arrow it defines should point from left to right (resp. right to left).
Thus as we pass along a given edge the arrows will not be all aligned.  This means that the set of roots corresponding to root vectors of  
$\q^*$ in a given edge, will not form a Dynkin diagram of type $A$.   Had this been the case and if in addition there were to be just one single edge passing through all the points of $[1,\,n]$, then $\eta$ viewed as an element of $\mathfrak g^*$ would already be regular nilpotent.

Let $S_E$ be the subset of $S$ defined by the arcs of an edge $E$.  Let $\eta_E$ be the sum of root vectors whose weights are the elements of the subset $S_E$. Then we need to find an element $y_E$ of $\mathfrak g^*$ which maps to $\eta_E$ under the restriction map $\rho:\mathfrak g^* \rightarrow \mathfrak q_\Lambda^*$.   Call a source (resp. sink) a point on $L$ in which arrows depart (resp. arrive).   The sources and sinks of $E$ form its turning points. Eliminating all internal sources and sinks of $E$, is called a straightening of the edge $E$, if the arcs of the resulting edge $E^*$ defines an element $y_E \in \mathfrak g^*$ with the same image as $\eta_E$.   Let $e$ denote the number of vertices of $E$.  The straightened edge $E^*$ corresponds to a root system of type $A_{e-1}$.

\subsection{Intervals.}\label{0.3}

Define a simple interval $I$ of $E$ to be any interval of $E$ lying between consecutive turning points. Define a compound interval $I$ to be the join of several simple intervals and to be odd if their number is odd. The value $\iota_I$ of a (compound) interval $I$ is defined to be the root defined by the arc joining the starting and finishing points of $I$.  Let $K$ denote the set of roots such that the corresponding root vectors lie in the kernel of the restriction map $\rho$.  In order to straighten $E$ we must show that $\iota_I  \in  K \cup -K$. An interval with this property is called nil.  A  fundamental observation (which significantly simplified our original proofs) is that an \textit{odd} compound interval must start at a source and end at a sink or vice-versa.

\subsection{Components of the Double Partition.}\label{0.4}

  A key point in our analysis is that we only have to describe how a given edge passes through sets of the form $J^+ \cap J^-$ with $J^+ \in \mathscr J^+$,  $J^- \in \mathscr J^-$, called a component of the double partition (or a double component).

 Observe that an interval fails to be nil if and only if its end-points lie in the same double component.

 Suppose that the centres $c^\pm$ of $J^\pm$ are distinct.  In this case we show that all the turning points in $J^+ \cap J^-$ are either all sources or all sinks. (This is an easy consequence of our partial replacement of the Dynkin involutions by anti-Toeplitz coming from single blocks - called anti-Toeplitz blocks). Consequently an odd compound interval cannot both start and end at such a component.

 When $c^+=c^-$, we call the component $J^+\cap J^-$ equicentral.  In this case it is possible to describe rather explicitly how an edge behaves in $J^+\cap J^-$.  Here we can assume that $J^+\neq J^-$ as otherwise we are reduced to the ``trivial" case when $\mathfrak q =\mathfrak g$.  Then we show that this intersection admits at most one source and one sink with both being end-points of an edge (not necessarily the same one) which cannot be confined in the intersection (\ref {3.8.4}).  Thus the only bad edge is one in which both its end-points lie in the same equicentral component.

  In fact in the above we have to be a little more careful about signs.  Thus we assign a value $\epsilon_I \in \{\pm 1\}$ to an odd interval $I$ and show (Lemma \ref {4.3}) that (apart from the above exception) one has $\epsilon_I\iota_I \in -K$.

\subsection{Straightening the edges.}\label{0.5}

In order to straighten an edge $E$ we must construct a system $\Pi^*_E$ of type $A_{e-1}$ such that

(a).  Every element of $S_E$ lies in $\mathbb N \Pi^*_E$ and

(b). Every element of $\Pi^*_E$ which is not in $S_E$ lies in $K$.

Then  $E^*$ defined by the elements of $\Pi^*_E$ is the required straightened edge.

The difficulty in the above construction comes from condition (b) above.   In order that it be satisfied we need to find sufficiently many elements lying in $K$.  This is where we need our preliminary material showing that an odd interval value is nil.

Given the result described in the last paragraph of \ref {0.4}, we obtain the required edge straightening in the following easy and elegant manner, representing a major simplification of the method used in \cite{FJ1}  for the coprime two-block case. At each internal turning point we delete one directed line meeting it. Outside the biparabolic case one may always choose the line to represent a root vector of the Levi factor of  
$\q^*$. Suppose this turning point is a sink (resp. source). Then we replace this directed line by an arc, with an arrow pointing in the same direction, joining the second point meeting the deleted line, to a source (resp. sink). In this there is a very simple algorithm to avoid that two such arcs arrive at the same turning point. This construction and the fact that for an odd interval $I$ one has $\epsilon_I \iota_I \in -K$, ensures that (b) above is satisfied. The original undeleted lines and the new arcs define 
$\Pi_E^*$.  Property (a) results from a very simple order relation on arcs.  The arcs of the straightened edge define $\Pi^*_E$.

Details are given in Section \ref {6}.

\subsection{Joining straightened edges.}\label{0.6}

This presented little difficulty in the coprime two-block case. Here new ideas are needed.

Let $E^*$ be the straightening of the edge $E$ and let $y_{E}$ denote the sum of the root vectors formed from the roots in $\Pi^* _E\cup S_E$.  Then condition (b) of \ref {0.5} implies that $\eta_E$ is the image of $y_{E}$ under the restriction map $\rho$.   Condition (a) implies that $y_{E}$ is conjugate to a Jordan matrix of size $e$.   In the case when $e=n$ it follows that $y_{E}$ is regular nilpotent.  However this numerical condition is practically never satisfied, indeed $S$ itself is given by a union of edges.  Thus $y:=\sum y_{E}$ is just a sum of Jordan blocks.  To obtain just one block of size $n$ we must add further root vectors to $y$, and to ensure that the resulting element has the same image as $y$ under $\rho$, these root vectors must lie in $K$.    The construction is in Sections \ref {8} and \ref {9}.   A key point in the proof is that a straightened edge value $\iota_{E^*}$ lies in $-K$.  (In this one shows that in the bad edge $E$ mentioned above in \ref{0.4}, $E$ is never equal to its straightening.)   Then a relatively simple and elegant algorithm allows one to join straightened edges. In the parabolic case this consists of aligning edges by the natural order on their start points, and then sequentially joining end and start points along the alignment.  A little extra effort is needed to deal with edges reduced to a single point and for the biparabolic case.

Details are given in Section \ref {9}.

\section{The Basic Combinatorial Objects.}\label{1}

\subsection{Arcs.}\label{1.1}

Fix $n$ an integer $>1$ and for all $1\le i,\,j\le n$, denote by $x_{i,\,j}$ the $n$-order matrix  with $1$ on the intersection of the $i^{th}$ row and the $j^{th}$ column.
 Let $\h$ be the span of the $x_{j,\,j}-x_{j+1,\,j+1}$ with $1\le j\le n-1$. It is a Cartan subalgebra of $\g=\s\l_n(\bf k)$.

Adopt the Bourbaki notation \cite [Planche I]{B} for the root system of $(\g,\,\h)$. In particular denote by $(\varepsilon_i)_{1\le i\le n}$, the dual basis of the basis $(x_{j,\,j})_{1\le j\le n}$ 
and still denote by $\varepsilon_i$ the restriction to the Cartan subalgebra $\h$.

Given $r\leq s$ positive integers, set $[r,\,s]:=\{r,\,r+1,\,\ldots,\,s\}$, which we refer to as a connected subset of the integers. Call $c_{[r,\,s]}:=\frac{(r+s)}{2}$, the centre of $[r,\,s]$.

We view $[1,\,n]$ as a set of $n$ equally spaced points lying on a horizontal line $L[1,\,n]$ (or simply $L$).

Take $r,\,s \in [1,\,n]$.  If $r \neq s$ we denote by $(r,\,s)$ the root $\varepsilon_r-\varepsilon_s$.  Diagrammatically we present $(r,\,s)$ as a (non-trivial) directed arc joining the points $r,\,s \in L$ with an arrowhead pointing from $r$ to $s$ placed on the arc.  If $r<s$ (resp. $r>s$), the root is said to be positive (resp. negative). 
All the directed arcs above (resp. below) $L$ will represent positive (resp. negative) roots 
(this rule is conserved up to the end of Section \ref {5}).
Let $\Delta^+$ (resp. $\Delta^-$) denote the set of positive (resp. negative) roots and set $\Delta=\Delta^+\cup \Delta^-$.

If $r=s$ then the point $r$ on $L$ may also be represented by an arc joining $r$ to itself, called a trivial arc.\par
If $(r,\,s)$ is a directed arc above or below $L$, then the points $r$ and $s$ are called the end-points of this arc. We may also say that the arc $(r,\,s)$ meets the point $r$ or meets the point $s$ or joins the points $r$ and $s$.
More precisely if the directed arc $(r,\,s)$ is non-trivial, then $r$ is called the starting point and $s$ the finishing point of the directed arc $(r,\,s)$.

We may also view the simple root  $\alpha_i=\varepsilon_i-\varepsilon_{i+1}$ as the half-integer point at $i+\frac{1}{2}$ on $L$.  Set $\pi:=\{\alpha_i\}_{i=1}^{n-1}$.

A root $\gamma \in \Delta$ may be written uniquely as a linear combination of the elements of $\pi$.  Those which occur with non-zero coefficient we refer to as forming the support of $\gamma$, denoted as $\text {supp} \ \gamma$. Given $r,s \in [1,\,n]$ with $r\neq s$, the support of an arc joining $r,s$ on $L$  is defined to be $\text {supp} \ (\varepsilon_r-\varepsilon_s)$ and in particular the empty set if $r=s$.

Let $x_\gamma$ denote a root vector of weight $\gamma \in \Delta$.  If $\gamma = \varepsilon_i-\varepsilon_j$, then $x_\gamma$ may be represented as the $n$-order matrix $x_{i,j}$.

\subsection{Link Patterns.}\label{1.2}

A set of arcs joining each point of $[1,\,n]$ on $L[1,\,n]$ to a unique second point on $L[1,\,n]$ (possibly the same one) and lying above (resp. below) $L$ is called a link pattern.  It may be identified with an involution $\sigma^+$ (resp. $\sigma^-$) of $[1,\,n]$ ($\sigma^\pm$ may coincide or not with the involutions $\sigma^\pm$ defined in Section \ref{2} but from Section \ref{2} up to the end of Section \ref{5} they will coincide with them.)

Observe that giving $\sigma^+$ (resp. $\sigma^-$) is also equivalent to giving a set of pairwise orthogonal positive (resp. negative) roots.  These are represented by arcs above (resp. below) $L$ with arrows pointing from left to right (resp. right to left). This convention will be kept up to the end of Section \ref{5}.

\subsection{Linear order.}\label{1.2bis}

Recall that an ordered set $(S,\,\leq)$ is said to be linearly (resp. strictly linearly) ordered if for all distinct elements $s,s^\prime \in S$ one has either $s\leq s^\prime$ or $s^\prime \leq s$, or both (resp. but not both).


\color{black}
\subsection{Meanders, Loops and Edges.}\label{1.3}

Retain the notation of \ref {1.2}. Consider the group $\Gamma:=<\sigma^+,\,\sigma^->$, as a Coxeter group with length function.

Let $g^+(t)$ (resp. $g^-(t)$) be the unique element of $\Gamma$ of length $t-1$ starting (for $t>1$) at $\sigma^+$ (resp. $\sigma^-$) on the right.

A meander is a $\Gamma$-orbit of a point in $[1,\,n]$.  It is called an edge if it is also a $<\sigma^+\sigma^->$ orbit and a loop otherwise. 
An edge reduced to one point of $[1,\,n]$ (which is fixed under both $\sigma^+$ and $\sigma^-$) is called a trivial edge.
A  non-trivial edge $E$ has two points fixed under either $\sigma^+$ or $\sigma^-$, one we call its starting point $a$ and one we call its finishing point $b$.
Observe that (when $E$ is non-trivial) $\varepsilon_a-\varepsilon_b$ is a root, called 
the value of the edge $E$ and denoted by $\iota_E$.

Given an edge $E$ (with starting point $a$ and finishing point $b$) we set $e = |E|$.  Define a function $\varphi\,: \,[1,\,e] \mapsto E$ as follows.  Suppose that $a$ is a fixed point of $\sigma^+$ (resp. $\sigma^-$) and set $\varphi(t)=g^-(t)a$ (resp. $\varphi(t)=g^+(t)a$) for $t>1$ and $\varphi(1)=a$.

Observe that $\varphi(e)=b$.

Viewing the points of $E$ as a subset of $[1,\,n]$ this defines a \textit{new} linear order on this subset.

For each $t\in[1,\,e]$, the neighbour(s) of $\varphi(t)$ are defined as to be the points $\varphi(t-1)$ and $\varphi(t+1)$ (when they exist).

\subsection{Turning Points, Sources and Sinks.}\label{1.4}

The elements of a meander and in particular of an edge lie in $[1,\,n]$ and so can be viewed as integers.

A turning point of an edge $E$ (with starting point $a$ and finishing point $b$, possibly equal if $E$ is trivial) is defined to be either an element of $\{a,\,b\}$ or one of the form $\varphi(i)\,:\,i \in [2,\,e-1]$ (if $e\ge 3$) with $\varphi(i)-\varphi(i-1)$ having a different sign to $\varphi(i)-\varphi(i+1)$.
In other words, a turning point of the edge $E$ is either a point of $E$ fixed under $\sigma^\pm$ (that is equal to $a$ or $b$), or a point $t\in E$ such that $t$ is not a fixed point under $\sigma^\pm$ and $t-\sigma^+(t)$ has a different sign to $t-\sigma^-(t)$. 

 The latter is called an internal turning point of $E$ and the former an end-point of $E$.

Let $T$ (resp. $T_0$) denote the set of turning (resp. internal turning) points of $E$.

Given a sequence $(i_1,\,i_2),\,(i_2,\,i_3),\,\ldots,\,(i_k,\,i_{k+1})$ of directed non-trivial arcs, then we say that the arrows on these arcs are aligned since they all point from $i_j$ to $i_{j+1}$\, : \,$j=1,\,2,\,\ldots,\, k$. We may also say that the arrows on the directed arcs $(i_k,\,i_{k+1}),\,(i_{k-1},\,i_k),\ldots,\,(i_2,\,i_3),\,(i_1,\,i_2)$ are aligned.

Suppose that $\varphi(t)$ is not an end-point of $E$. Now place arrowheads on the arcs defined by the elements of $\sigma^\pm$ as prescribed in \ref{1.1},\,\ref {1.2}, that is by drawing the arcs defined by $\sigma^+$ (resp. $\sigma^-$) above (resp. below) $L$ with an arrowhead pointing from left to right (resp. from right to left). Then the arrows on the two directed arcs meeting $\varphi(t)$ are aligned
 if and only if $\varphi(t)$ is not a turning point.  
 
 A turning point (internal or not) of a non-trivial edge $E$ is called a source (resp. sink) of $E$ if the arc(s) point away from (resp. towards) the turning point.
 
 Equivalently for $t\in[1,\,e]$, $\varphi(t)$ is a source (resp. a sink) of $E$ if $\varphi(t)$ is the starting (resp. finishing) point of the directed arc(s) meeting $\varphi(t)$. 

Observe that sources and sinks alternate as one passes along the points $\varphi(t)\,:\,t\in [1,\,e]$ of the edge $E$ with $t$ increasing.

\subsection{Intervals and Interval Values.}\label{1.5}

Let $\varphi(s),\,\varphi(t)\in T$ be turning points of an edge $E$ of cardinality $e>1$ with $1\le s<t\le e$. The subset $I_{s,\,t}:=\{s,\,s+1,\,\ldots,\,t-1\}=[s,\,t-1]$ is called an interval.
If $\varphi(t)$ is the immediate successor to $\varphi(s)$ in $T$, with respect to the new linear order defined by $\varphi$, see \ref{1.3},
it is called a simple interval.  Otherwise it is called a compound interval.

Set $\beta_i=\varepsilon_{\varphi(i)}-\varepsilon_{\varphi(i+1)}\,:\, i \in [1,\,e-1]$.
Then $\beta_i$ is called a neighbouring value (or a neighbour) of $\varphi(i)$ and of $\varphi(i+1)$ and, if $i\in[1,\,e-2]$, $\beta_i$ and $\beta_{i+1}$ are said to be neighbouring values.

Let $I=I_{s,\,t}$ ($s<t$) be an interval.

The sum
$$\iota_I:=\sum_{i\in I}\beta_i,$$
is called a simple (resp. compound) interval value if $I$ is
simple (resp. compound). 
The points $\varphi(s)$ and $\varphi(t)$ are called the end-points of the interval $I$,
more precisely $\varphi(s)$ is called the starting point of $I$ and $\varphi(t)$ its finishing point.

 The set $\text {Supp} \ \iota_I:=\{\beta_i\,: \,i \in I\}=: \text {Supp} \ I$ is
called the $\beta$-support of $\iota_I$ or of $I$. It is a subset of $\Delta$.  Set $\Phi(I)=\{\varphi(j)\,|\,j \in [s,\,t]\}\subset [1,\,n]$. The interval $I$ (or interval value $\iota_I$) is called internal if $\varphi(s),\,\varphi(t)\in T_0$.


\subsection{Dynkin Diagram Involutions.}\label{1.6}

Fix $x \in \mathbb R$ and let $\sigma_x$ denote the involution $\sigma_x\,:\,y\mapsto 2x-y$ of $\mathbb R$.  In what follows we shall always assume that $x \in \frac {1}{2} \mathbb Z$.  Moreover for any connected subset $[r,\,s]$ of the integers, we define the restriction $\sigma_x$ to $[r,\,s]$ to be $\sigma_x(j)$ if $j,\sigma_x(j) \in [r,\,s]$ and the identity otherwise.  Then $c_{\sigma_a}:=a$, is called the centre of $\sigma_a$.  It is a fixed point of the restriction of $\sigma_a$ to $[r,\,s]$ if $a \in [r,\,s]$.

Given a connected subset $J=[r,s]$, of $[1,\,n]$, with centre $c_J$, see \ref{1.1}, we define $\kappa_J$ to be the identity on $[1,\,n]\setminus J$ and the restriction to $J$ of $\sigma_{c_J}$ on $[r,\,s]$.  Of course $\kappa_{[1,\,n]}$ is just the diagram automorphism of $\s\l_n(\bf k)$ which reverses order in the Dynkin diagram.

\subsection{Partitions.}\label{1.7}
If $[1,\,n]=\sqcup J_i$ (disjoint union) where each $J_i$ is a connected subset of $[1,\,n]$, we call 
$\mathscr J=\{J_i\}$ a partition of $[1,\,n]$
(here we may include singleton subsets of $[1,\,n]$). 
Every connected subset $J_i$ in $\mathscr J$ is called a connected component of $\mathscr J$. Let $c_i$ denote the centre of $J_i$.  Define the involution $\kappa_{J_i}$ on $[1,n]$ as the restriction of $\sigma_{c_i}$ to $J_i$, extended by the identity on $[1,n]\setminus J_i$.  These involutions commute for different $i$ and we let $\kappa_\mathscr J$ denote their product which is again an involution. This involution (and the corresponding arcs)  will be referred to as the integer involution (integer arcs) defined by $\mathscr J$.

Each connected set $J_i$ of cardinality $>1$ defines a non-empty subset $\pi_i$ of $\pi$ through the short central lines in $L$ joining nearest neighbours (cf. \ref {1.1}).  
By analogy with \ref{1.1} we may call $\pi_i$ the support of the set $J_i$. Of course one may view the (sum of the) corresponding root vectors as forming a Jordan block of size $|J_i|$. 
 The root vectors defined by the arcs of $\kappa_{J_i}$ lie on the anti-diagonal above the main diagonal of this block and form the Kostant cascade $\mathscr K_i$ in this block.  Set $\mathscr K =\cup \mathscr K_i$.

Now in fact we shall choose two partitions 
$\mathscr J^\pm=\{J_i^\pm\}$
 of $[1,n]$, let $\kappa_{\mathscr J^\pm}$ (or simply, $\kappa^\pm$) be the corresponding involutions and $\pi_i^\pm$ the corresponding subsets of $\pi$ (for $J_i$ such that $\mid\!J_i\!\mid>1$).  
 Set $\pi^\pm= \cup\pi_i^\pm$.  
The subsets $\pi_i^\pm$ are called the connected components or subsets of $\pi^\pm$.The arcs defined by $\kappa^+$ (resp. $\kappa^-$) are viewed as lying above (resp. below) $L$.

The pair $\mathscr J^\pm$ is referred to as a double partition of $[1,\,n]$.

\subsection{Biparabolic Subalgebras.}\label{1.8}

To the subset $\pi^\prime$ of $\pi$, we associate the corresponding standard parabolic subalgebras $\mathfrak p^\pm_{\pi^\prime}$ of 
$\g=\s\l_n(\bf k)$ whose roots form the subsets $\Delta^\pm_{\pi^\prime}:=\Delta^\pm \cup (\Delta \cap \mathbb Z \pi^\prime)$.

Fix a double partition $\mathscr J^\pm$ of $[1,\,n]$ and let $\pi^+,\,\pi^-$ denote the corresponding pair of subsets of $\pi$.
Set $\mathfrak q_{\pi^+,\,\pi^-}=\mathfrak p^-_{\pi^+}\cap \mathfrak p^+_{\pi^-}$. It is called a biparabolic subalgebra of $\mathfrak g$. We may also write it as $\mathfrak q_{\mathscr J^+,\,\mathscr J^-}$.

Our standing hypothesis throughout this paper is that $\pi^+\cup \pi^- =\pi$. Were this to fail we could have just replaced $\pi$ by $\pi^+\cup \pi^-$. With this convention $\hat{\pi}^\pm:=\pi \setminus \pi^\pm = \pi^\mp \setminus (\pi^+ \cap \pi^-)$ and so $\hat{\pi}^+\cap \hat{\pi}^-=\emptyset$.

 Unless otherwise stated we shall also assume throughout that $\pi^+ \cap \pi^- \varsubsetneq \pi$, which exactly excludes the ``trivial" case $\mathfrak q_{\pi^+,\,\pi^-}=\mathfrak g$.

We use the convention that the biparabolic algebra $\mathfrak q$ assigned to the double partition  $\mathscr J^\pm$ means $\mathfrak q_{\mathscr J^-,\,\mathscr J^+}$ with its dual $\mathfrak q^*$ identified with $\mathfrak q_{\mathscr J^+,\,\mathscr J^-}$ (which is also the opposed subalgebra $\q^-$ of $\q$) through the Killing form on 
$\s\l_n(\bf k)$.

Let 
$\mathfrak q_{\Lambda}$ be the canonical truncation of $\mathfrak q$ (\ref{00.3}).

Assume that we have constructed an adapted pair $(h,\,\eta)$ for $\mathfrak q_{\Lambda}$, as in \cite {J1}.  In this $\eta \in \mathfrak q_{\Lambda}^*$ and can be written in the form
$$\eta = \sum_{\alpha \in S}x_\alpha,\eqno {(*)}$$
for some $S$ in the set of roots of $\q^-$ and $h\in\h_{\Lambda}$.

 Observe that $R=(\mathbb N\pi^+\cup-\mathbb N\pi^-)\cap \Delta$ is the set of roots of $\mathfrak q_{\mathscr J^+,\,\mathscr J^-}$ and then that $-R$ is the set of roots of $\mathfrak q_{\mathscr J^-,\,\mathscr J^+}$. In particular we must have $S \subset R$ in $(*)$ above. Set $K=\Delta\setminus R$.  The kernel of the map $\mathfrak g^* \rightarrow \mathfrak q_{\mathscr J^-,\,\mathscr J^+}^*$ is just the span of the root vectors $x_\alpha\,: \,\alpha \in K$.  Set $R_*=(R\cup-R)\setminus R=-R \setminus (-R\cap R)$. The nilradical of $\mathfrak q_{\mathscr J^-,\,\mathscr J^+}$ is the span of the root vectors of $R_*$, so the latter is additively closed in $\Delta$.  Clearly $R_* \subset K$.  Equality holds in the parabolic case, that is if either $\pi^+$ or $\pi^-$ equals $\pi$.  Again $R\cup -R = \Delta$ in the parabolic case, but otherwise $R\cup -R$ may not be even additively closed in $\Delta$.

  One has $M:=R\cap -R=\mathbb Z(\pi^+\cap\pi^-)\cap \Delta$.  It is the set of roots of the Levi factor of $\mathfrak q$.  Thus the set of roots of $\mathfrak q$ is the disjoint union of $R_*$ and $M$.  Set $R_*^\pm=R_*\cap \Delta^\pm,\, K^\pm=K\cap \Delta^\pm$.  It is immediate that

\begin {lemma}

\
\begin{enumerate}

\item[(i)]  $R^+_* \subset \{\beta \in \Delta^+$ with a positive coefficient of some $\alpha \in \hat{\pi}^+=\pi^-\setminus(\pi^+\cap \pi^-)\} = K^+$.

\item[(ii)]  $R^-_*\subset \{\beta \in \Delta^-$ with a negative coefficient of some $\alpha \in \hat{\pi}^-=\pi^+\setminus(\pi^+\cap \pi^-)\} = K^-$.

\item[(iii)] $R_*\cap -R_*=\emptyset$.

\item[(iv)]  $\Delta=M \sqcup (K\cup -K)$.

\item[(v)] $R=M\sqcup (R\cap (K\cup -K))=M\sqcup(R\cap-K)$.

\item[(vi)] $-R_*= R\setminus M$.

\item[(vii)] $(M+K)\cap\Delta\subset K$.

\item[(viii)] $(M+(-K))\cap\Delta\subset -K$.

\item[(ix)] $(M+R_*)\cap\Delta\subset R_*$.

\item[(x)] $(M+(-R_*))\cap\Delta\subset -R_*$.

\end{enumerate}
\end {lemma}

\textbf{Remark}.   Recall that $K=R_*$ in the parabolic case.  A technical difficulty in the biparabolic case is that $K\cap -K$ need not be empty.
For example, in $\s\l_4(\bf k)$, if $\pi^+=\{\alpha_1,\,\alpha_3\}$ and $\pi^-=\{\alpha_2\}$ then $-(\alpha_1+\alpha_2+\alpha_3)\in K^-\cap(-K^+)\subset K\cap(-K)$.

\subsection{Components.}\label{1.9}

We call the pair $\mathscr J^\pm$ the double partition associated  to the biparabolic subalgebra $\mathfrak q_{\mathscr J^-,\mathscr J^+}$.
Recall that 
$\mathscr J^\pm=\{J_i^\pm\}$.  We call a non-empty
subset  of the form $J^+_i \cap J^-_j$ a component of the double partition or a double component. In this we shall generally drop subscripts. The analysis of the passage of a meander through a double component  is a central feature of our work.

Let $u,\,v\in[1,\,n]$. Observe that $\varepsilon_u-\varepsilon_v\in M$ if and only if $u$ and $v$ belong to the same component of the double partition.

\section{Modified Involutions.}\label{2}

\subsection{Anti-Toeplitz Involutions.}\label{2.1}

In order to construct an adapted pair (\ref{00.2}) it is necessary to modify the involutions $\kappa_{\mathscr J^\pm}$ defined in \ref {1.7}. This is carried out by first introducing a set $\mathscr M$ of markings of the Kostant cascade $\mathscr K$.  (Actually $\mathscr M$ is just a subset of $\mathscr K$; but we prefer this more picturesque terminology as the markings will be used on the figures.)

The above modification is the same for arcs above and arcs below $L$ so we shall just consider a single partition $\mathscr J$ and arcs above $L$.  Moreover the modification is carried out on each $J_i$ of cardinality $>1$ separately so in a first step we shall just assume $\mathscr J$ is reduced to $[1,\,n]$.

Set $\kappa:=\kappa^+$ and $\beta'_i = \varepsilon_i-\varepsilon_{\kappa(i)}\,: \, i \in [1,\,[n/2]]$. Then $\mathscr K = \{\beta'_i\}_{i=1}^{[n/2]}$, is the Kostant cascade.  The corresponding root vectors lie on the main antidiagonal of the $n\times n$ matrix defining $\s\l_n(\bf k)$ and above the main diagonal.  Fix a subset $\mathscr M \subset [1,\,[n/2]]$.  
We will say that (the arc representing) $\beta'_i\in\mathscr K$ is marked if $i\in\mathscr M$.
We shall also view (cf \ref{3.1}, \ref {3.6}) a marking as fixing a subset of $\overline{\pi}=\pi/\iota$ by assigning to $\beta'_j=\sum_{i=j}^{n-j}\alpha_i \in \mathscr K$, its ``end" simple roots $\alpha_j,\,\alpha_{n-j}$ identified through the involution $\iota$.

 Let $[r,\,s-1]$ be a connected component of $\mathscr M$ ($1\le r\le s-1\leq [n/2]$), that is that $[r,\,s-1]\subset\mathscr M$ and, if $s\leq [n/2]$, then $s\not\in\mathscr M$.
 If $n$ is even and $s-1=n/2$, set $s^-=s-1$, otherwise set $s^-=s$.
Then set $\mathscr M_{r,\,s}=[r,\,s^-]\cup \kappa([r,\,s^-])$ and call it an extended component of $\mathscr M$.

The restriction $\kappa_{\mathscr M_{r,\,s}}$ of $\kappa$ to $\mathscr M_{r,\,s}$ is an involution whose arcs correspond to root vectors in an anti-diagonal block. Conversely this anti-diagonal block defines $\mathscr M_{r,\,s}$. In \cite [5.2]{J1} we considered replacing part of this anti-diagonal block by an anti-Coxeter block (based on the notion of a Coxeter element) when $s-1<[n/2]$ and an anti-Jordan block (based on the notion of Jordan block) when $s-1=[n/2]$.  This replacement led to a new set of arcs some of which will cross and therefore not be linearly ordered. These crossings lead to some unavoidable technical complications.

For convenience consider the link pattern defined by the involution $\kappa_{\mathscr M_{r,\,s}}$ as a set of arcs drawn below $L$.

Then the new set of arcs (drawn above $L$) define a link pattern and hence an involution $\sigma_{\mathscr M_{r,\,s}}$ of $\mathscr M_{r,\,s}$.  A key property of the above construction is that the meander defined as the orbit of the group $<\kappa_{\mathscr M_{r,\,s}},\,\sigma_{\mathscr M_{r,\,s}}>$ acting on $\mathscr M_{r,\,s}$ is a single loop (resp. edge) for $s-1<[n/2]$
(resp. $s-1=[n/2]$).

We can take this further without introducing additional complications as follows. We start with a word of motivation.

Recall that an $m \times m$ Toeplitz matrix is one in which the entries of the $i^{th}$ and $(m-i)^{th}$ diagonals are all equal.  We define an $(m,\,i)$ anti-Toeplitz block as an $m \times m$ matrix block in which all entries are either zero, or equal one if they lie on the $i^{th}$ and $-(m-i)^{th}:i=0,1,\ldots,m-1$ anti-diagonals. Here the main anti-diagonal is taken to mean the $0^{th}$ anti-diagonal.

An anti-Toeplitz involution $\sigma_{\mathscr M_{r,\,s}}$ is an involution on $\mathscr M_{r,\,s}$ such that

$T(i)$.  The meander defined by $<\kappa_{\mathscr M_{r,\,s}},\,\sigma_{\mathscr M_{r,\,s}}>$ is a single loop (resp. edge) for $s-1<[n/2]$ (resp. $s-1=[n/2]$) passing through all points of $\mathscr M_{r,\,s}$.

$T(ii)$. The root vectors corresponding to its non-trivial arcs lie in the anti-diagonal block defined by $\mathscr M_{r,\,s}$, that is every non-trivial arc of the meander defined in $T(i)$ has one of its end-points in $[r,\, s^-]$ and its second end-point in $\kappa([r,\,s^-])$.

An example of an anti-Toeplitz involution (for $1\le r\le s-1<[n/2]$) obtains from an anti-Toeplitz block $(m,\,i)$ 
: $m=s+1-r\leq [n/2]$, with the non-zero entries defining the arcs of the involution, satisfying the following conditions. First that $i$ is coprime to $m$.  
 Second that it replaces the  antidiagonal block defining $\mathscr M_{r,\,s}$.  This is necessary and sufficient for $T(i)$ and $T(ii)$ to hold.  For $i=0$ this is just the antidiagonal block it is deemed to replace and for $i=1,\,m-1$ it is an anti-Coxeter block.
 
 When $s-1<[n/2]$, an anti-Toeplitz involution associated to an anti-Coxeter block is called an anti-Coxeter involution.\par


When $s-1<[n/2]$, the anti-Toeplitz involution $\sigma_{\mathscr M_{r,\,s}}$ has no fixed points (by $T(i)$).

Now consider the case $s-1=[n/2]$.  Here there is a degeneration corresponding to our previous replacement \cite [5.2]{J1} of an anti-Coxeter block by an anti-Jordan block.

In this case an anti-Toeplitz involution associated to an anti-Jordan 
 block is called an anti-Jordan involution.

Suppose $n$ is odd. Then $(n+1)/2$ is a fixed point of $\kappa_{\mathscr M_{r,\,s}}$ and consequently $\mathscr M_{r,\,s}$ has odd cardinality. Thus $\sigma_{\mathscr M_{r,\,s}}$ must have a fixed point and exactly one by $T(i)$.  Notice further that the edge (required by $T(i)$) must join the fixed point of $\kappa_{\mathscr M_{r,\,s}}$ to the fixed point of $\sigma_{\mathscr M_{r,\,s}}$ (of course in general not by a single arc).

Suppose $n$ is even.   
Then $\kappa_{\mathscr M_{r,\,s}}$ has no fixed point, thus by $T(i)$, $\sigma_{\mathscr M_{r,\,s}}$ must have two fixed points and hence defines one  non-trivial arc less than these defined by $\kappa_{\mathscr M_{r,\,s}}$. This will give 
\cite [Lemma 5.5]{J1} (see \ref{3.2}).

Finally observe that the subsets $\mathscr M_{r,\,s}$ corresponding to the connected components of $\mathscr M$ are disjoint.  Thus we may define an involution $\sigma$ as the product of the $\sigma_{\mathscr M_{r,\,s}}$
and of the restriction of $\kappa$ to $[1,\,[n/2]]\setminus\mathscr M$.
We may also (see \ref{2.4}) add subscripts $J$ and $\mathscr M$ to $\sigma$ to recall that it is defined with respect to $J$ (in this case $[1,\,n]$) and the marking $\mathscr M$. (Of course there are still many possible choices of the resulting $\sigma_{J,\,\mathscr M}$.) One may remark that the link pattern defined by the new involution $\sigma$ admits crossings but these are constrained within each extended component of $\mathscr M$.

\subsection{Batches of Arcs.}\label{2.2.1}

 Retain the notations of \ref {2.1} and in particular that $[r,\,s-1]$ is a connected component of the set $\mathscr M$ of markings.  Let $A_{r,\,s}$ (resp. $A^\prime_{r,\,s}$) be the set of non-trivial arcs defined by $\kappa_{\mathscr M_{r,\,s}}$ (resp. $\sigma_{\mathscr M_{r,\,s}}$).

To each non-trivial arc in $A_{r,\,s}$ (resp. in $A^\prime_{r,\,s}$) drawn below $L$ (resp. above $L$),  assign an arrow so that it represents a negative (resp. positive) root, as specified by \ref{1.2}. 
  
Then 
$T(i)$ ensures that the arrows are aligned along the meander defined by the arcs of $A_{r,\,s}$ and $A^\prime_{r,\,s}$.
\par
 Suppose the meander is a loop. Then the map $\theta\,:\, A^\prime_{r,\,s} \rightarrow A_{r,\,s}$ obtained by taking the successive element in the direction of the arrows, is independent of where we start and is bijective.\par
 Suppose the meander is an edge. Then the arrows define a starting point $a$ and 
 a finishing point $b$. If $n$ is even, then starting from the arc $\bigl(\kappa_{\mathscr M_{r,\,s}}(a),\,\sigma_{\mathscr M_{r,\,s}}(\kappa_{\mathscr M_{r,\,s}}(a))\bigr)$ and taking successive elements as before gives a bijection $\theta$ of $A^\prime_{r,\,s}$ onto $A_{r,\,s}\setminus (a,\,\kappa_{\mathscr M_{r,\,s}}(a))$. 
 If $n$ is odd, then there is a unique arc lying in $A^\prime_{r,\,s}$ whose end-point is $a$ or $b$. Then if this end-point is $a$, starting from $a$ and taking successive elements (in the direction of the arrows) gives a bijection of $A^\prime_{r,\,s}$ onto $A_{r,\,s}$ and if this end-point is $b$, starting from $b$ and taking successive elements in the opposite direction of the arrows gives a bijection of $A^\prime_{r,\,s}$ onto $A_{r,\,s}$.

We call either of these two sets a \textit{batch of arcs} and we speak of $\theta$ as defining a correspondence between the batch $\mathscr B'$ defined by $A_{r,\,s}^\prime$ and the batch $\mathscr B$ defined by $A_{r,\,s}$.

\subsection{Fixed Points.}\label{2.3}

Let $\sigma=\sigma_{[1,\,n],\,\mathscr M}$ be an anti-Toeplitz involution. Let $r$ be maximal such that 
$[r,\,s-1]\,:\,1\le r\le s-1\le[n/2]$
 is a connected component of $\mathscr M$ and set ${\mathscr M}^{max}_{r,\,s}=[r,\,s-1]\sqcup \kappa[r,\,s-1]$.

\begin {lemma}

\begin{enumerate}

\item[(i)] If $s-1<[n/2]$ then $\sigma$ has no fixed points except $(n+1)/2$ if $n$ is odd.  If $s-1=[n/2]$ it has one fixed point if $n$ is odd and two fixed points if $n$ is even. In both these last two cases all the fixed points lie in ${\mathscr M}^{max}_{r,\,s}$.

\item[(ii)]  Suppose $t \in [1,\,n]$ is not a fixed point of $\sigma$, then 
either $t$ or $\sigma(t)$ is equal to $(n+1)/2$ or $t$ and $\sigma(t)$ lie on opposite sides of $(n+1)/2$.

\end{enumerate}
\end {lemma}

\begin {proof} Assertion (i) (resp. (ii)) is an immediate consequence property of $T(i)$ (resp. $T(ii)$).
\end {proof}

\subsection{Modified Meanders.}\label{2.4}

If $J$ is a connected subset of $[1,\,n]$ of cardinality $>1$ we define $\sigma_{J,\,\mathscr M}$  on $J$ as in \ref {2.1} by replacing $[1,\,n]$ by $J$, $\mathscr M$ by its intersection with $J$,
and extending $\sigma_{J,\,\mathscr M}$ to $[1,\,n]$ by the identity on $[1,\,n]\setminus J$.  
It can occur that the intersection of $\mathscr M$ with $J$ is empty. In this case we set $\sigma_{J,\,\mathscr M}=\kappa_J$, the restriction of $\kappa$ to $J$.
If $J$ is reduced to one point of $[1,\,n]$, then $\sigma_{J,\,\mathscr M}$ is defined to be the identity on $J$.
Given a partition
$\mathscr J=\{J_i\}$ of $[1,\,n]$ we let $\sigma_{\mathscr J,\,\mathscr M}$ be the product of the $\sigma_{J_i,\,\mathscr M}$, which is again an involution.  Given two such partitions $\mathscr J^+,\,\mathscr J^-$ we write $\sigma_{\mathscr J^\pm,\,\mathscr M}$, simply as $\sigma^\pm$. We call these the modified (integer) involutions and the arcs they define, the modified integer arcs. 
Recall that every modified integer arc is a directed arc with its arrow pointing from left to right (resp. right to left) if it is an arc above (resp. below) $L$, as specified in \ref{1.1} or \ref{1.2}.
The arcs defined by $\kappa^\pm$ will be called unmodified integer arcs.
Of course a particular modified integer arc 
may remain unchanged under the above modification. 
Such an arc will be called a trivially modified integer arc and the remaining arcs non-trivially modified integer arcs.

 The meanders described in \ref {1.3} obtained thanks to these modified involutions $\sigma^\pm$  will be called the modified integer meanders. Again a particular modified integer meander may rest unchanged after modification.

\subsection{Fixed Points Revisited.}\label{2.5}

We explicit some easy (but important) conclusions from \ref {2.3} in the framework of \ref {2.4}.

Recall \ref {1.9} and let $J^+\cap J^-$ be a component of the double partition.  Recall our standing hypothesis \ref {1.8}. The condition $\pi^+\cup \pi^- =\pi$ implies that $J^+\neq J^-$ unless $\pi^+=\pi^-=\pi$ which we also excluded.  Let $c^+$ (resp. $c^-$) denote the centre of $J^+$ (resp. $J^-$).  Just in this section we let $\sigma^+$ (resp. $\sigma^-$) denote its restriction to $J^+$ (resp. $J^-$).

Assume that $c^+$ is integer.  Then $\sigma^+$ has just one fixed point.  It is $c^+$ if the arc defined by $\kappa^+$ joining the points $c^+ \pm 1$ is unmarked. Otherwise it may be $c^+ \pm s$, if the arcs defined by $\kappa^+$ joining the points $c^+ \pm t\,:\, t=1,\,2,\,\ldots,\,s$ are all marked.

Assume that $c^+$ is half-integer but not integer.  Then $\sigma^+$ has no fixed points unless the arc defined by $\kappa^+$ joining $c^+\pm \frac{1}{2}$ is marked. In this latter case $\sigma^+$ has exactly two fixed points.  Moreover these lie on opposite sides of $c^+$ and at a distance $\leq s+ \frac{1}{2}$ from $c^+$ given that the arcs defined by $\kappa^+$ joining the points $c^+ \pm (t+\frac{1}{2})\,:\, t=0,\,1,\,\ldots,\,s$ are all marked and the arc joining the points $c^+\pm(s+\frac{3}{2})$ is not marked.
 
Of course a similar result with $+$ replaced by $-$ is also valid.

\section{Construction of Markings and Some Consequences.}\label{3}

\subsection{Half-Integer Arcs.}\label{3.1}

Fix a double partition 
$\mathscr J^\pm=\{J_i^\pm\}$ of $[1,\,n]$.  Let $c_i^\pm$ be the centre of $J^\pm_i$ viewed as point on $L$ (these may include half-integer points). Let $\kappa^\pm$ denote the corresponding involutions viewed as arcs above and below $L$.

From this data we can recover the subsets $\pi^\pm$ of $\pi$ as in \ref {1.7}.

 Each connected component $\pi^\pm _i$ is equipped with a Dynkin diagram automorphism $\iota_i^\pm$ interchanging ends. Let $\iota^+$ (resp. $\iota^-$) denote the product of the commuting involutions $\iota^+_i$ (resp. $\iota^-_i$).  In terms of link patterns the $\iota^\pm$ may be represented as arcs above and below $L$ lying between the integer arcs and joining half-integer points of $L$, where specifically the point $i+\frac{1}{2}$ on $L$ represents $\alpha_i$. These arcs (as well as their extensions below) will be referred to as the half-integer arcs.

 The arc joining $i+\frac{1}{2}$ to $j+\frac{1}{2}$ for $i\leq j$ is viewed as representing the root $\alpha_i+\alpha_{i+1}+\ldots+\alpha_j$.

The $\iota_i^\pm$ are extended to involutions on an extended simple root system $\widetilde{\pi}$ as in \cite [3.1]{J1}. This may be described through our present diagrammatic presentation as follows.

Consider $\alpha \in \pi$.  Since $\pi^+ \cup \pi^- = \pi$ by hypothesis, $\alpha$ belongs to say $\pi^+$.  Assume that $\alpha \notin \pi^-$, so then $\iota^-$ is not defined on $\alpha$. Yet we may still follow $\alpha$ along $\ldots \iota^-\iota^+ \alpha$ until we reach $\alpha$ again or a new element $\alpha^\prime \in \pi$ at which the process cannot be further continued.

In the first case we set $\iota^-(\alpha)=\alpha$ and call the corresponding arc a fictitious one.  In the second case if it is $\iota^-$ which is not defined on $\alpha^\prime$, we set $\iota^-(\alpha)=\alpha^\prime$ and call the corresponding arc a fictitious one. Otherwise we adjoin \cite [3.1]{J1} a fictitious simple root $\widetilde{\alpha}$,  and set $\iota^-(\alpha)=\widetilde{\alpha},\,\iota^+(\widetilde{\alpha})=\alpha^\prime$ and call both of the corresponding arcs fictitious. (See \cite [11.2.1]{J1} for an example of the latter.)

In the above fashion we obtain meanders defined by the half-integer arcs, presented in \cite [3.2]{J1} as orbits of the group $<\iota^+,\,\iota^->$. We refer to them as the half-integer meanders. As in \ref {1.3} these may be subdivided into edges and loops.  
On a loop we mark any and exactly one arc \textit{except} a fictitious arc.  On an edge we exactly mark one non-fictitious arc (of the two possible) joining a fixed point of one of the involutions $\iota^\pm$ to itself. (This non-fictitious arc was called an end in \cite [3.3]{J1}.)

Now suppose this rule gives a mark on the half-integer arc above (resp. below) $L$ joining $i+\frac{1}{2}$ to $j+\frac{1}{2}$.  Of course we may assume that $i\leq j$.  Furthermore $i,\,j+1$ will both lie in the same connected component of $\mathscr J^+$ (resp. $\mathscr J^-$).   Then we translate this to a mark on the integer arc above (resp. below) $L$ joining $i,j+1$.  (The logic here is that half-integer arc and its neighbouring integer arc surrounding it represent the same root, namely $\alpha_i+\alpha_{i+1}+\ldots+\alpha_j$.)
Actually every half-integer arc we will consider will represent a root $\beta'_i=\varepsilon_i-\varepsilon_{\kappa(i)}$ of the Kostant cascade, that is will join the half-integer points $i+\frac{1}{2}$ and $j+\frac{1}{2}$ where $j=\kappa(i)-1=\iota(i)$. Thus we may also call an half-integer arc, an half-integer arc of the Kostant cascade (here $\kappa$ and $\iota$ designate respectively $\kappa^+$ and $\iota^+$ if the half-integer arc is above $L$, or $\kappa^-$ and $\iota^-$ if it is below $L$).

This construction gives the set of markings on the integer arcs.

\textbf{Remarks}.   1) It is of importance to stress that it is not possible to directly mark the unmodified integer arcs.  Moreover there are two good reasons for this.  First, at a turning point
of a half-integer meander, the integer arcs assigned as above to the half-integer arcs, lie respectively above and below the latter and so are not joined on $L$. Consequently they need not belong to the same integer meander.
Second, the integer and half-integer meanders measure different things.  More precisely the number of orbits of $<\kappa^+\kappa^->$ (resp. of $<\iota^+\iota^->$) equals $\ell(\mathfrak q) +1$ (resp. $\ell(\mathfrak q_{\Lambda})$) - see \cite [A4,\,A5] {J0.5}.\par
2) The rule of modification of integer arcs given in \ref{2.1} translates as follows.
All the unmodified integer arcs encircling the marked half-integer arcs (below and above each of them) have to be non-trivially modified and only these ones.

\subsection{Edge Decomposition.}\label{3.2}

From these markings we proceed as in Section \ref{2} to construct $\sigma^\pm$, that is to say the
modified integer arcs from which the modified integer 
meanders are obtained. We regard the modified integer non-trivial arcs, above and below $L$, as giving subsets $S^+ \subset \Delta^+,\, S^- \subset \Delta^-$ of $R$. Up to the sign convention defined by the arrows on arcs (\ref{1.2}) one has  $S^\pm=\{\varepsilon_i-\varepsilon_{\sigma^\pm(i)}\,:\, i \in [1,\,n],\,i\neq\sigma^\pm(i)\}$.

We need to know that the elements of $S:=S^+\cup S^-$ are linearly independent. For this we generalize \cite [Prop. 5.7]{J1} to the present situation.  We remark that this result was one of two crucial observations in \cite {J1} to construct an adapted pair.  The second was \cite [Thm. 8.6]{J1}.  Unfortunately the latter does not in general go over to arbitrary anti-Toeplitz involutions.  Consequently $\eta$ as defined in \ref {1.8} $(*)$ need not be regular.

Following \cite [2.4]{J1} we let 
$\mathfrak h_{\Lambda}$
 denote the common kernel of the weights of the semi-invariant elements of $S(\mathfrak q)$ in the Cartan subalgebra $\mathfrak h$ of the biparabolic algebra $\mathfrak q$.

\begin {prop} $S|_{\mathfrak h_{\Lambda}}$ is a basis of $\mathfrak h_{\Lambda}^*$.  In particular the elements of $S$ are linearly independent.
\end {prop}

\begin {proof} This results from \cite [Lemmas 5.2, 5.5 and 5.6]{J1}. It suffices to briefly discuss why the first of two of these still hold giving only more extensive details for the third.

The basis of \cite [Lemma 5.2]{J1} is \cite [5.2, Eq.(4)]{J1}.  This only applies to connected components of $\mathscr M$ for which $s-1 <[n/2]$ in the notation of \ref {2.1} (where we only mark arcs above $L$ whose end-points lie in $[1,\,n]$).  In this case the meander is a loop, moreover as we noted in \ref {2.2.1} when we take the arcs in $A_{r,\,s}$ (resp. $A^\prime_{r,\,s}$) to correspond to negative (resp. positive) roots, then the arrows are aligned along a direction of the loop. Under this condition a loop sum is obviously zero. Thus we obtain the sum rule
$$\sum_{\beta \in A_{r,\,s}}\beta =-\sum_{\gamma \in A^\prime_{r,\,s}}\gamma.$$

From this \cite [5.2, Eq.(4)]{J1} immediately follows.

In our construction of an anti-Toeplitz involution we ensured that the number of arcs it defined was the same as in the special case described in \cite [Sect. 5]{J1} using the anti-Coxeter and anti-Jordan blocks.  (Recall the penultimate paragraph of \ref {2.1}.) Thus \cite [Lemma 5.5]{J1} still holds.

Let us now show that \cite [Lemma 5.6(i)]{J1} also still holds.  This reduces to the case when either $\mathscr J^+$ or $\mathscr J^-$ equals $\{[1,\,n]\}$ (consider only the first case).  Then in the present notation
if we identify the set $\mathscr M\subset[1,\,[n/2]]$ with a subset of marked simple roots it defines (see beginning of \ref{2.1})
 we must show that
$$\mathbb Z \mathscr M \subset \mathbb Z \mathscr K + \mathbb Z S \eqno {(*)}.$$

Recall the notations of \ref {2.1} and \ref{2.2.1} and fix a connected component $[r,\,s-1]$ of $\mathscr M$.  By definition $\mathbb Z\mathscr K$ (resp. $\mathbb Z S$) contains the set $\mathbb Z A_{r,\,s}$ (resp $\mathbb Z A^\prime_{r,\,s}$) since $A_{r,\s}$ (resp. $A'_{r,\,s}$) is a set of arcs defined by $\kappa_{\mathscr M_{r,\,s}}$ (resp. $\sigma_{\mathscr M_{r,\,s}}$) viewed as roots.  Then $(*)$ is reduced to showing that
$$\alpha_t \in \mathbb Z A_{r,\,s} + \mathbb Z A^\prime_{r,\,s}, \forall t \in [r,\,s-1] \eqno {(**)}.$$

By a relabeling of the simple roots we can assume that $r=1$.

$T(i)$ means that the right hand side of $(**)$ is a free $\mathbb Z$ module of rank $2s-1$ unless $s=[n/2]+1$ in which case it has rank $n-1$.  In the latter case this means that the right hand side is just $\mathbb Z \pi$, since the containing Lie algebra is $\s\l_n(\bf k)$. In this case $(**)$ holds trivially. In the former case, $T(ii)$ implies that the right hand side of $(**)$ is contained in
$$\bigoplus_{t \in [1,\,s-1]}(\mathbb Z \alpha_t \oplus \mathbb Z \iota^+(\alpha_t)) \oplus \mathbb Z (\alpha_s+\ldots +\iota^+(\alpha_s)),$$
which is also a free $\mathbb Z$ module of rank $2s-1$.  Again since the containing Lie algebra is $\s\l_n(\bf k)$, so that all the coefficients of a simple root in a positive root are $\leq 1$, the equality of ranks implies that this is exactly the right hand side of $(**)$.  From this $(**)$ follows.

The proof of \cite [Prop. 5.6(ii)]{J1} is unchanged, being simply that $\mathscr M$ contains a representative of every $<\iota^+,\,\iota^->$ orbit (by \ref{3.1}).  Since the last part of \cite [Prop. 5.6]{J1} follows from its first two parts, it also extends to the present situation.

To complete the proof of the proposition we must show that $\mathfrak h^\perp_{\Lambda}+ \sum_{\alpha \in S}\bf k\alpha = \mathfrak h^*$ and so is a direct sum by \cite [Lemma 5.5]{J1}.  This follows exactly as in \cite [5.7]{J1} using \cite [Lemmas 5.2, 5.6(iii)]{J1}.

\end {proof}

\subsection{Sign Changes.}\label{3.3}

Since loop sums are zero, the above proposition immediately gives the following rule

$R(i)$.  None of the modified integer meanders can be loops.

Consequently $S$ forms a disjoint union of non-trivial edges.

Recall that the second element $\eta$ of an adapted pair constructed in \cite{J1} is given as in \ref{1.8}$(*)$. For every non-trivial edge $E$ (that is for every non-trivial modified integer meander), denote by $S_E$ the subset of $S$ formed by the roots corresponding to the non-trivial modified integer arcs whose end-points belong to the set of points of $E$. Set $\eta_E=\sum_{\alpha\in S_E}x_{\alpha}$. 
Let ${\bf E}(S)$ denote the set of non-trivial edges defined by $S$. 
Then we have
$$\eta=\sum_{E\in {\bf E}(S)}\eta_E.$$

To simplify notation we shall just assume that $S$ corresponds to a single non-trivial edge $E$. This assumption will be dropped in Section \ref{8} when edges are joined.  Recover the notations of \ref {1.3} - \ref {1.5} (in particular the notation of the $\beta_i$ given in \ref{1.5} for a given edge $E$ of cardinality $e$).  It is clear that $|S|=e-1$.  Again it is clear that $(\beta_i,\,\beta_{i+1})=-1$, for all $i,\,i+1 \in I_{1,\,e}$ whilst the remaining scalar products of distinct elements vanish. Thus $\Pi:=\{\beta_i\}_{i \in [1,\,e-1]}$ is a simple root system of type $A_{|S|}$.

Again it is clear that there exist $\epsilon_i \in \{\pm 1\}\,:\,i \in [1,\,e-1]$ such that $S=S^+\cup S^-=\{\epsilon_i\beta_i : i \in [1,\,e-1]\}$. (Recall that $S\subset R$ with $S^+\subset\Delta^+$ and $S^-\subset\Delta^-$, \ref{3.2}).

Following our previous work \cite [Lemma 2.5] {FJ1}, we describe how the $\epsilon_i$  behave, through the following

\begin {lemma} For all $i \in [1,\,e-2]$ one has $\epsilon_i\epsilon_{i+1}=-1$, if $\varphi(i+1)$ is an internal turning point of $E$, and equals $1$ otherwise.
\end {lemma}

\begin {proof}  Straightforward.
\end {proof}

Recall the definition of sources and sinks given in \ref{1.4} and that the roots of $R$ are represented above (resp. below) $L$ with an arrowhead pointing from left to right (resp. right to left), see \ref{1.1}, \ref{1.8} and \ref{2.4}. Observe also that $\epsilon_i\beta_i\in R$ for all  $i\in I_{1,\,e}$. Then we get the following :

\begin{cor}
Let $\varphi(r)$ be a turning point of $E$ not equal to the finishing point of $E$. If $\varphi(r)$ is a source (resp. a sink) then $\epsilon_r=1$ (resp. $\epsilon_r=-1$).

\end{cor}

\begin {proof}  Straightforward.
\end {proof}

\subsection{An Exclusion of Markings.}\label{3.6}

Recall \ref {2.1} that marking a subset of $\mathscr K$ is equivalent to marking a subset of $\overline{\pi}:=\pi/\iota$.

More generally in the notation of \ref {3.1}, we set $\overline{\pi}^\pm:=\pi^\pm/\iota^\pm$.  Since only one (half-integer) arc of each $<\iota^+,\,\iota^->$ orbit is marked, we have the following rule.

$R(ii)$. Suppose $\alpha \in \pi^+ \cap \pi^-$.  If $\alpha$ is marked as an element  of $\overline{\pi}^+$, it is not marked as an element of $\overline{\pi}^-$ and vice-versa.

\subsection{A Further Exclusion of Markings.}\label{3.7}

 Retain the notation of \ref {3.1} and consider a double partition $\mathscr J^\pm=\{J_i^\pm\}$ of $[1,\,n]$. When $|J^\pm_i|=1$ we call the corresponding point on $L$ a singleton
and $J_i^\pm$ a trivial component of $\mathscr J^\pm$.  The arcs above (resp. below) $L$ join the singletons to themselves and 
were called (\ref{1.1})
 trivial (integer) arcs. By \ref{2.4} trivial integer arcs of trivial components of $\mathscr J^\pm$ are not marked 
and are trivially modified. (By contrast a half-integer arc joining 
a half-integer point to itself represents an element which may be marked, but not simultaneously above and below $L$ - this condition is $R(ii)$.)   An integer arc (modified or not) joining distinct points $i <j$ of $L$ represents an element of $\pm R$.  Recall (\ref {1.1}) that its support is defined to be the support of the root $\varepsilon_i-\varepsilon_j$ it represents.  The supports of the arcs above (resp. below) $L$ lie in the connected subsets of $\pi^+$ (resp. $\pi^-$).

 Inclusion of supports gives an order relation on the set of non-trivial integer unmodified (or modified) arcs above (resp. below) $L$.  The set of all unmodified arcs whose support lies in a fixed $\pi^+_i$ (resp. $\pi^-_i$) are strictly linearly ordered.   As noted in the end of \ref {2.1} this is not true of the set of modified arcs because they may cross.

 We now give a third rule on exclusion of markings complementing $R(i)$, $R(ii)$.  This additional rule obtains from the necessary existence of a gap between certain markings of the integer arcs (defined by the $\kappa_{J_i^\pm}$).  It is important because it limits the crossings of the modified integer arcs (defined by the $\sigma^\pm$).  Details are given for arcs above $L$.  Arcs below $L$  behave similarly.
 In this ``can" does \textit{not} mean ``must".

Assume that $\pi^-_i$ is not empty, that is to say that $J_i^-$ is not a trivial component of $\mathscr J^-$.

The proofs of $R(iii)_e^+$ and $R(iii)_o^+$ below are illustrated by figures (see Figures \ref{Figure 1}, \ref{Figure 2} and \ref{Figure 3}) in which $c^+$ denotes the centre of a component $J^+\in\mathscr J^+$ whose intersection with the set of integer points of $J_i^-$ encircling its centre $c_i^-$ is non-empty. (To  simplify the figures, we shall denote $J_i^-$ by $J^-$ and $c_i^-$ by $c$. The half-integer arcs will be represented by dashed lines and the unmodified integer arcs by solid lines).

\subsubsection{}\label{3.7.1}

First assume that the centre $c^-_i$ of $J_i^-$ is half-integer (and not integer), equivalently that $|J^-_i|$ is even.  This means that $c^-_i$ is a fixed point of $\iota^-$.  In particular the half-integer meander meeting $c_i^-$ is an edge and then only an end of this meander can be marked, by \ref {3.1} and this implies that only a half-integer point of this meander fixed by $\iota^-$ or by $\iota^+$ can be marked.


$R(iii)_e^+$.  If the unmodified integer arcs above $L$ meeting the points $c^-_i\pm \frac{1}{2}$ are $(a)$ distinct and $(b)$ strictly linearly ordered, then only the smaller one can be marked.

\begin{figure}[!h]
\centering
\input{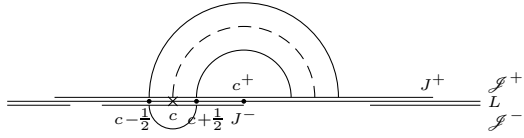}
\caption{\scriptsize{\it 
The case $c+\frac{1}{2}=c_i^-+\frac{1}{2}\le c^+$.}
The half-integer arc (dashed line) is not marked 
and so the outermost integer arc cannot be marked.}
\label{Figure 1}
\end{figure}

\begin{proof}
Indeed $c^-_i$ cannot also be a fixed point of $\iota^+$ by $(a)$.  Thus the half-integer arc above $L$ meeting $c^-_i$ is not an end so cannot be marked.  Then by $(b)$, of the integer arcs above $L$ only the smaller one can be marked (see Figure \ref{Figure 1}).
\end{proof}


\subsubsection{}\label{3.7.2}

Assume that $c^-_i$ is integer, equivalently that $|J^-_i|$ is odd and $\geq 3$.

$R(iii)_o^+$.  If the unmodified integer arcs above $L$ meeting the points $c^-_i, \,c^-_i\pm 1$ are distinct and strictly linearly ordered, then at most one of the larger two can be marked. If those meeting $c_i^--1,\,c_i^-$ coincide, then if the one meeting $c_i^-+1$ is strictly larger, it cannot be marked.

\begin{figure}[!h]
\centering
\input{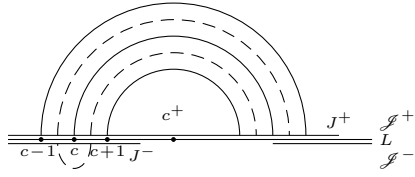}
\caption{\scriptsize{\it 
The case $c+1=c_i^-+1\le c^+$. }
One of the dashed arcs above $L$  is not marked. }
\label{Figure 2}
\end{figure}

\begin{figure}[!h]
\centering
\input{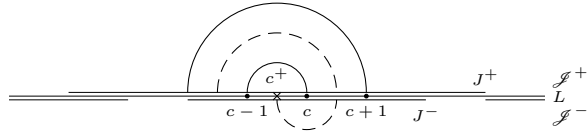}
\caption{\scriptsize{\it
 The case $c^+=c_i^--\frac{1}{2}=c-\frac{1}{2}$.}
 The dashed arcs are not marked.}
 \label{Figure 3}
 \end{figure}

\begin{proof}
Indeed the hypothesis $|J^-_i|\geq 3$ implies that $c_i^-\pm \frac {1}{2}$ are joined by a half-integer arc below $L$.  Thus the half-integer arcs above $L$ meeting these two points cannot both be marked if they are distinct (see Figure \ref{Figure 2}).  In the second case $c_i^--\frac{1}{2}$ is a fixed point of $\iota^+$ and so the half-integer meander meeting it is an edge. Thus the half-integer arc above $L$ meeting $c_i^-+\frac{1}{2}$ is not an end and so cannot be marked (see Figure \ref{Figure 3}).
\end{proof}

\subsubsection{}\label{3.7.3}
The corresponding assertions for the integer arcs below $L$ are designated by interchanging  the superscripts $+$ and $-$.

\subsection{More on Fixed Points.}\label{3.8.1}

Retain the notation and conventions of \ref {2.5}.   Write $J^\pm=[\ell^\pm,\,r^\pm]$ and $\mathscr I=J^+\cap J^-$.

\begin {lemma}   Assume that $c^+<c^-$.

\begin{enumerate}

\item[(i)]  Suppose $c^-\leq r^+$. Then there is no fixed point of $\sigma^+$ in 
$\mathscr I\cap [c^-,\,r^+]$, which is not a fixed point of $\sigma^-$.

\item[(ii)]  Suppose $\ell^-\leq c^+$. Then there is no fixed point of $\sigma^-$ in 
$\mathscr I\cap[\ell^-,\,c^+]$, which is not a fixed point of $\sigma^+$.
\end{enumerate}

\end {lemma}

\begin {proof}  Recall the rules given in \ref {3.1} describing the rules to define markings and the rules given in \ref {2.5} describing fixed points.

$(i)$ Let $a$ be a fixed point of $\sigma^+$ in 
$\mathscr I\cap [c^-,\,r^+]$, which is not a fixed point of $\sigma^-$.

Assume $c^-$ to be half-integer (and not integer).

Then $a \geq c^-+ \frac{1}{2}$ and by $R(iii)^+_e$ taking $c^-=c^-_i$ the (unmodified) integer arc above $L$ meeting $c^-+ \frac{1}{2}$  cannot be marked, which gives a contradiction with \ref {2.5}.

Assume $c^-$ to be integer.

Suppose that $a \geq c^- +1$. Then 
by $R(iii)^+_o$ taking $c^-=c^-_i$ the (unmodified) integer arcs above $L$ meeting $c^-+1,\,c^-$  cannot both be marked, which gives a contradiction with \ref {2.5}.

Finally suppose that $a=c^-$.  Since by hypothesis $a=c^-$ is not a fixed point of $\sigma^-$ it follows from \ref {2.5} that the (unmodified) integer arc below $L$ joining $c^- -1,\,c^- +1$ must have been marked, which in turn means that the half-integer arc below $L$ joining $c^- -\frac{1}{2},\,c^- +\frac{1}{2}$ must have been marked.  However in this case the half-integer arc above $L$ meeting $c^- -\frac{1}{2}$ cannot have been marked (because it belongs to the same half-integer meander) and in turn the (unmodified) integer arc above $L$ meeting $c^-$ cannot have been marked.  Since $a=c^-$, this gives a contradiction with \ref {2.5}.

The conclusion $(i)$ obtains.

Consider $(ii)$.

Let $a$ be a fixed point of $\sigma^-$ in 
$\mathscr I\cap [\ell^-,\,c^+]$ which is not a fixed point of $\sigma^+$.

Assume $c^+$ to be half-integer (and not integer).  Then the half-integer arc below $L$ meeting $c^+$ cannot be marked because it lies on an edge with end-point $c^+$. Thus the integer arc below $L$ meeting $c^+-\frac{1}{2}$ cannot be marked (this actually results from $R(iii)^-_e$ with $c_i^+=c^+$), hence 
 $a\le c^+-\frac{1}{2}$ cannot be a fixed point of $\sigma^-$ by \ref {2.5}. 

Assume $c^+$ to be integer.



If $a\le c^+-1$ then by $R(iii)^-_o$ taking $c_i^+=c^+$,  the (unmodified) integer arcs below $L$ meeting $c^+-1,\,c^+$  cannot both be marked, which gives a contradiction with \ref {2.5}.\par
If $a=c^+$  then, since $a$ is not a fixed point of $\sigma^+$, it follows from \ref {2.5} that the (unmodified) integer arc above $L$ joining $c^+ -1,\,c^+ +1$ must have been marked, which in turn means that the half-integer arc above $L$ joining $c^+-\frac{1}{2},\,c^+ +\frac{1}{2}$ must have been marked.  However in this case the half-integer arc below $L$ meeting $c^+ +\frac{1}{2}$ cannot have been marked (because it belongs to the same half-integer meander) and in turn the (unmodified) integer arc below $L$ meeting $c^+$ cannot have been marked, which gives a contradiction with \ref {2.5}. Hence $(ii)$.

\end {proof}

\textbf{Remark}. If $c^-\le r^+$, then the above proof shows that the only fixed point of $\sigma^+$ in 
$\mathscr I\cap[c^-,\,r^+]$ may be $c^-$ and in this case $c^-$ is integer and fixed under both $\sigma^+$ and $\sigma^-$.\par
Similarly
if $\ell^-\leq c^+$, then the only fixed point of $\sigma^-$ in $\mathscr I\cap[\ell^-,\,c^+]$ may be $c^+$ and in this case $c^+$ is integer and fixed under both $\sigma^+$ and $\sigma^-$.

\subsection{More on turning points lying in double components.}\label{3.4.5}


Consider a component $\mathscr I=J^+\cap J^-$ of the double partition (see \ref{1.9}) with $J^\pm\in\mathscr J^\pm$. Recall the notational convention introduced in the last paragraph of \ref {2.1} and set $\sigma^+ =\sigma_{J^+}, \sigma^-=\sigma_{J^-}$.  Let $c^+,\,c^-$ be the centres of $J^+$ and $J^-$ respectively. Recall also the definitions in \ref{1.1} and in \ref{1.4}.

Our goal is to describe the turning points of an edge $E$ in $\mathscr I$ and to determine whether they are sources or sinks.

\subsubsection{Observations.}\label{3.4.5.0}

Consider a directed arc $(r,\,s)$ defined by $\sigma^+$, whose end-points lie in $J^+$.
Then if this arc is non-trivial, such an arc points from left to right.

If $c^+$ is not an end-point of $(r,\,s)$ and if $(r,\,s)$ is non-trivial, then moreover by \ref{2.3}(ii) $r$ and $s$ lie on opposite sides of $c^+$.

If $c^+$ is an end-point of $(r,\,s)$, then $c^+$ is integer and is either a fixed point of $\sigma^+$ (and $r=s=c^+$ and the arc is trivial) or the half-integer arc above $L$ joining the points $c^+\pm\frac{1}{2}$ is marked by \ref{2.5}.

Thus every element of $J^+$ not fixed under $\sigma^+$ and strictly to the left (resp. strictly to the right) of $c^+$ must be the starting (resp. the finishing) point of a non-trivial directed arc defined by $\sigma^+$.

Similar assertions hold if $\pm$ are interchanged.

\subsubsection{}\label{3.4.5.1}
From the above observations we obtain the

\begin{lemma}
Assume $c^+<c^-$.

\begin{enumerate}
\item[(i)] An internal turning point of an edge $E$ in $\mathscr I$ must lie in $\mathscr I\cap[c^+,\,c^-]$ and is a sink.

\item[(ii)] Every element of $[c^+,\,c^-]\cap\mathscr I$ must be a turning point of some edge (possibly trivial) with just the following possible exceptions : $c^+$ (resp. $c^-$) when it 
is integer and the half-integer arc above (resp. below) $L$ joining the points $c^+\pm \frac {1}{2}$ (resp. $c^-\pm\frac{1}{2}$) is marked.

\item[(iii)]
An end-point of a non-trivial edge $E$ lying in $\mathscr I$ is a sink.
\end{enumerate}
\end{lemma}




\begin{proof}
Lemma \ref {3.8.1} implies that an end-point of a non-trivial edge $E$ lying in $\mathscr I$ fixed by $\sigma^+$ (resp. $\sigma^-$) lies strictly to the left of $c^-$ (resp. strictly to the right of $c^+$). Then $(iii)$ follows also from the above observations.  \end{proof}





\subsubsection{}\label{3.4.5.2}

In exactly the same manner we obtain the

\begin{lemma}
Assume $c^+>c^-$. 

\begin{enumerate}
\item[(iv)] An internal turning point of an edge $E$ in $\mathscr I$ must lie in $\mathscr I\cap [c^-,\,c^+]$ and is a source.

\item[(v)] Every element of $[c^-,\,c^+]\cap\mathscr I$ must be a turning point of some edge (possibly trivial) with just the following possible exceptions : $c^+$ (resp. $c^-$) when it is integer and the half-integer arc above (resp. below) $L$ joining the points $c^+\pm \frac {1}{2}$ (resp. $c^-\pm\frac{1}{2}$) is marked.

\item[(vi)] An end-point of a non-trivial edge $E$ lying in $\mathscr I$ is a source.
\end{enumerate}
\end{lemma}

\subsubsection{}\label{3.4.5.3}

\begin{lemma}

Assume $c^+=c^-:=c$.

\begin{enumerate}
\item[(vii)] There are no internal turning points of an edge in $\mathscr I$.

\item[(viii)] If $c$ is integer, then $c$ is an end-point of a non-trivial edge.

\item[(ix)] The combined number of fixed points of $\sigma^+$ and $\sigma^-$ in $\mathscr I$ is at most two and if this number is two, then one is a sink and the second is a source of some edge (not necessarily the same one).
\end{enumerate}
\end{lemma}

\begin{proof}

Through \ref{3.4.5.0} the only possible internal turning point of an edge in $\mathscr I$ is $c$ when $c$ is integer.
Thus suppose $c$ is integer. Then either the arc defined by $\kappa^+$ joining $c\pm 1$ is marked, or the arc defined by $\kappa^-$ joining $c\pm 1$ is marked, but not both (by $R(ii)$).\par
In the first case the arc defined by $\kappa^-$ joining $c\pm 1$ is unmarked and $c$ is the unique fixed point of $\sigma^-$ in $J^-$. At the same time $\sigma^+$ has a unique fixed point in $J^+$ distinct from $c$. Moreover by $T(ii)$, if $c$ is the starting (resp. finishing) point of an arc defined by $\sigma^+$, then the fixed point $f$ of $\sigma^+$ lies to the left (resp. right) of $c$. Then $c$ is a source (resp. sink) whereas if $f$ is also in $J^-$, then it is a sink (resp. source).\par
The second case is similar.\par
Suppose that $c$ is half-integer and not integer. Then by $R(ii)$ and \ref{2.5}, $\sigma^+$ (resp. $\sigma^-$) has two fixed points in $J^+$ (resp. $J^-$) on opposite sides of $c$ and $\sigma^-$ (resp. $\sigma^+$) has no fixed points in $J^-$ (resp. $J^+$). Then by \ref{3.4.5.0} one is a source and the other is a sink.

\end{proof}

\subsubsection{Terminology.}\label{3.4.5.4}

 1) When $c^+ =c^-$, we say that the double component $\mathscr I=J^+\cap J^-$ is equicentral.
\par\noindent
2) We also say that the double components $\mathscr I_1:=J^+_1\cap J_1^-,\,\mathscr I_2:=J^+_2\cap J_2^-$ are linked  if either $J_1^+=J_2^+$ or $J^-_1=J^-_2$.  In the parabolic case all double components are linked. If $\mathscr I_1$ and $\mathscr I_2$ are not linked, and if $f_i\in \mathscr I_i$, $i\in\{1,\,2\}$, then trivially $\varepsilon_{f_1}-\varepsilon_{f_2}\in K\cap -K$.

 \subsubsection{Remark.}\label{3.4.5.5}

The case of interest in the sequel is when $c^+=c^-$ and the combined number of fixed points in $\mathscr I$ is exactly two.
(See also Proposition \ref {3.8.4}(i) for more details).\par
Since equicentral components may contain one fixed point which is a sink and  another which is a source, they lead to technically more complex situations (cf. \ref{9.6} and \ref{9.10}).\par

 \subsection{More on equicentral components.}\label{3.8.4}

We will give in this subsection some more properties on equicentral components which we will use later.

 Retain the notation and hypotheses of \ref {2.5}  
and in addition suppose that the component $\mathscr I=J^+\cap J^-$ of the double partition is equicentral and set $c=c^+=c^-$.  Under our standing hypotheses (\ref {1.8}) on $\pi^\pm$ we can assume without loss of generality that $J^+$ overlaps $J^-$ by the same non-zero amount on both sides and that $J^-$ is not reduced to a point. Thus we can write $J^\pm=[\ell^\pm,\,r^\pm]$ with $t:=\ell^--\ell^+=r^+-r^- \in \mathbb N^+$. One has $r^\pm=2c-\ell^\pm$.  Then $\mathscr I = [\ell^-,\,r^-]$ viewed as a subset of $L$ and called the internal region. Let $\mathscr E$ denote its complement in $L$, called the external region.  
We will show that a non-trivial edge of the (modified) integer meander (defined by $S$) cannot be confined within $\mathscr I$. 
 Of course confinement of an edge in $\mathscr I$ can only occur if both its end-points lie in $\mathscr I$. In the following proposition we shall prove a little more.

 A half-integer meander in the internal region consists of an arc above $L$ completed by an arc below $L$, exactly one of which must be marked 
by \ref{3.1} and $R(ii)$.  On the other hand the half-integer arc above $L$ joining the points $\ell^--\frac{1}{2},r^-+\frac{1}{2}$, which according to the algorithm of  
 \ref{3.1}
  has to be completed by a fictitious arc below $L$, must also be marked 
(\ref{3.1}).

 We introduce a numerical labeling to describe the marked half-integer arcs meeting the half-integer points of $L$ lying between $\ell^+$ and $r^+$.

 Let $\ell_1:  \ell^+ <  \ell_1 \leq \ell^-$ be minimal such that half-integer arcs above $L$ meeting the points $t-\frac{1}{2}: \ell_1 \leq t \leq \ell^-$ are all marked.

 Then let $\ell_1<\ell_2<,\ldots, <\ell_{s+1} = [c+\frac{3}{2}] $ be chosen so that the half-integer arcs meeting the points $t-\frac{1}{2}$ above (resp. below) $L$ are marked for all $t \in [\ell_{2i-1},\,\ell_{2i}-1]$ (resp. $t \in [\ell_{2i},\,\ell_{2i+1}-1]$) : $i=1,\,2,\,\ldots, [(s+1)/2]$ (resp. $[s/2]$).  Then the (unmodified) integer arcs  meeting the points $t-1$ above (resp. below) $L$ are marked for all $t \in [\ell_{2i-1},\,\ell_{2i}-1]$ (resp. $t \in [\ell_{2i},\,\ell_{2i+1}-1]$) : $i=1,\,2,\,\ldots, [(s+1)/2]$ (resp. $[s/2]$).

 Set  $r_t=2c-\ell_t$, and
 $\mathscr F_t =[\ell_t-1,\,r_t+1]$, $\mathscr F^\prime_t =[\ell_t,\,r_t]$, for all $t =1,\,2,\ldots,\,s$. Notice that $\mathscr F^\prime_t$ just omits the extremal points of $\mathscr F_t$
and that, if $s\ge 2$, and $t\in\mathbb Z$ satisfies $2\le t+1\le s$, then for all
 $u\in\mathscr F_t\setminus\mathscr F_{t+1}$, the integer arc  meeting $u$, above $L$  if $t$ is odd (resp. below $L$ if $t$ is even) is marked and then for all $u\in\mathscr F_t\setminus \mathscr F^\prime_{t+1}$, the integer arc meeting $u$, above $L$ if $t$ is odd (resp. below $L$ if $t$ is even)  is non-trivially modified. 
 
Then if $s\ge 2$, for all $2\le t\le s$, both integer arcs above and below $L$ meeting the point $\ell_t-1$ are non-trivially modified. 

 Set $\mathscr F= \mathscr F_s$. 
Observe that $\mathscr F\subset\mathscr F_t$ for all $1\le t\le s$.
The non-trivial integer arcs meeting the points of $\mathscr F$ are either all marked above $L$ or all marked below $L$ (depending on the parity of $s$ : if $s$ is odd (resp. even) they are all marked above $L$ (resp. below $L$)) and then these arcs are all non-trivially modified above (resp. below) $L$ if $s$ is odd (resp. even).

 Let $F$ be the union of the fixed points of $\sigma^+$ and of the fixed points of $\sigma^-$ lying in $\mathscr I$.  For all $f \in F$ let $M_f$ be the (modified) integer meander with end-point $f$.  (Some of these meanders may coincide, that is to say $M_f,\,M_{f^\prime} : f,\,f^\prime \in F$, might be joined in $\mathscr E$.)

 \begin {prop}
 \begin{enumerate}
 \item[(i)]  $F$ does not contain a  point
fixed under both $\sigma^+$ and $\sigma^-$.  Moreover $F \subset \mathscr F$ and $|F| \leq 2$, with equality if $\mathscr F \subset \mathscr I$.

\item[ (ii) ] For any $f \in F$, the points of $M_f$ do not lie entirely in $\mathscr I$.

Take $f\in F$. Starting from $f$ and following the meander $M_f$, we obtain a sequence of points $f_1=f,\,f_2,\,f_3,\,\ldots,\,f_k$ with $f_k$ the first point of $M_f$ not in $\mathscr I$. Set $\overline{M}_f=\{f_i\}_{i=1}^{k-1}=M_f\cap \mathscr I$.

\item[ (iii)] If $F\neq\emptyset$, then $\mathscr I=\cup_{f\in F}\overline{M}_f$.

 \item[(iv)] The first point at which $M_f$ leaves $\mathscr I$ is a turning point, that is in the above notation $f_k$ is a turning point.
 \end{enumerate}
 \end {prop}

\begin{figure}[!h]
\centering
\input{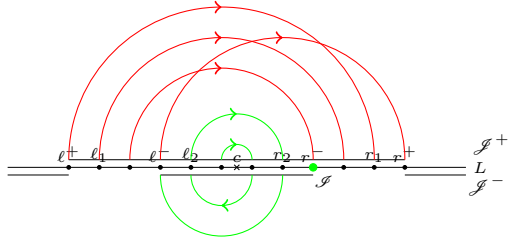}
\caption{\scriptsize{\it 
The case $s=2$ with $c$ half-integer (and not integer).} 
The  two fixed points under the modified involution
 $\sigma^-$ in $J^-$ are $c-\frac{1}{2}$ and $c+\frac{5}{2}$.
In green are represented the arcs of $N_2$ and one observes that $N_2$
 is a union of the trivial edge joining $r_2+1=c+\frac{5}{2}$ to itself
and another edge joining $\ell_2-1$ to $c-\frac{1}{2}$. 
In red are represented the non-trivial arcs of $N_1\setminus N_2$.
 We have $\mathscr I\subset M_{c-\frac{1}{2}}\cup M_{c+\frac{5}{2}}$.}
 \label{Figure 4}
 \end{figure}

\begin{figure}[!h]
\centering
\input{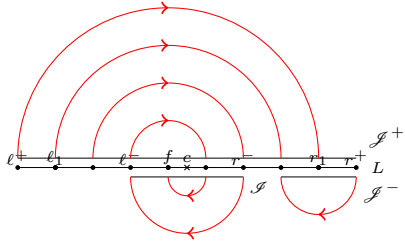}
\caption{\scriptsize
{\it The case $s=1$, $\mid F\mid =1$ with $c$ half-integer (and not integer).}
 The fixed points under $\sigma^+$ in $J^+$ are 
 $c-\frac{1}{2}$ and $c+\frac{9}{2}$ (only $c-\frac{1}{2}$ belongs to $\mathscr I$).
  Observe that $\mathscr I=\overline{M}_{c-\frac{1}{2}}$.}
  \label{Figure 5}
\end{figure}

 \begin {proof}

 $(i)$  Observe that $\mathscr F \subset \mathscr I$ if and only if $s >1$ (indeed if $s>1$, $\ell_s\ge\ell_2\ge\ell^-+1$ since the half-integer arc above $L$ meeting $\ell^--\frac{1}{2}$ is marked).

 Suppose $s >1$ and that $s$ is odd (resp. even).  Then the marked (unmodified) integer arcs 
meeting the points of $\mathscr F$ lie 
above (resp. below) $L$. 
 Count just the fixed points of $\sigma^{\pm}$ lying in $\mathscr I$ and apply \ref{2.5}. Suppose $c$ is integer. Then $\sigma^+$ (resp. $\sigma^-$) has just one fixed point and this lies in $\mathscr F\setminus\{c\}$ whilst $c$ is the only fixed point of $\sigma^-$ (resp. $\sigma^+$). If $c$ is half-integer and not integer, then $\sigma^+$ (resp. $\sigma^-$) has two fixed points, whilst $\sigma^-$ (resp. $\sigma^+$) has no fixed points. Hence $(i)$ in this case.


 Suppose $s=1$.  In this case similar considerations apply but the fixed points in $\mathscr F$ under $\sigma^+$ or $\sigma^-$ need not all lie in $\mathscr I$ and so only an inequality may hold for the cardinality of $F$.

 We shall prove $(ii)$ and $(iii)$ by an appropriate reverse induction. 
 
 Consider \textit{just} the modified
 integer arcs joining
 points of $\mathscr F_t$ (for $1\le t\le s$) and let $N_t$ denote the union of the meanders defined by these arcs (above and below $L$).
 
 It is convenient to first assume that $s>1$. 


 We prove by decreasing induction on $t$ that, for $t\ge 2$, $N_t$ is a union of two edges (possibly trivial) passing through all points of $\mathscr F_t$, one joining (not necessarily by a single arc) the first fixed point to one extremal point (say $\ell_t-1$) of $\mathscr F_t$ and the other joining the second fixed point 
 to the second extremal point  $r_t+1$ of $\mathscr F_t$.

 To carry out the induction step, it is enough to show that if we join the extremal points of $\mathscr F_t$ by an arc, viewed as being below (resp. above) $L$ depending on 
 whether $t$ is odd (resp. even), then the meander becomes a single edge passing through all points of $\mathscr F_t$ joining the two fixed points.  If $t=s$ this results from the property $T(i)$ of an anti-Toeplitz involution, taking account of the arcs below $L$ if $s$ is odd, resp. above  $L$ if $s$ is even, which are  trivially modified.  If $2\le t<s$, then by the induction hypothesis, we can dispense with the pair of edges defined by the meander $N_{t+1}$ and replace it by an arc,  viewed as being below (resp. above) $L$ depending 
on whether $t+1$ is even (resp. odd), joining the extremal points of $\mathscr F_{t+1}$.  Then we just need to show that the resulting meander is a loop joining all the points of $\mathscr F_t \setminus \mathscr F^\prime_{t+1}$.

  However this assertion again results from the property $T(i)$ of an anti-Toeplitz involution, taking account of the arcs below (resp. above) $L$ which are trivially modified.


Consider now the passage from step $t=2$ to step $t=1$. If the integer points lying in $\mathscr F_1\setminus \mathscr F'_2$ are joined by arcs below $L$ defined by $\kappa^-$, then $T(i)$ of an anti-Toeplitz involution implies that there is a single loop (formed by the arcs below $L$ defined by $\kappa^-$ and the arcs above $L$ defined by $\sigma^+$) passing through all points of $\mathscr F_1\setminus\mathscr F'_2$. Note that $\mathscr I\subsetneq\mathscr F_1$ and that the modified integer arcs below $L$ meeting points in $\mathscr I$ coincide with those defined by the above loop if these points lie in $\mathscr F_1\setminus \mathscr F_2$ or with arcs below $L$ belonging to the set of arcs of $N_2$ otherwise. 
 Since $N_2$ is a union of two edges passing through all points of $\mathscr F_2$ we deduce that the part of $N_1$ formed by modified integer arcs above and below $L$ meeting points in $\mathscr I$  is still a union of two edges passing through all points of $\mathscr I$
 and precisely is contained in $M_{f_1}\cup M_{f_2}$ where $F=\{f_1,\,f_2\}$. Moreover the first arc at which $M_{f_i}$ ($i=1,\,2$) leaves $\mathscr I$ is an arc above $L$ and joins a point in $\mathscr I$ to a point in $\mathscr F_1\setminus\mathscr I$.
 Hence $(ii)$ and $(iii)$ in this case.
 


 Finally suppose that $s=1$ and $F$ is non-empty.  
 Then by a similar argument as above, 
$T(i)$ implies that the part of $N_1$ formed by modified integer arcs above and below $L$ meeting points of $\mathscr I$ is one edge or a union of two edges passing through all points of $\mathscr I$ and precisely is contained in $\cup_{f\in F}M_f$. Moreover the first arc at which one meander $M_f$ ($f\in F$) leaves $\mathscr I$ is an arc above $L$ and joins a point of $\mathscr I$ to a point of $\mathscr F_1\setminus\mathscr I$. 

 Hence $(ii)$ and $(iii)$ in this case too.

 For $(iv)$ we can suppose that $M_f$ first meets 
 $\mathscr F_1\setminus\mathscr I$ at $f_k$ lying to the left of 
$\mathscr I$. Equivalently that $f_k< c$.  Let $c^-$ be the centre of a component $J_{\ell}^-$ of $\mathscr J^-$ containing $f_k$. We can assume that $|J_{\ell}^-|>1$, otherwise $f_k$ is a fixed point of $\sigma^-$, hence a turning point. Then by \ref {3.7} and the definition of $\ell_1$, it follows that $c^-\leq \ell_1-1 \leq f_k$. Moreover 
if $f_k=c^-$ then 
 the latter is an integer and the half-integer arc below $L$ joining $c^-\pm \frac{1}{2}$ is not marked (since this arc belongs to the same half-integer meander as the half-integer arc above $L$ meeting $\ell_1-\frac{1}{2}=c^-+\frac{1}{2}$ which is marked by definition of $\ell_1$) forcing $c^-$ to be a fixed point of $\sigma^-$ (by \ref{2.5}), hence a turning point.  On the other hand if $c^-<f_k$, then the required assertion follows from \ref {3.4.5.2}$(v)$ applied to the pair $J_{\ell}^-$,\,$J^+$.
(See Figures \ref{Figure 4} and \ref{Figure 5}).

 \end {proof}





\textbf{Remark}.   If $c$ is integer then it is a fixed point of either $\sigma^+$ or of $\sigma^-$.  Thus $F$ cannot be empty in this case. However if $c$ is half-integer then $F$ can be empty.  This arises for the parabolic in $\mathfrak {sl}_6(\bf k)$ whose Levi factor has just $\alpha_3$ as its simple root.



\begin{cor}
Let $E$ be a non-trivial edge having no internal turning point. Then $E=\Phi(I)$ with $I=I_{1,\,e}$ a simple interval. We define $\epsilon_I=\epsilon_1=\epsilon_{e-1}$. Then $\epsilon_I\iota_I\in-R_*$.

\end{cor}

\begin{proof}

The first assertion of the corollary is obvious. Moreover $\varphi(1)$ is a source (resp. a sink) and $\varphi(e)$ is a sink (resp. a source) since all directed arcs forming the edge $E$ are aligned. This means that $\epsilon_I\iota_I=\varepsilon_{\varphi(1)}-\varepsilon_{\varphi(e)}$ (resp. $\epsilon_I\iota_I=\varepsilon_{\varphi(e)}-\varepsilon_{\varphi(1)}$) since $\epsilon_1=1$ (resp. $\epsilon_1=-1$) by Corollary \ref{3.3}. 
In particular $\epsilon_I\iota_I$ is a root. Moreover since the signs $\epsilon_i$ ($i\in I$) do not change in the simple interval $I$ (and are all equal to $\epsilon_I$) and since every $\epsilon_i\beta_i\in R$ ($1\le i\le e-1$), one has $\epsilon_I\iota_I\in R$ since $R$ is additively closed in $\Delta$. Finally $\epsilon_I\iota_I\not\in M$. Indeed if $\epsilon_I\iota_I\in M$ then $\varphi(1)$ and $\varphi(e)$ would belong to the same component $\mathscr I$ of the double partition (\ref{0.4}), which by \ref{3.4.5} would have to be equicentral. Then by Proposition \ref{3.8.4} $(ii)$ and $(iv)$ if one end-point of $E$ lies in $\mathscr I$, then the other one must lie in the external region $\mathscr E=[1,\,n]\setminus\mathscr I$ since $E$ has no internal turning point. It follows that $\epsilon_I\iota_I\in R\setminus M=-R_*$ by Lemma \ref{1.8}$(vi)$.

\end{proof}

\section{Compound Intervals.}\label{4}

Let $S$ be the set of arcs defined as in \ref {3.2}. Throughout this section we fix an edge $E$ of $S$.  The points of $E$ form a subset of $[1,n]$ whose cardinality we denote by $e$. In view of Corollary \ref{3.8.4}, we can assume that the edge $E$ has at least one internal turning point. This forces $\mathfrak q_{\pi^+,\,\pi^-}$ to be a proper biparabolic subalgebra, equivalently that $\pi^+\cap \pi^- \varsubsetneq \pi$ by \ref{3.4.5.3}(vii).

\subsection{Arcs, Order Relations and Gaps.}\label{4.1}

Recall the definitions in \ref{1.2bis}.


Given any (integer) arc above or below $L$ recall the notion (\ref {1.1}) of its support. Recall $\hat{\pi}^\pm$ defined in \ref {1.8} and observe that under the standing hypothesis of \ref {1.8} one has $\hat{\pi}^+ \cap \hat{\pi}^- = \emptyset$.

Let $A^+$ denote the set of unmodified integer arcs above $L$ defined by $\kappa^+$.  Let $<$ be the order relation on $A^+$ defined by inclusion of supports. Then any two elements with non-empty support are either strictly linearly ordered or have disjoint supports.  Let $\subseteq$ be the order relation on $A^+$ defined by the intersection of supports with $\hat{\pi}^-$. The linear order defined by $\subseteq$ may fail to be strict.  Similar definitions apply to the unmodified integer arcs $A^-$ below $L$.  These order relations may be also defined on the modified integer arcs; but they may fail to be linear due to crossings.  

  A gap above (resp. below) $L$ is defined to be a connected component $\hat{\pi}^+_j$ (resp. $\hat{\pi}^-_j$) of $\hat{\pi}^+$ (resp. $\hat{\pi}^-$).  We say that an arc attains a given gap if its support has a non-empty intersection with that gap. An arc above (resp. below) $L$ can never attain a gap above (resp. below) $L$; but an arc above (resp. below) $L$ can attain a gap below (resp. above) $L$.  (Except in \ref {8.6} and onwards.)

\subsection{Nil and Boundary Values.}\label{4.2}

Call a root $\alpha \in \Delta$ nil if supp $\alpha \cap (\hat{\pi}^+\cup\hat{\pi}^-) \neq \emptyset$, equivalently (Lemma \ref {1.8}$(i)$,\, $(ii)$) that $\alpha \in K \cup -K$.  By Lemma \ref {1.8}$(v)$, the nil elements of $R$ lie in $R\setminus M=-R_*$, which are the roots of the nilradical of $\mathfrak q^*$. This fact is the origin of our terminology.

Thus the subset $-R_*$ of nil roots is additively closed in $\Delta$, whereas the full set $K\cup -K$ of nil roots is not. 

Call $t \in [1,\,e-1]$ a nil point if $\beta_t$ (definition \ref{1.5}) is a nil root, and a boundary point if $\varphi(t)$ or $\varphi(t+1)$ is a turning point.  Call $\beta_t$ a nil (resp. boundary) value if $t$ is a nil (resp. boundary) point.  Since $\epsilon_t\beta_t \in S \subset R$, it follows that $\beta_t$ is a nil value if and only if $-\epsilon_t\beta_t \in R_*$.
The arc above or below $L$ which represents a nil value $\beta_t$ is called a nil arc.

\subsection{Odd Compound Intervals.}\label{4.3}

A compound interval $I$ is called odd if it is the union of an odd number of simple intervals. A simple or a compound interval $I$ is called nil if $\iota_I \in K\cup -K$.\par
Let $I=I_{r,\,s}$ ($1\le r<s\le e$) of the edge $E$ which is simple or an odd compound interval.\par
We will first show that, when $I$ is {\it nil}, it verifies necessarily $\epsilon_I\iota_I\in -K$
where $\epsilon_I$ is defined to be equal to the common sign $\epsilon_r=\epsilon_{s-1}$ (see \ref{3.3}).\par
Then we will show, when $I$ does not contain both end-points of the edge $E$, that $I$ is nil : it may be not true for $I$ such that $\Phi(I)=E$ since in this case the end-points of $E$ may belong to the same equicentral component of the double partition (actually we will show in \ref{8.4} that the latter is not possible if we also assume that the end-points of $E$ coincide with the end-points of the straightened edge $E^*$ associated to $E$ - see proof of Corollary \ref{8.4}).\par

In the index one case (\cite{FJ1}), $\hat\pi^-$ is reduced to one element $\alpha$ and the Levi factor of $\q$ is formed by two coprime blocks - that is why we call this case the coprime two-block case.
We have shown in \cite[2.8]{FJ1} that every simple interval $J$ is such that $\iota_J$ has $\alpha$ in its support.
Thus $\alpha$ occurs with alternate signs in each simple interval of a compound interval $I$  and then  remains in the support of $\iota_I$ if $I$ is odd. Hence an odd compound interval $I$ is nil in the coprime two-block case.

The fact that a simple interval being nil implies that every odd compound interval is nil could have failed even in the parabolic case because of crossings of the modified arcs which may then not be linearly ordered.

In our original proof lasting over some several pages, we were able to overcome this difficulty.
However by focusing attention on sources and sinks rather than the possible cancellation we were able to reduce the proof to the few lines given below.

\subsubsection{}\label{4.3.1}

\begin{lemma} Let $I=I_{r,\,s}$ be a simple or an odd compound interval such that $\Phi(I)\subset E$ (possibly $\Phi(I)=E$) which is nil.
Then
$$\epsilon_I\iota_I\in -K.$$

\end{lemma}

\begin{proof} Since $I$ is a simple or an odd compound interval, one element among $\varphi(r)$ and $\varphi(s)$ is a source and the other one is a sink of the edge $E$. Moreover since $I$ is assumed to be nil, $\varphi(r)$ and $\varphi(s)$ do not belong to the same component of the double partition (by Lemma \ref{1.8}$(iv)$ and \ref{1.9}).
Denote by $\mathscr I_r=J_r^+\cap J_r^-$ (resp. $\mathscr I_s=J_s^+\cap J_s^-$) the component of the double partition to which $\varphi(r)$ (resp. $\varphi(s)$) belongs. If these components are not linked (\ref{3.4.5.4}) then $\varepsilon_{\varphi(r)}-\varepsilon_{\varphi(s)}\in K\cap -K$ and in particular $\epsilon_I\iota_I\in -K$.

Assume now that $\mathscr I_r$ and $\mathscr I_s$ are linked and that $\varphi(r)<\varphi(s)$. First assume that $J_r^-=J_s^-$ (and then $J_r^+\neq J_s^+$). Then $\varepsilon_{\varphi(r)}-\varepsilon_{\varphi(s)}\in K^+$ by Lemma \ref{1.8}$(i)$. Denote by $c_r^+$  (resp. $c_s^+$) the centre of $J_r^+$ (resp. $J_s^+$) and by $c^-$ the centre of  $J_r^-=J_s^-$.
Then $c_r^+<c_s^+$ and by \ref{3.4.5} the only cases that can occur are when $c_r^+<c^-<c_s^+$ or $c_r^+=c^-<c_s^+$ or $c_r^+<c^-=c_s^+$ since in other cases $\varphi(r)$ and $\varphi(s)$ would be both sources or both sinks, which is not possible (see above).
When $c_r^+<c^-<c_s^+$, $\varphi(r)$ is a sink and $\varphi(s)$ is a source by Lemmas \ref{3.4.5.1}$(i)$, $(iii)$, \ref{3.4.5.2}$(iv)$ and $(vi)$. If $c_r^+=c^-<c_s^+$ then $\varphi(s)$ is a source by Lemma \ref{3.4.5.2} and by Lemma \ref{3.4.5.3}$(vii)$ $\varphi(r)$ is an end-point of $E$ which is necessarily a sink. If $c_r^+<c^-=c_s^+$ then $\varphi(r)$ is a sink by Lemma \ref{3.4.5.1} and by Lemma \ref{3.4.5.3}$(vii)$ $\varphi(s)$ is an end-point of $E$ which is a source.
In all these cases one has $\epsilon_r=\epsilon_{s-1}=-1$ by Corollary \ref{3.3} since moreover $\varphi(r)$ is not the finishing point of $E$. Then $\epsilon_I\iota_I=\epsilon_r(\varepsilon_{\varphi(r)}-\varepsilon_{\varphi(s)})\in-K^+$.\par
 
 We use a similar argument if $\varphi(r)>\varphi(s)$ with the same hypothesis that $J_r^-=J_s^-$. In this case we obtain that $\epsilon_r=\epsilon_{s-1}=1$. Hence again
 $\epsilon_I\iota_I\in-K^+$. Finally if $J_r^+=J_s^+$ we prove similarly that $\epsilon_I\iota_I\in-K^-$.

\end{proof}

\subsubsection{}\label{4.3.2}

\begin {lemma} Let $I=I_{r,\,s}$ be an odd compound interval, or a simple one of the edge $E$ such that $\Phi(I)$ does not contain both end-points of $E$ (that is $\Phi(I)\subsetneq E$). Then $\epsilon_I \iota_I \in -K$.
\end {lemma}

\begin{proof}

Suppose that $I$ is non-nil. By Lemma \ref{1.8}$(iv)$, this means that $\iota_I=\varepsilon_{\varphi(r)}-\varepsilon_{\varphi(s)}\in M$ and by \ref{1.9} the latter implies that $\varphi(r)$ and $\varphi(s)$ belong to the same component $\mathscr I=J^+\cap J^-$ of the double partition.
Since $I$ is odd (possibly simple), one element among $\varphi(r)$ and $\varphi(s)$ is a source and the other one is a sink. Then, by Lemmas \ref{3.4.5.1}$(i)$, $(iii)$, \ref{3.4.5.2}$(iv)$ and $(vi)$, $\mathscr I$ must be equicentral. However this is excluded by Lemma \ref{3.4.5.3}$(vii)$ through the hypotheses of the Lemma.
Hence $I$ is nil and the conclusion follows from the previous lemma.
\end{proof}

\begin{cor}  Let $I$ be a simple interval of $E$. Then $\epsilon_I\iota_I\in-R_*$.

 \end{cor}
 
 \begin{proof}

 Write $I=I_{r,\,s}\,:\,1\leq r<s\leq e$.  Since by Lemma \ref {3.3}, the $\epsilon_i$ do not change sign in $I_{r,\,s}$, the elements of its $\beta$-support either all lie in $R$ or in $-R$. 
 Since $R$ is additively closed in $\Delta$ and $E$ is an edge, one has $\iota_I \in R$ (resp. $\iota_I \in -R$) and then $\epsilon_I\iota_I\in R$ by definition of $\epsilon_I$. To conclude it suffices to note that $R\cap -K=-R_*$ by Lemma \ref{1.8}$(v)$ and $(vi)$ and to apply the previous Lemma.
 
 \end{proof}

\section{Construction of Nil Elements.}\label{5}

A crucial issue in our analysis is to ensure that condition b) of \ref {6.1} is satisfied. This was already delicate (\cite {FJ1}) in the coprime two-block case.  Here we develop some basics results towards this goal.

\subsection{}\label{5.1}

We continue to consider a single edge $E$ retaining the notation of \ref {3.3}.  We assume that $E$ is a 
non-trivial edge having at least one internal turning point.

Let $[r,\,s]$ be a connected subset of $[1,\,n]$.  Set $\pi_{[r,\,s]}:=\{\alpha_i\mid\! i,\,i+1 \in [r,\,s]\}$. In particular $\pi=\pi_{[1,\,n]}$.

\subsection{}\label{5.3}

Choose $t \in ]1,\,e[$ such that $\varphi(t)$ is an internal turning point of $E$.  Choose $t_+>t$ (resp. $t_-<t$) be minimal (resp. maximal) so that $\varphi(t_+)$ (resp. $\varphi(t_-)$) are turning points (not necessarily internal).  Let $I_-=I_{t_-,\,t}$ and $I_+=I_{t,\,t_+}$, be the simple intervals they define and $\iota_-$ and $\iota_+$ the corresponding interval values.  Set $\beta_-=\beta_{t-1},\,\beta_+=\beta_t$, which are boundary values.  Set $\epsilon_-=\epsilon_{t-1},\, \epsilon_+=\epsilon_t$. We obtain the following lemma, which  is a generalization of \cite [Lemma 3.2]{FJ1}.

\begin {lemma}  Suppose that $\beta_-$ (resp. $\beta_+$) is not nil.  Then $\epsilon_- (\beta_- +\iota_+)$ (resp. $\epsilon_+ (\beta_+ +\iota_-)$) belongs to $K$.
\end {lemma}

\begin {proof}   Suppose for example that $\beta_-$ is not nil.  Since  $\varphi(t)$ is a turning point $\epsilon_-\epsilon_+=-1$ (by Lemma \ref{3.3}). Thus $\epsilon_-(\beta_- +\iota_+)=\epsilon_-\epsilon_+(\epsilon_+\beta_- +\epsilon_+\iota_+)\in R_* \subset K$, by Corollary \ref {4.3.2} and Lemma \ref{1.8}$(x)$.
\end {proof}

\subsection{}\label{5.4}

In the situation studied in \cite {FJ1}, $\beta_-$ and $\beta_+$ cannot both be nil values.  This holds more generally for a parabolic subalgebra, that is when $\pi^+=\pi$.  Indeed if $\beta_-$ is given by an arc above the horizontal line then $\beta_+$ must be given by an arc below the horizontal line.   Thus for $\beta_+$ to be a nil value it must have a non-zero coefficient in the set $\pi^-\setminus (\pi^+\cap \pi^-)$, which for a parabolic, is the empty set.

Obviously these considerations fail in the biparabolic case.  Nevertheless we have the following slightly surprising result.   Retain the notations and conventions of \ref {5.3}.

Let $J=[\ell,\,r]$ be a connected subset of $[1,\,n]$.  We write $\ell(J)=\ell$ (resp. $r(J)=r$) called its leftmost (resp. rightmost) element.

\begin {lemma}  Suppose that $\beta_-$ and $\beta_+$ are both nil.  Then $\epsilon_-(\beta_-+\iota_+)$ and $\epsilon_+(\beta_++\iota_-)$ belong to $K$.
\end {lemma}

\begin{figure}[!h]
\centering
\input{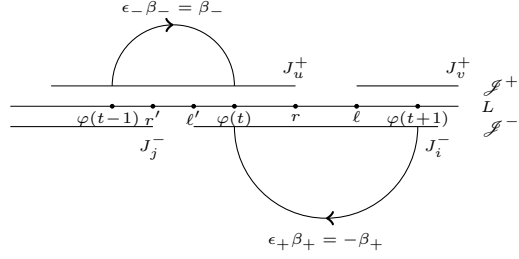}
\caption{\scriptsize{\it 
 Illustration of the proof of Lemma \ref{5.4}.}
Here we denote $r(J_u^+)$ by $r$, resp. $r(J_j^-)$ by $r'$ and $\ell(J_v^+)$ by $\ell$, resp.
$\ell(J_i^-)$ by $\ell'$.
The support of $\beta_+$ has a non-empty intersection 
with $\pi_{[r,\,\ell]}$
 hence attains a gap above $L$
  whilst the support of $\beta_-$ has an empty intersection with it 
  because it cannot attain a gap above $L$ by \ref{4.1}. 
 Similarly the support of $\beta_-$ has a non-empty intersection 
 with $\pi_{[r',\,\ell']}$, whilst the support of $\beta_+$ 
 has an empty intersection with it.}
 \label{Figure 8}
 \end{figure}

\begin {proof}  We can assume without loss of generality that $\beta_-$ comes from an arc above the horizontal line and hence $\beta_+$ comes from an arc below the horizontal line.

With respect to the natural order on the horizontal line, that is on $[1,\,n]$, we can assume without loss of generality that $\varphi(t-1) < \varphi(t)$, and hence since $\varphi(t)$ is a turning point that $\varphi(t) < \varphi(t+1)$. This means that the arrow on the arc defining $\beta_-$, resp. $\beta_+$, points from left to right (here we consider the elements $\beta_i$ and not $\epsilon_i\beta_i$).

Let $J_u^+$ (resp. $J_i^-$) be the connected component of $\mathscr J^+$ (resp. $\mathscr J^-$) containing $\varphi(t)$.  Then let $J_v^+$ (resp. $J_j^-$) be the connected component of $\mathscr J^+$ (resp. $\mathscr J^-$) of cardinality $>1$ immediately to the right (resp. left) of $J_u^+$ (resp. $J_i^-$) if there exists. Recall the above notations 
and set $\ell(J_v^+)=n$ if $J_v^+$ does not exist (resp. $r(J_j^-)=1$ if $J_j^-$ does not exist).
Observe that, since $\beta_+$ (resp. $\beta_-$) is nil, $r(J_u^+)<n$ (resp. $\ell(J_i^-)>1$).

Since $\varphi(t) \leq r(J_u^+)$, it follows that
$$\text {supp} \ \beta_- \cap \pi_{[r(J_u^+),\,\ell(J_v^+)]} = \emptyset.\eqno {(*)}$$

The arrow on every arc below the horizontal line corresponding to an element of $\beta_s\,:\, s \in I_+$ points from left to right.  On the other hand no element of $\beta_s\,:\, s \in I_+$ above the horizontal line can cross the interval on $[1,\,n]$ separating $J^+_u$ and $J_v^+$ (resp. $J_u^+$ and $n$), if $J_v^+$ exists (resp. otherwise) that is to say that no element of $\beta_s\,:\, s \in I_+$ above $L$ can attain the gap $\pi_{[r(J_u^+),\,\ell(J_v^+)]}$ of $\hat\pi^+$ above $L$ (\ref{4.1}).  Hence
$$\text {supp} \ \iota_+ \supset \text {supp} \ \beta_+ \cap \pi_{[r(J_u^+),\,\ell(J_v^+)]}. \eqno {(**)}$$

On the other hand $\beta_+=\varepsilon_{\varphi(t)}-\varepsilon_{\varphi(t+1)}$, whilst $\varphi(t) \leq r(J_u^+)$.  Then the assumption on $\beta_+$ being nil forces the right hand side of $(**)$ to be non-empty.  

Take $\alpha$ in the right hand side of $(**)$.  Then $\alpha \in \pi^-\setminus (\pi^-\cap \pi^+)$.  Yet $\epsilon_+\beta_+ \in S$, so $\alpha$  occurs with a negative coefficient in $\epsilon_+\beta_+ \in S$. Yet $\epsilon_-\epsilon_+=-1$, so by $(*)$ it occurs with a positive coefficient in $\epsilon_-(\beta_-+\iota_+)$.  Hence $\epsilon_-(\beta_-+\iota_+) \in K$ by Lemma \ref{1.8}$(i)$.

Similarly $\epsilon_+(\beta_++\iota_-) \in K$, via the pair $J_i^-,J_j^-$. (See Figure \ref{Figure 8}).

\end {proof}

\subsection{}\label{5.6}

Retain the above notation and hypotheses.  From Lemma \ref {4.3.2} we immediately obtain the following extension of Lemma \ref {5.3}.  Let $I:=I_{r,s}$ be an odd compound interval. 

\begin {prop}  Suppose $\beta_{r-1}$ (resp $\beta_s$) is defined and non-nil (so in particular $I$ is not the whole edge). Then $\epsilon_{r-1}(\beta_{r-1}+\iota_I)\in K$ (resp. $\epsilon_s(\beta_s + \iota_I) \in K$).
\end {prop}

\begin{proof}
Follows from Lemma \ref{4.3.2} and from Lemma \ref{1.8}$(vii)$.

\end{proof}

\section{The Construction of a New Simple System.}\label{6}

\subsection{The goal.}\label{6.1}

Recall (see \ref{2.4} and \ref {3.3}) that after modifying the link patterns the resulting meanders    
decompose into a finite union of non-trivial edges and as in \ref {3.3} we just consider one non-trivial edge $E$ and retain the notation already introduced.

Then $\Pi:=\{\beta_i\}_{i \in [1,\,e-1]}$ is a simple root system of type $A_{e-1}$.  Unfortunately it is not what we need ! 
(except if $E$ has no internal turning point).

This is because the set describing the eigenvector of the adapted pair is $S_E=\{\epsilon_i\beta_i\}_{i\in [1,\,e-1]}$ rather than $\Pi$.  (Here $S_E$ is a subset of $S$ defined in \ref {3.2} and was defined in \ref{3.3}.) What we need to do is to find a new simple root system $\Pi^*=\{\beta^*_i\}_{i \in [1,\,e-1]}$ with the following two properties :

\begin{enumerate}

\item[a)]  $\epsilon_i\beta_i\ \in \mathbb N\Pi^*$, for all $i \in [1,\,e-1]$,

\item[b)]  If $\beta^*_i \neq \epsilon_i\beta_i$, then $\beta^*_i \in K$.

\end{enumerate}

It is the second condition which causes all our present misery. Indeed recall \ref {3.2} and let $\mathfrak h^\prime_{\Lambda}$ be a subspace of $\mathfrak h_{\Lambda}$ of dimension $|S_E|$ separating the elements of $S_E$. Then we could satisfy a) through the time-honoured method of fixing positive integers $c_i: i \in [1,e-1]$ for which the relation $h(\epsilon_i\beta_i)=c_i$ defines a regular element of $\mathfrak h^\prime_{\Lambda}$.  Then $\{\alpha \in \Delta\mid h(\alpha)>0\}$ is a set of positive roots for $\Delta$ and so defines a simple root system $\Pi^*$, which satisfies a) by construction. However imposing b) can result in there being no solution at all (for example if $K=\emptyset$ and not all the $\epsilon_i$ have the same sign) and even when a solution exists the choice of the $c_i$ is extremely delicate.


Since the case when $E$ has no internal turning point needs nothing to do more,
we now assume that the edge $E$ has at least one internal turning point (as in Sections \ref{4} and \ref{5}) and then $e>2$.

\subsection{Methodology.}\label{6.2}

In many ways our solution resembles that given in \cite {FJ1}.  However it is complicated by the fact that it is much more difficult to ensure no cancellations in odd compound intervals and indeed this was 
why all the preliminary material in earlier sections was needed.  On the other hand since we are forced to be more systematic we understand much better our solution found in \cite {FJ1} and indeed what we present here will be now much more natural.

 The data describing the edge $E$ is presented in the form of a central vertical line $V$ with $e$ equally spaced points with 
$\varphi(1)$ (resp. $\varphi(e)$) at the top (resp. bottom) end.   The short central line joining  $\varphi(i),\,\varphi(i+1)$, is given an arrowhead pointing downwards (resp. upwards) if $\epsilon_i$ equals $1$ (resp. $-1$), and represents $\epsilon_i\beta_i$.

What we shall do in order to obtain $\Pi^*$ from $\Pi$ will closely resemble what we did to obtain $\Pi$ from $\pi$.  Thus we shall mark certain elements of $\Pi$, that is to say some of the $\beta_i$.  In analogy with \cite [3.3]{FJ1} these markings will then define a left (resp. a right) signature and from this data we will construct left (resp. right) directed external arcs joining pairs of points on the vertical central line $V$ and at the same time remove the same number of the short central lines.

A directed arc going from the points $i$ to $j$ on $V$ will be designated as $(i,\,j)$ and drawn with an arrowhead pointing in the direction of $j$. The same convention will apply to the short central lines which remain.  Thus 
$(\varphi(i),\,\varphi(i+1))$, resp. $(\varphi(i+1),\,\varphi(i))$, will designate $\epsilon_i\beta_i$ if $\epsilon_i$ equals $1$, resp. $-1$ (observe that in the first case, $\epsilon_i\beta_i=\varepsilon_{\varphi(i)}-\varepsilon_{\varphi(i+1)}$ and in the second case, $\epsilon_i\beta_i=\varepsilon_{\varphi(i+1)}-\varepsilon_{\varphi(i)}$).  They will be referred to as short directed arcs.  

Observe also that the notation of a directed arc as above is somehow consistent with the notation of \ref{1.1} of a directed arc below or above the horizontal line $L$ though here the directed arcs refer to $V$ whilst the arcs in \ref{1.1} refer to $L$.
As in \ref{1.1} we will call the points $i$ and $j$ of $V$ of a directed arc $(i,\,j)$ defined as above the end-points of the arc $(i,\,j)$. More generally all the terminology defined for directed arcs above or below the horizontal line $L$ may be used here for directed arcs whose end-points will belong to the vertical central line $V$.\par
 Moreover in the previous Sections the directed arcs above and below the horizontal line $L$ represented the $\epsilon_i\beta_i$ for $1\le i\le e-1$, and here in this Section the short directed arcs (on the vertical central line $V$) will also represent the $\epsilon_i\beta_i$ for $1\le i\le e-1$. To go further with the analogy, recall that each $\epsilon_i\beta_i$ represented by a directed arc above or below the horizontal line $L$ came from an unmodified integer arc (encircling an half-integer arc of the Kostant cascade) which was  trivially modified (that is not modified) if the corresponding half-integer arc was not marked or was non-trivially modified if the corresponding half-integer arc was marked, so that the directed arcs obtained by this procedure represent the elements $\epsilon_i\beta_i$ of the set $S_E$ which are the weights of the vectors whose sum is equal to the element $\eta_E$ as defined in \ref{3.3}.
 
Now generalize the definition given in \ref{1.3} of a meander and more precisely of an edge, which will designate a set of points $\{i_1,\,,\ldots,\,i_p\}$ 
(for $p\ge 2$) of $[1,\,n]$ such that there exists exactly one directed arc meeting $i_{j}$ and $i_{j+1}$ for all $1\le j\le p-1$ and no directed arc meeting $i_u$ and $i_t$ if $\mid t-u\mid\neq 1$.
Such an edge may be represented on the vertical central line $V$ or on the horizontal line $L$.
Observe that an edge defined as above and represented on the horizontal line $L$ has its arcs not necessarily alternating above and below $L$ whilst the arcs of an edge defined as in \ref{1.3} do. The points $i_1$ and $i_p$ are called the end-points of the edge and when the directed arcs have all their arrows aligned (\ref{1.4}) we call $i_1$ (resp. $i_p$) the starting point and $i_p$ (resp. $i_1$) the finishing point of the edge if these arrows point from $i_1$ to $i_2$ (resp. from $i_2$ to $i_1$) : that is to say when the arrows of the edge are all aligned, the starting point will correspond to a source of this edge and the finishing point to a sink, in the sense of \ref{1.4}.

In this Section we will remove some of the short directed arcs (those which will be marked) and replace them by directed external arcs representing roots which will be shown to belong to $K$, so that all the new directed arcs (that is the directed external arcs we add together with the short directed arcs which remain) have their arrows all aligned (\ref{1.4}) and form a single edge in the sense given above, called the straightened edge $E^*$ of $E$. This straightened edge $E^*$ passes through the same points as $E$ but by a different path and then gives a system of simple roots of $\s\l_e(\bf k)$ such that there exists a regular nilpotent element $y_E$ of $\s\l_e(\bf k)^*$ whose restriction to $\q_{\Lambda}$ is equal to $\eta_E$ (see also \ref{9.1} for more details).

In Sections \ref{1}-\ref{3} we presented $\Pi$ 
 as obtained from $\pi$ by arcs above and below $L$, whilst we now extract $\Pi^*$ from $\Pi$ by an analogous but different procedure using arcs to the left and to the right of $V$.
 We show that a path following the new directed arcs will be a single edge (denoted by $E^*$ and passing through the same points as $E$) and hence define the new simple root system $\Pi^*$ (for an $\mathfrak {sl}_e(\bf k)$ subalgebra of $\mathfrak g$).  The directed arcs which are not short will be shown to define elements belonging to $K$. This will establish b). It is the main technical difficulty. Finally as in \cite [3.7]{FJ1} a rather simple order relation will be introduced to verify a).

\textit{The boundary values that were changed in \cite [3.1]{FJ1} correspond to the marked $\beta_i$.  This will also be true here in the present more general context.}

It seems to us quite remarkable that a procedure not very different to \cite {FJ1} used in the very special index one case applies in the general biparabolic case with only some technical modifications.

\subsection{Isolated Values.}\label{6.3}

Following \cite [2.11]{FJ1} we define $t \in [1,\,e-1]$ to be an isolated point if both $\varphi(t)$ and $\varphi(t + 1)$ are turning points. If $t$ is an isolated point we call $\beta_t$ an isolated value.  
It follows from 
Corollary \ref{4.3} that an isolated value is nil.
Suppose that we are in the parabolic case.  By the remark in \ref {5.4} the neighbouring values of an isolated value must be non-nil.  In particular isolated values cannot be neighbours.

However in the biparabolic case there can be arbitrary long sequences of neighbouring isolated values.

\subsection{Choice of Markings and Their Consequences.}\label{6.4}


Let $A$ (resp. $B$) be the subset of the set $T$ of turning points of $E$ which are sources (resp. sinks), see \ref {1.4}. (A rather different definition of these sets was given in \cite [2.4]{FJ1} based on ``connectedness", which is not possible in the present more general situation, though it is a crucial fact that a vestige (\ref {3.4.5}) of this property remains.)

Notice that the above decomposition implies that turning points of type $A$ alternate with those of type $B$ (by \ref{1.4} again).  In particular an interval 
whose end-points are one of type $A$ and one of type $B$
 is always odd.  This brings to relevance Lemma \ref {4.3} which rather trivially fails for even compound intervals.  Then the way that the rather delicate property b) of \ref {6.1} is satisfied is by a judicial use of this lemma effected by the markings described below.

For each internal turning point $\varphi(t)$ we shall mark one or both of its boundary values according to the following rules which we first state informally.

\begin{enumerate}

\item[(i)] If neither boundary value is nil then either one or the other is marked, but \textit{not} both.

\item[(ii)] If exactly one boundary value is nil, then the opposite (non-nil) one is marked. 

\item[(iii)] If both boundary values $\beta_{t-1},\,\beta_t$ are nil (which can only happen in the biparabolic case) then one or both may be marked (in a manner to be specified below and in the next subsection, see $\textbf {(M)}$).

\end{enumerate}

Notice that the marking of a $\beta_i$ can be viewed as a marking of the short central line it represents.

In the parabolic case a marked value (which is non-nil in this case) is  the neighbour of exactly one internal turning point.
This is also true in the biparabolic case except in the case of a sequence strictly longer than 1
of isolated values. Except in this latter case \textbf{the marked value is assigned to the internal turning point in question.} In the latter case we assign a marked value to a unique internal turning point as follows.

Suppose we are in the biparabolic case and that we have a maximal sequence of $r+1>1$ internal turning points $\varphi(t),\,\varphi(t+1),\,\ldots,\,\varphi(t+r)$ which are neighbours.  We must mark exactly  
$r+1$ elements of the set $\{\beta_{t+j}\}_{j=-1}^{r}$.  Since $\{\beta_{t+j}\}_{j=0}^{r-1}$ are all nil values, we can assume without violating the above rules that $\{\beta_{t+j}\}_{j=-1}^{r-2}$ are all marked.  If $\beta_{t+r}$ is non-nil, then we are forced by (ii) above to mark it.  Otherwise we mark $\beta_{t+r-1}$.

For $j=-1,\,0,\,\ldots,\,r-2$, we assign the marked value $\beta_{t+j}$ to the turning point below it, namely to $\varphi(t+j+1)$.  We assign the last marked value to $\varphi(t+r)$.  This marked value may be above or below $\varphi(t+r)$.

 Since the arrows point in opposite directions on the short 
central lines meeting an internal turning point, one of these lines has to be deleted. Then we replace the marked lines by arcs to align all arrows. The precise algorithm is given in the next section.  One may consider that it is essentially canonical, that is few other choices are possible.

\subsection{The Signature Rule.}\label{6.5}

First we add a final precision to rule (iii) given in \ref {6.4}.  Then we shall describe how the left and right 
directed external arcs are constructed.

Suppose $\varphi(t)$ is an internal turning point.

If the marked value assigned to $\varphi(t)$ lies above (resp. below) $\varphi(t)$, that is to say is $\beta_{t-1}$ (resp. $\beta_t$), we write $sg(t)=1$ (resp. $sg(t)=-1$) and we say that the signature at $t$ is positive (resp. negative).

The subset of internal turning points of type $A$, resp. $B$, having positive signature will be denoted by $A_+$, resp. $B_+$. The subset of internal turning points of type $A$, resp. $B$, having negative signature will be denoted by $A_-$, resp. $B_-$.

As in \cite {FJ1} we construct $\beta_i^*$ by adding to a marked $\beta_i$, an odd compound interval value on the opposite side of the internal turning point for which the marked value is a neighbour.  In the parabolic case, the marked value is non-nil and so the resulting $\beta^* _i$ is nil by Proposition \ref{5.6}.  However in the biparabolic case it can happen that the marked value is also nil and then if the compound value is not simple, the sum can be non-nil.  However in this case we have the choice of marking instead the second neighbouring value and we do this to ensure that only a simple interval value is added.

The precise rule is very simple.

$\textbf{(M)}$. If both boundary values $\beta_{t-1},\,\beta_t$ are nil, then we shall mark $\beta_{t-1}$ and assign it to the internal turning point $\varphi(t)$ that is $sg(t)=1$.  Note that this rule is compatible with the rule given in \ref {6.4} if 
$\beta_{t-1}$ is an isolated value or if $\beta_t$ is an isolated value.

We have now specified completely which short 
central lines are marked and have assigned to each marked short central line a unique internal turning point, which is always at one end of the marked line.  Conversely every internal turning point has been assigned to a unique marked short 
central line and lies at one of its ends. This defines a bijection between marked 
short central lines and internal turning points. With some exceptions in the biparabolic case this bijection is simply neighbour assignment, that is to say every marked short central line has a unique end-point which is an internal turning point (to which it is assigned) and every internal turning point is the end-point of a unique marked short central line.

 In addition to the above, signature has been assigned.  In terms of this we now describe how the $\beta^*_i$ (equivalently the left and right directed external arcs) are constructed.

Observe that, if $\varphi(t)$ is an internal turning point having negative signature, then our construction implies that $\beta_t$ is non-nil. (For this recall $\bf (M)$).

Step $(i)$.

 Let $\varphi(t)$ be an internal turning point of type $A$ having negative signature.  This means that $\beta_t$ is marked (and assigned to $\varphi(t)$). In this case we remove the short vertical central line defined by $\beta_t$ and replace it by a directed external arc, \textit{with the arrow pointing in the same direction as the deleted line}, joining the lower end-point of $\beta_t$, namely $\varphi(t+1)$, to the turning point $\varphi(s)$ of $B$ defined as follows. Suppose that there is a turning point of type $A_-$ strictly above $\varphi(t)$. Let $\varphi(r)$ be the first such turning point as we go up the central line $V$ and let $\varphi(s)$ be the first turning point of type $B$ below $\varphi(r)$. Otherwise let $\varphi(s)$ be the first turning point of type $B$, on going down $V$, which lies necessarily strictly above $\varphi(t)$.


Then $s<t$ and this operation replaces $\epsilon_t\beta_t$ by 
$\beta^*_t:=\epsilon_t(\beta_t+\iota_{I_{s,\,t}})$.  In our presentation these new arcs are drawn on the left of the vertical central line $V$.

 In the above $I_{s,\,t}$ may be a compound interval, but because of the alternation in turning point type (that is $A$ or $B$) it is always an odd compound interval.

 \begin {lemma}

 \begin{enumerate}

 \item[(i)] In the above construction $\beta^*_t \in K$.
 \end{enumerate}
 \end {lemma}

 \begin {proof}

 Recall that $\iota_{I_{s,\,t}}$ is an odd compound interval value (possibly simple).


By the remark above Step (i), $\beta_t$ is non-nil. Then the assertion follows from Proposition \ref {5.6} and Lemma \ref{5.3}. 


 \end {proof}

 Step $(ii)$.

  At every internal turning point $\varphi(t)$ of  type $A$ having positive signature we remove the short vertical central line defined by $\beta_{t-1}$ and replace it by a directed external arc drawn on the left of $V$, \textit{with the arrow pointing in the same direction as the deleted line}, joining the upper end-point of 
$\beta_{t-1}$, namely $\varphi(t-1)$, to the first turning point of type $B$ below $\varphi(t)$, denoted by $\varphi(s)$, which is defined because $\varphi(t)$ is an internal turning point. Then $s>t$ and this operation replaces $\epsilon_{t-1}\beta_{t-1}$ by $\beta^*_{t-1}:=\epsilon_{t-1}(\beta_{t-1}+\iota_{I_{t,\,s}})$.

  \begin {lemma}

\begin{enumerate}
 \item[(ii)] In the above construction $\beta^*_t \in K$.
 \end{enumerate}
 \end {lemma}

 \begin {proof}

  Since $\iota_{I_{t,\,s}}$ is a simple interval value, the assertion follows from Lemmas \ref{5.3} and \ref{5.4} since in this case, $\beta_{t-1}$ may be non-nil or if $\beta_{t-1}$ is nil (and marked) then $\beta_t$ is also nil by (ii) of \ref{6.4} and the rule describing the marking of a sequence of consecutive isolated values (see end of \ref{6.4}). 

  \end {proof}

  Step $(iii)$.  Exactly the same procedure is applied to the internal turning points of type $B$ drawing directed external arcs on the right of the vertical central line $V$.

  \bigskip
  
{\bf Terminology.}  
  \medskip
  
For every internal turning point $\varphi(t)$ the directed external arc replacing $\epsilon_t\beta_t$ if $sg(t)=-1$ or $\epsilon_{t-1}\beta_{t-1}$ if $sg(t)=1$ is said to be {\it assigned} to $\varphi(t)$ and denoted by $\mathscr A_t$ and the end-point of $\mathscr A_t$ which we have denoted above by $\varphi(s)$ and which is a turning point of opposite type to $\varphi(t)$ is said to be {\it crossed} by $\mathscr A_t$. Then if $sg(t)=1$ (resp. $sg(t)=-1$) the directed external arc $\mathscr A_t$ meets the points $\varphi(t-1)$ (resp. $\varphi(t+1)$) and $\varphi(s)$ but only crosses $\varphi(s)$ in the above sense. The other point which is joined by $\mathscr A_t$ will be said to be {\it simply met} by $\mathscr A_t$.
  
 \bigskip 

  \begin {lemma}  Consider a turning point $\varphi(s)$ of type $A$ (resp. type $B$).

 \begin{enumerate}
\item[(iii)] Then via the above construction at most two directed arcs (external or not) meet $\varphi(s)$ and at most one directed external arc crosses it. More precisely:

 \item[(iv)] If exactly one directed external arc 
meets $\varphi(s)$ via the above construction then either this directed external arc does not cross $\varphi(s)$ (in the above sense) and then both short directed arcs meeting $\varphi(s)$ (if $\varphi(s)$ is internal) are deleted (and replaced) and the direction of the arrow on this external arc is away from (resp. towards) $\varphi(s)$ or 
this directed external arc crosses $\varphi(s)$ (in the above sense) and then
 the direction of the arrow on this external arc is towards (resp. away from) $\varphi(s)$.
 In both cases this directed external arc lies to the right (resp. left) of $V$.

\item[(v)]  If two directed external arcs 
meet $\varphi(s)$ via the above construction then  they lie on the same side of $V$ and their arrows are aligned and only one crosses $\varphi(s)$ in the above sense.  Moreover in this both short 
directed arcs meeting $\varphi(s)$ (if $\varphi(s)$ is internal) are replaced (and deleted).
\end{enumerate}
 \end {lemma}

\begin {figure}[!h]
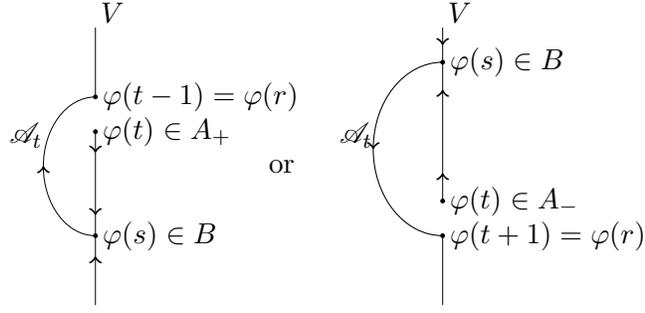

\centering


\begin{texgraph}[name=fig9,export=pgf]
  view(-4,22,-9,8.5),Marges(0.5,0,0,0.5),
  size(6.5), Width:=2,
  Horizont:=left,

DotSize:=1+i,
LabelSize:=small,

Seg(-8*i,2*i),
Seg(4*i,8*i),

LabelDot(-4*i,"$\varphi(s)\in B$","E",1,0.1),

LabelDot(2*i,"$\varphi(t)\in A_+$","E",1,0.1),

LabelDot(4*i,"$\varphi(t-1)=\varphi(r)$","E",1,0.1),

Width:=8,
flecher(Get(Seg(-8*i,-4*i)),[0.7]),

flecher(Get(Seg(-4*i,2*i)),[0.2+i]),

flecher(Get(Seg(-4*i,2*i)),[0.8+i]),

flecher(Get(ellipticArc(-4*i,0,4*i,3,4,-1)),[0.5]),

Width:=2,

ellipticArc(-4*i,0,4*i,3,4,-1),

LabelDot(8*i,"$V$","NE",0,0.1),

LabelDot(-4,"$\mathscr A_t$","N",0),

LabelDot(10,"or","E",0,0),
Seg(20+8*i,20-2*i),

Seg(20-4*i,20-8*i),

LabelDot(20+8*i,"$V$","NE",0,0.1),

LabelDot(20+6*i,"$\varphi(s)\in B$","E",1,0.1),

LabelDot(20-2*i,"$\varphi(t)\in A_-$","E",1,0.1),

LabelDot(20-4*i,"$\varphi(t+1)=\varphi(r)$","E",1,0.1),

ellipticArc(20+6*i,20+i,20-4*i,4,5,1),

LabelDot(15,"$\mathscr A_t$","N",0),

Width:=8,

flecher(Get(ellipticArc(20+6*i,20+i,20-4*i,4,5,1)),[0.5]),

flecher(Get(Seg(20+8*i,20+6*i)),[0.5]),

flecher(Get(Seg(20+6*i,20-2*i)),[0.2+i]),

flecher(Get(Seg(20+6*i,20-2*i)),[0.8+i]),

\end{texgraph}

\caption{\scriptsize{\it 
Illustration of Lemma \ref{6.5}$(iv)$.}
A short central line has been deleted 
and replaced by an external arc $\mathscr A_t$
 which crosses $\varphi(s)$ in the above sense. 
 At a further step at least one of the short
 central lines
  meeting $\varphi(s)$ is deleted 
 (if $\varphi(s)$ is internal).}
 \label{Figure 9}
 \end{figure}

\begin{figure}[!h]
\centering
\input{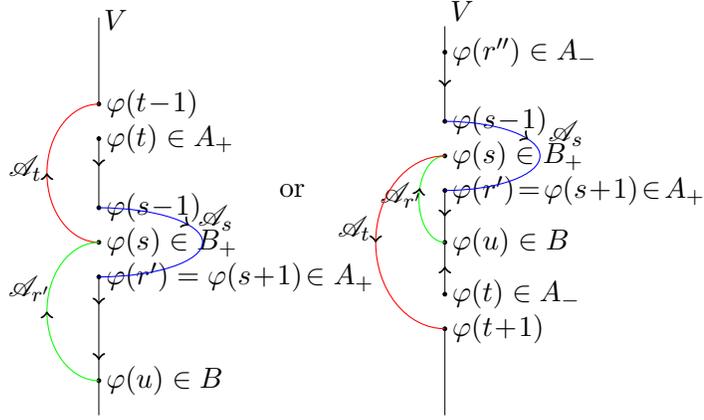}
\caption{\scriptsize{\it 
Illustration of Lemma \ref{6.5}$(v)$ with $\varphi(s)\in B$.}
 In the Figure to the left, resp. to the right, $\varphi(t)\in A_+$,
 resp. $\varphi(t)\in A_-$.
Both short directed arcs meeting $\varphi(s)$ 
 (if $\varphi(s)$ is internal)
 are deleted 
 and replaced
  by the external arcs 
$\mathscr A_s$ in blue and $\mathscr A_{r'}$ in green.
 The short central line joining $\varphi(t)$ and $\varphi(t-1)$,
  resp. $\varphi(t)$ and $\varphi(t+1)$ is deleted and replaced 
 by the external arc $\mathscr A_t$ (drawn in red) which crosses $\varphi(s)$. 
 The external arc $\mathscr A_{r'}$ does not cross $\varphi(s)$.
At a further step at least one of the short directed arcs
 meeting $\varphi(u)$ is deleted 
(if $\varphi(u)$ is internal).}
\label{Figure 10}
\end{figure}

 \begin {proof}

 These statements are easy consequences of our construction; but we shall spell out the details in the case when $\varphi(s)$ is of type $B$. In this case the external arc(s) appear on the left of $V$.

 $(iii)$. 
 Via the above construction the turning point $\varphi(s)$ cannot be crossed together by a directed external arc 
 $\mathscr A_t$ with $sg(t)=1$ and a directed external arc $\mathscr A_{t'}$ with $sg(t')=-1$. Actually this explains why the rule for negative signature is more complicated than for positive signature, in particular why the directed external arc assigned to a point of negative signature has to jump over all points of the same type but of positive signature before crossing the point of opposite type to it.   
 
 It is also not possible that two directed external arcs $\mathscr A_t$ and $\mathscr A_{t'}$ meet simply $\varphi(s)$. Indeed in this case $\varphi(s)$ should be internal and a neighbour of both $\varphi(t')=\varphi(s-1)\in A_-$ and $\varphi(t)=\varphi(s+1)\in A_+$. Then $\beta_{s-1}$ (resp. $\beta_s$) should be a marked isolated value assigned to $\varphi(s-1)$ (resp. to $\varphi(s+1)$). But this is in contradiction with the rules given above since 
$\beta_{s-1}$ should have to be marked and assigned to $\varphi(s)$ (by $\bf(M)$). 
 
 Finally the above construction  implies that at most two directed external arcs can meet $\varphi(s)$ (one which will cross it and the other which will meet it simply) and in case of two we will show that no short directed arc will meet it.
Since at least one of the short directed arcs meeting $\varphi(s)$ must be deleted (if $\varphi(s)$ is an internal turning point), it follows that there will be at most two directed arcs (external or not) meeting $\varphi(s)$.  

 $(iv)$.
 Assume that there is exactly one directed external arc meeting $\varphi(s)$ and denote it by $\mathscr A$.
Suppose that $\mathscr A$ crosses $\varphi(s)$. It means that $\mathscr A$ is assigned to an internal turning point $\varphi(t)\in A$ (and $\mathscr A=\mathscr A_t$, which lies to the left of $V$) and joins a neighbour $\varphi (r)$ of the turning point $\varphi(t) \in A$ (which may be $\varphi(t-1)$ or $\varphi(t+1)$ depending on signature) to the turning point $\varphi(s) \in B$.  Since $\varphi(t)$ is a source, our construction means that the arrow on the external  arc $\mathscr A_t$ 
 points from $\varphi(s)$ to this $\varphi(r)$, as required. (See Figure \ref{Figure 9}).
 
Suppose that $\mathscr A$ does not cross $\varphi(s)$. Then it meets simply $\varphi(s)$ and $\varphi(s)$ is the neighbour of an internal turning point $\varphi(r')\in A$. That is $\varphi(s)=\varphi(r'-1)$ (resp. $\varphi(s)=\varphi(r'+1)$) and $\mathscr A=\mathscr A_{r'}$ so to the left of $V$. Moreover 
 $\varphi(r')\in A_+$ (resp. $\varphi(r')\in A_-$) and
$\beta_{r'-1}$ (resp. $\beta_{r'}$) is an isolated value marked and assigned to $\varphi(r')$ (so in particular this phenomenon can only occur in the biparabolic case). But because of rules in case of isolated values, we cannot have $\varphi(r')\in A_-$ together with 
 $\varphi(s)=\varphi(r'+1)\in B$.
 Then $\varphi(s)=\varphi(r'-1)\in B$ (if $\varphi(s)$ is internal, necessarily $\varphi(s)\in B_+$) and $\varphi(r')\in A_+$ and the directed arc $\mathscr A=\mathscr A_{r'}$ points towards $\varphi(s)$, as required. Moreover (if $\varphi(s)$ is internal) both short directed arcs meeting $\varphi(s)$ must be deleted and replaced by external arcs $\mathscr A_{s}$ (to the right of $V$) and $\mathscr A_{r'}$ to the left of $V$. (See Figure \ref{Figure 13} on the right).

$(v)$.  By our construction two external arcs can only
 have $\varphi(s)$ as an end-point, 
 if $\varphi(s)$ is crossed by an external arc $\mathscr A_t$ and met simply by a second external arc $\mathscr A_{r^\prime}$ for $\varphi(r^\prime)\in A$
such that the short 
directed arc joining $\varphi(s)$ and $\varphi(r^\prime)$ represents an isolated value $\beta$ which is marked and assigned to $\varphi(r^\prime)$. 
Then by rules of markings in case of isolated values,  $\varphi(r^\prime)\in A_+$ and $\varphi(s)=\varphi(r^\prime-1)\in B_+$ necessarily (if $\varphi(s)$ is internal).
By our construction both arcs $\mathscr A_t$ and $\mathscr A_{r'}$ lie on the left of $V$, and $\mathscr A_{r'}$ has an arrow in the same direction as that on the short 
directed arc
 representing $\beta$ and replaces it. Thus the arrow on $\mathscr A_{r^\prime}$  points towards $\varphi(s)$, whilst  the arrow on the external arc $\mathscr A_t$ points away from $\varphi(s)$.
(See Figure \ref{Figure 10}).

 \end {proof}

 \subsection{The Single Edge Property.}\label{6.6}

 The following is immediate from Lemma \ref {6.5}$(iii)$,\,$(iv)$ and $(v)$ above.

\begin {cor}  The directed arcs representing $\{\beta^*_i\}_{i=1}^{e-1}$ have their arrows aligned and form a single edge $E^*$.

\end {cor}

\textbf{Remark}.  Of course this means that $\Pi^*:=\{\beta^*_i\}_{i=1}^{e-1}$ is a simple root system of type $A_{e-1}$.  We call $E^*$ the straightened edge associated to $E$.  As a meander (in the new sense given in \ref{6.2}) it passes through the same points of $E$ but on a different path.

\subsection{The Order Relation.}\label{6.7}

By \ref{6.6}, $\Pi^*=\{\beta_i^*\}_{i=1}^{e-1}$ is a simple root system of type $A_{e-1}$.  By Lemma \ref {6.5}$(i)$, $(ii)$ it satisfies b) of \ref {6.1}.  To show that a) holds it is enough to find a
strict linear order relation $\leq$ on $[1,\,e-1]$ so that the matrix taking $\epsilon_i\beta_i$ to $\beta_i^*$ is triangular with ones on the diagonal.

Set $\mathscr S:=\{i \in 
[1,\,e-1]\mid \beta_i^*-\epsilon_i\beta_i=0\}$.  These of course correspond to the unmarked elements.

The elements of $\mathscr S$ are deemed to be the smallest elements.

It remains to define an order relation on the marked elements.  Since each marked element $m \in [1,\,e-1]$ corresponds to just one external arc $A_m$, this is equivalent to giving an order relation on the external arcs.

Define the support of an external arc joining 
$\varphi(i)$ to $\varphi(j)$ with $j>i$ as the set $[i,\,j-1]$.

Define an order relation on the external arcs by the marked elements contained in their support. We claim that this gives the required order relation.

Consider just the external arcs on the left of $V$.

External arcs corresponding to simple (resp. compound) interval values have disjoint supports.

 The support of an external arc $A_m$ corresponding to a compound interval value either contains the support of an external arc or is disjoint from it.  Moreover $m$ does not lie in the support of any other external arc.

 In other words just concerning external arcs on the left of $V$, the external arcs corresponding to simple interval values correspond to smaller elements and the ones corresponding to compound interval values to larger ones.

 A similar property holds for the external arcs on the right of $V$.

 Finally the marked element corresponding to an external arc $A^\ell$ on the left of $V$ cannot be contained in the support of an external arc $A^r$ on the right of $V$ \textit{and at the same time} the marked element corresponding to $A^r$ be contained in the support of $A^\ell$.  Indeed this could only be if both $A^\ell$ and $A^r$ correspond to compound (and not simple) intervals \textit{and} arising from marked elements of defining a signature of different sign at the turning point to which they are associated (as specified in \ref{6.5}).  However in our construction only negative signature leads to compound values which are not simple. Thus this situation does not arise.

 Our claim is proved.    Combined with the remarks in the first paragraph of this section we obtain the

 \begin {prop}  $\Pi^*$ is a simple root system of type $A_{e-1}$ satisfying a) and b) of \ref {6.1}.
 \end {prop}

 \section{Joining Edges.}\label{8}

 \subsection{The problem.}\label{8.1}

Recall that Proposition \ref {3.2} means that $S$ may be viewed as defining a set of non-trivial edges.  In \ref {3.3} we noted that each non-trivial edge $E$ up to a specified sign change defines a connected set of Dynkin type $A_{|E|-1}$.   A major effort was then expended in the subsequent sections 
(\ref{3}-\ref{6}) to show that these signs could be eliminated nevertheless respecting conditions a) and b) of \ref {6.1}. We call this process a straightening of the edge $E$ and the new edge defined by the set $\{\beta^*_i\}$, a straightened edge. As in the index $1$ case studied in \cite {FJ1}, this means that we can find a nilpotent element in $\mathfrak g^*$ (where here we recall that $\mathfrak g = \mathfrak {sl}_n(\bf k)$) of the same regularity as a Jordan block of size $|E|$ whose restriction to $\mathfrak q$ is exactly the element of $\mathfrak q^*$ defined by $E$.  In the hypothetical case when $|E|=n$, this means the regular element $\eta \in \mathfrak q^*$ of the adapted pair $(h,\eta)$ is the restriction of a regular nilpotent element of $\mathfrak g^*$, so realizing our basic goal.

We remark that by Proposition \ref {3.2} one has $\dim \mathfrak h_{\Lambda}= |S|=|E|-1$ (if there is a unique  non-trivial edge $E$) and so $|E|=n$ implies that $\mathfrak h_{\Lambda}= \mathfrak h$.  By definition of $\mathfrak h_{\Lambda}$, this implies that any semi-invariant in $S(\mathfrak q)$ is invariant which by \cite [7.9]{J0} implies that $\mathfrak q = \mathfrak g =\mathfrak {sl}_n(\bf k)$.

In order to reach our goal outside the trivial case $\mathfrak q = \mathfrak g =\mathfrak {sl}_n(\bf k)$, it is necessary to augment $S$ to a set of $n-1$ elements defining a single edge of cardinality $n$.  There are two possible ways to do this, either before straightening or after.  The first approach was adopted in \cite {FJ1} in which we introduced the notion of an exceptional value. However in general it is not too clear what a corresponding set of exceptional values should be, though one can easily make inspired guesses in some special cases. 
Let $V$ be an $\ad\, h$ stable complement to $(\ad\,\q_{\Lambda})(\eta)$ in $\q_{\Lambda}^*$. It is shown in \cite[2.2.4]{FJ3} that we can take $V=\sum_{\alpha\in T}{\bf k}x_{\alpha}$ for some subset $T$ of $R$. Moreover in \cite[Prop. 2.2.5]{FJ3}, we have shown that $S\cup T$ spans $\h^*$.
 Despite this result, the subset $T$ is too big and in general it is even inappropriate to take a suitable subset of $T$, as for example when $\mathfrak q$ is a Borel subalgebra.  An alternative is simply to define a complement to $S$ by joining the non-trivial edges defined by $S$.  Here one may remark that the end-points of these edges are the fixed points of $\sigma^+$ and $\sigma^-$  
but not of both and are relatively easy to locate.  Nevertheless we failed to find a general rule to join these fixed points so that the straightening procedure could be carried out - a process which is particularly delicate.

 \subsection{Straightened Edge Values.}\label{8.2}

 The second approach to the above problem is to augment $S$ after straightening, that is to add elements to $\Pi^*$ so that it becomes a system of type $A_{n-1}$.  This in itself would be a trivial matter. However we must ensure that these elements do not make their appearance in the expression for $\eta$, equivalently that the corresponding root vectors lie in the kernel of the restriction map $\mathfrak g^* \rightarrow \mathfrak q^*$.  This means simply that the added roots must lie in $K$.  This is the analogue of condition b) of \ref {6.1} which we saw was rather non-trivial to satisfy.

 A key point in satisfying condition b) of \ref {6.1} was to have sufficiently many nil elements.  This is also true in the present analogous situation, so our first step which will occupy several subsections is to find such nil elements.  Here we first examine the ``value" of a straightened edge ultimately showing (Corollary \ref {8.4}) that it belongs to $-K$.

 Let $E$ be a non-trivial edge defined by elements of $S$ and recall \ref {1.3} that we call one of its end-points its starting point $a$ and the second its finishing point $b$. They can be interchanged.  
Recall that the straightened edge $E^*$ associated to $E$ (when it is represented with respect to the horizontal line $L$) passes through the same points of $L$ as the edge $E$ 
where some paths are altered, that is by replacing a marked value $\beta_i$ by a new value $\beta_i^*$ (which corresponds on the vertical central line $V$ to replace a marked short central line by an external arc following the construction of Section \ref{6}).
Recall also (Corollary \ref{6.6}) that the arrows of the directed arcs forming $E^*$ are all aligned, hence (\ref{6.2}) we can define the starting and the finishing point of $E^*$. Let $a^*$ (resp. $b^*$) denote the starting (resp. finishing) point of $E^*$.
We will show in the Proposition below that this source (resp. sink) for $E^*$ is also a source (resp. a sink) for the edge $E$.


 Observe that $a^*,\,b^*$ are integer points on $L$.  We call $\varepsilon_{a^*}-\varepsilon_{b^*}$ the value $\iota_{E^*}$ of the straightened edge $E^*$ or in general a straightened edge value.
(Observe that this terminology is consistent with that used in \ref{1.3} for the value of an edge.)

 Obviously $\iota_{E^*}=\sum_{i=1}^{e-1} \beta^*_i$, however this is not so useful.

 A key point is that $\iota_{E^*}$ can be expressed as an interval value (in the sense of \ref{1.5}), see Lemma \ref{8.3} below.  
 
In the Proposition below we represent the points of $E$ on a vertical central line $V$ as in Section \ref{6}, 
 replacing some short directed arcs on $V$ (those which are marked) by directed external arcs.
In the following the order (defining first, last, ...) is given by passing from the top to the bottom of $V$. 
 Recall (\ref{6.4}) that $A$ (resp. $B$) denotes the type of a turning point of the edge $E$, that is  a source (resp. a sink) of $E$. Recall also the notations of $A_+$, $A_-$, $B_+$ and $B_-$ (\ref{6.5}).
 
 

 \begin {prop} Assume $\mid E\mid=e >1$. Let $E^*$ be its straightened edge, with starting point $a^*$ and finishing point $b^*$.

 \begin{enumerate}

 \item[(i)] Then $a^*,\,b^*$ are turning points of $E$, with $a^*$ of type $A$ and $b^*$ of type $B$.

 \item[(ii)] If $B_-\neq \emptyset$, then $a^*$ is the turning point of type $A$ which lies immediately 
after the last element of $B_-$ . Otherwise $a^*$ is the first turning point in $A$.

\item[ (iii)] If $A_-\neq \emptyset$, then $b^*$ is the turning point of type $B$ which lies immediately after the last 
element of $A_-$. Otherwise $b^*$ is the first turning point in $B$.

 \item[(iv)]  $a^*=a$ if and only if $a\in A$ and $B_-=\emptyset$. 

 \item[(v)]   $a^*=b$ if and only if $b\in A$ and the last internal turning point of type $B$ has negative signature.

 \item[(vi)] $b^*=a$ if and only if $a\in B$ and $A_-=\emptyset$.

\item[ (vii)] $b^*=b$ if and only if $b\in B$ and the last internal turning point of type $A$ has negative signature.
\end{enumerate}

 \end {prop}

\begin{figure}[!h]
\centering
\input{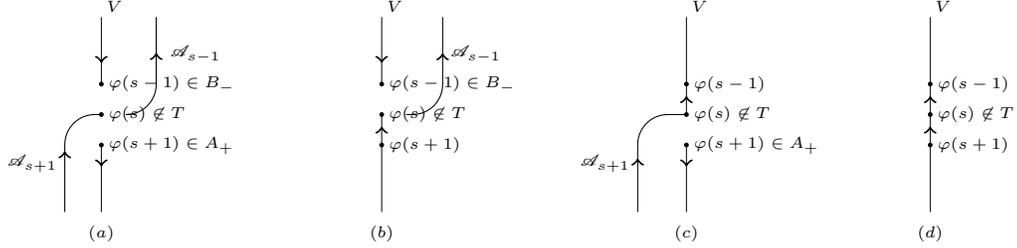}
\caption{\scriptsize{\it 
The case when $\varphi(s)$ is not a turning point of $E$.}
Case $(a)$ : both $\beta_{s-1}$ and $\beta_s$ are marked 
and then two directed external arcs
$\mathscr A_{s-1}$ and $\mathscr A_{s+1}$
replacing the short directed arcs corresponding
to $\beta_{s-1}$ and $\beta_s$
meet (simply) $\varphi(s)$.
Case $(b)$ : only $\beta_{s-1}$ is marked
and then only one directed external arc,
namely $\mathscr A_{s-1}$, meets (simply) $\varphi(s)$
and the short directed arc corresponding to $\beta_s$ is not deleted.
Case $(c)$ : only $\beta_s$ is marked and then only one directed external arc,
namely $\mathscr A_{s+1}$, meets (simply) $\varphi(s)$
and the short directed arc corresponding to $\beta_{s-1}$ is not deleted.
Case $(d)$ : neither $\beta_{s-1}$ nor $\beta_s$ are marked
and then the short directed arcs corresponding to $\beta_{s-1}$ and $\beta_s$
meet $\varphi(s)$ and are not deleted.
In all these cases, $\varphi(s)$ cannot be an end-point of the straightened edge $E^*$.
 }
 \label{Figure 11}
  \end{figure}

\begin{figure}[!h]
\centering
\input{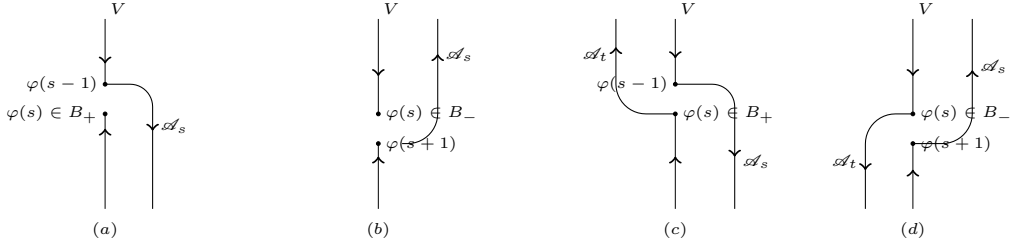}
\caption{\scriptsize{\it 
The case when $\varphi(s)$ is an internal turning point 
 with only one boundary value 
  marked.}
 This boundary value is deleted and 
 replaced by the external arc $\mathscr A_s$.
 Cases $(a)$ and $(b)$ : no external arc crosses $\varphi(s)$ 
 and then $\varphi(s)$ is an end-point of the straightened edge $E^*$.
 Cases $(c)$ and $(d)$ : an external arc $\mathscr A_t$ crosses $\varphi(s)$ 
 and then $\varphi(s)$ is not an end-point of $E^*$.}
\label{Figure 12}
\end{figure}

\begin{figure}[!h]
\centering
\input{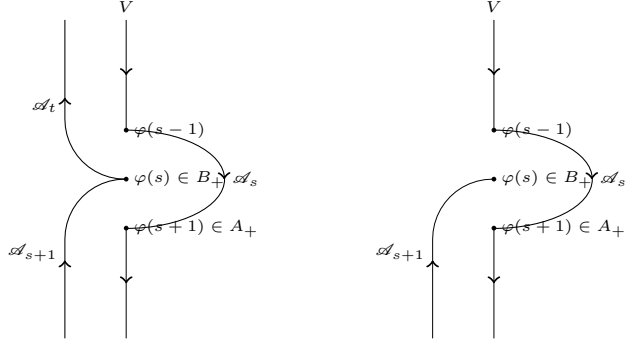}

\caption{\scriptsize{\it 
The case when $\varphi(s)$ is an internal turning point 
 with both boundary values 
  marked.}
It means that $\varphi(s)$ and $\varphi(s+1)$
 are both
  internal turning points
    and then that $\beta_s$ is an isolated value.
This implies that $\varphi(s)\in B_+$ and $\varphi(s+1)\in A_+$.
Then the short directed arcs corresponding 
 to the values $\beta_{s-1}$ and $\beta_s$
 are deleted and replaced resp. 
  by the directed 
  external arcs $\mathscr A_s$  
and $\mathscr A_{s+1}$.  
The Figure on the left (resp. on the right) represents  
the case when $\varphi(s)$ is crossed (resp. is not crossed),
 and then $\varphi(s)$  is not (resp. is) an end-point of $E^*$. }
 \label{Figure 13}
 \end{figure}

 \begin {proof}

 Recall the procedure to assign markings described in \ref {6.4} and \ref{6.5} {\bf(M)}.

  Suppose $\varphi(s)$ is not a turning point of $E$.

 If one neighbouring value of $\varphi(s)$ that is $\beta_{s-1}$ or $\beta_s$ is marked then it is also the neighbour to an internal turning point $\varphi(t)=\varphi(s\mp 1)$ of $E$ and the corresponding 
 short directed arc is replaced by an external arc 
 $\mathscr A_t$ assigned to $\varphi(t)\in T_0$ in the terminology of \ref{6.5} and which meets simply $\varphi(s)$. 
    Consequently $\varphi(s)$ cannot become an end-point of $E^*$. (See Figure \ref{Figure 11}).

  Suppose $\varphi(s)$ is a turning point of $E$.
With the terminology of \ref{6.5},
  we claim in this case that $\varphi(s)$ becomes an end-point of $E^*$ if and only if it is not crossed by a directed external arc. (See Figures \ref{Figure 12} and \ref{Figure 13}).

 Suppose first that there are no isolated points. Then $\varphi(s)$ has exactly one neighbouring value which is left unmarked.  This is unaltered whilst if $\varphi(s)$ is internal its second neighbouring value is marked and the corresponding 
short directed arc suppressed. This proves the claim when there are no isolated points.

  If there are isolated points then in the biparabolic case it can happen that both neighbouring values are marked.  However in this case one checks from our construction (see \ref {6.4}) that there is an external arc 
  $\mathscr A_{r'}$ assigned to $\varphi(r')=\varphi(s+1)$ such that the short directed arc $(\varphi(r'),\,\varphi(s))$ if $\varphi(s)\in B$ (resp. $(\varphi(s),\,\varphi(r'))$ if $\varphi(s)\in A$) represents an isolated value and such that  $\mathscr A_{r'}$ meets simply $\varphi(s)$ and then does not cross $\varphi(s)$ in the terminology of \ref{6.5}.
   Consequently the claim also follows in the presence of isolated points.

  To prove $(i)$ assume that $\varphi(s)$ is an end-point of $E^*$. Then by the above $\varphi(s)\in T$. Moreover we note that either there is exactly one unmarked 
  short directed arc 
meeting $\varphi(s)$ and this is left unchanged, or (in the presence of isolated values) one of the two marked 
  short directed arcs meeting $\varphi(s)$ is replaced by an external arc 
 meeting simply $\varphi(s)$ (terminology of \ref{6.5}) with an arrow pointing in the \textit{same direction}.  Thus if $\varphi(s)$ is a source (resp. sink) of $E$ and becomes an end-point of $E^*$, then it must become the starting (resp. finishing) point $a^*$ (resp. $b^*$) of $E^*$.

  (One may also observe that the number of turning points of type $B$ (resp. $A$) which are 
crossed by a directed external arc is exactly the number of internal turning points of type $A$ (resp. $B$).  Yet the difference between these latter numbers is always one !)

   To prove $(ii)$ we must determine which turning points of type $A$ are 
crossed in the terminology of \ref{6.5} by an external arc on the right assigned to an internal turning point $\varphi(t)$ of type $B$.  Here we note that if the signature at $\varphi(t)$ is negative then the corresponding external arc starts just below $\varphi(t)$, goes upwards, skipping over any immediate internal turning points of type $B$ which have positive signature, till it reaches a turning point of type $A$. Thus every turning point of type $A$ above an internal turning point of type $B$ having negative signature must be 
crossed by an external arc assigned to some internal turning point of type $B$.  On the other hand if at a turning point $\varphi(s)$ of type $A$ the 
internal turning point of type $B$ immediately above $\varphi(s)$ has negative signature whilst all the internal turning points of type $B$ below $\varphi(s)$ (if these ones exist) have positive signature, then $\varphi(s)$ becomes an end-point of $E^*$ since it is not crossed by a directed external arc.  Hence $(ii)$.

   (Strictly speaking the first observation in the paragraph above was unnecessary in view of $(i)$.  However it does show the consistency of our analysis.)

   The proof of $(iii)$ is similar to that of $(ii)$.

   It is clear that $(iv)$ and $(v)$ follow from $(ii)$ and $(vi)$,\,$(vii)$ from $(iii)$.

   \end {proof}

\clearpage

\subsection{Signs.}\label{8.3}

Recall the hypotheses and notations of \ref {8.2}.

In terms of the function $\varphi \,:\,[1,\,e] \mapsto E$, we can write $a^*=\varphi(s),\,b^*=\varphi(t)$, for some $s,\,t \in [1,\,e]$.

Set $I=I_{s,\,t}$ if $s<t$ and $I_{t,\,s}$ if $t<s$.

By Proposition \ref {8.2}$(i)$, $\iota_I$ is an odd compound interval value or a simple one.  Define $\epsilon_I$ as in \ref {4.3} or as in Corollary \ref{3.8.4}, that is $\epsilon_I=\epsilon_s$ if $s<t$ and $\epsilon_I=\epsilon_t$ if $t<s$.  

\begin {lemma}  $\iota_{E^*} = \epsilon_I\iota_I$.
\end {lemma}

\begin {proof}  Suppose 
$s<t$ (resp. $s>t$). Then $\iota_{E^*} = \iota_I$ (resp. $\iota_{E^*} = -\iota_I$). On the other hand 
by Corollary \ref{3.3} since $a^*=\varphi(s)$ is a source then $\epsilon_s=1$ unless $a^*=b$ and since $b^*=\varphi(t)$ is a sink then $\epsilon_t=-1$ unless $b^*=b$. Thus if $s<t$ then $\epsilon_I\iota_I=\epsilon_s\iota_I=\iota_I=\iota_{E^*}$ and if $s>t$ then $\epsilon_I\iota_I=\epsilon_t\iota_I=-\iota_I=\iota_{E^*}$.

\end {proof}

\subsection{Straightened Edge Values lie in $-K$.}\label{8.4}

First we need the following Lemma also useful in next Section (see Lemma \ref{9.3}).

\begin{lemma}

Let $E$ be a non-trivial edge and $E^*$ its associated straightened edge with starting point $a^*$ and finishing point $b^*$. When $E\neq E^*$ it cannot occur that $a^*$ belongs to an equicentral component of the double partition and simultaneously that $b^*$ belongs to an equicentral component of the double partition
 not necessarily equal to the first one.
\end{lemma}

\begin{proof} Assume that $a^*$ (resp. $b^*$) belongs to some equicentral double component $\mathscr I$ (resp. $\mathscr I'$)  (with $\mathscr I$ and $\mathscr I'$ equal or not). Since an equicentral double component does not contain any internal turning point (by
\ref{3.4.5.3}$(vii)$) $a^*$ and $b^*$ must be the end-points $a,\,b$ of $E$. 
Suppose that $a^*=a$ and $b^*=b$.
Recall Proposition \ref {3.8.4}$(iv)$ and denote by $t_a$ (resp. $t_b$) the first turning point of the meander $M_a$ (resp. $M_b$) as the meander leaves $\mathscr I$ (resp. $\mathscr I'$) starting from $a$ (resp. $b$). 
Since $E\neq E^*$, $E$ has at least one internal turning point and then $t_a$ and $t_b$ are internal turning points of $E$.

Set $t_b=\varphi(u)$. Then $1<u<e$ and by \ref{1.9} $\beta_u=\varepsilon_{\varphi(u)}-\varepsilon_{\varphi(u+1)}$ is nil since $\varphi(u)\in\mathscr E'=[1,\,n]\setminus\mathscr I'$ and $\varphi(u+1)\in\mathscr I'$.
  By the rule $(ii)$ of \ref {6.4} supplemented by rule ($\textbf{M}$) of \ref {6.5} the opposite neighbour to $\varphi(u)$ (that is $\beta_{u-1}$) is marked and assigned to $t_b$.  This means that $t_b$ has a positive signature and hence $b\neq b^*$ by Proposition \ref {8.2}$(vii)$.

The other case (that is $a=b^*$ and $b=a^*$) is similar interchanging $\mathscr I$ and $\mathscr I'$ and applying Proposition \ref {8.2}$(v)$.
\end{proof}

\begin {cor} Let $E$ be a non-trivial edge and $E^*$ its associated straightened edge.  Then $\iota_{E^*} \in -K$.  If $E=E^*$, then $\iota_{E^*} \in -K\cap R = -R_*$.
 \end {cor}

 \begin {proof}  If $E=E^*$, that is to say if $E$ is already straightened, then the assertion follows from  Lemma \ref{8.3} and Corollary \ref {3.8.4} since in this case $E$ has no internal turning point.  Otherwise it follows from Lemma \ref{8.3} and 
Lemma \ref{4.3.2} unless $\Phi(I)=E$, equivalently unless that $\{a,\,b\}=\{a^*,\,b^*\}$. 
Notice that this last condition immediately implies, via Proposition \ref {8.2}$(i)$, that $E$ has the same number of turning points of type $A$ as of type $B$.  Moreover in this case $\iota_{E^*}$ is non-nil if and only if $\varepsilon_a-\varepsilon_b \in M$ forcing both $a,\,b$ to lie in an equicentral component $\mathscr I$ of the double partition (by \ref{3.4.5}). 
In the notation of \ref {3.8.4} this means that 
$M_a=M_b(=E)$.  But by Lemma \ref{8.4} above this contradicts both $b =a^*$ and $b=b^*$, when as in the present case $E\neq E^*$.

 
  
When $\{a,\,b\}=\{a^*,\,b^*\}$ and $E\neq E^*$ then we have $\iota_{E^*}\in K\cup -K$. Let us show that $\iota_{E^*}\in -K$.
If $a^*=a$ and $b^*=b$, then $\iota_{E^*}=\varepsilon_a-\varepsilon_b=\epsilon_1(\varepsilon_a-\varepsilon_b)$ since by Corollary \ref{3.3}, $\epsilon_1=1$ because $a=\varphi(1)$ is a source in this case (by Proposition \ref{8.2}$(i)$).
  If $a^*=b$ and $b^*=a$, then $\iota_{E^*}=\varepsilon_b-\varepsilon_a=\epsilon_1(\varepsilon_a-\varepsilon_b)$ since by Corollary \ref{3.3}, $\epsilon_1=-1$ because $a=\varphi(1)$ is a sink in this case (by Proposition \ref{8.2}$(i)$). In both cases $\iota_{E^*}=\epsilon_I\iota_I$ with $I=I_{1,\,e}$ nil
  and the conclusion of the Corollary follows from Lemma \ref{4.3.1}.

  \end {proof}

 \textbf{Remark}.  
Recall the notations and hypotheses of Lemma \ref{8.4} and of its proof.
The neighbour to $t_a$ in the simple interval joining $t_a$ and $a$ (that is $\beta_{v-1}$ if $t_a=\varphi(v)$) is also nil.  However except in the parabolic case we cannot conclude that $t_a$ has negative signature. This is because of the slight asymmetry in rule ($\textbf{M}$).

 

 \subsection{Fully Fixed Points.}\label{8.5}

 An element of $[1,\,n]$ fixed under $<\sigma^+,\,\sigma^->$ is called a fully fixed point (for short).  A fully fixed point can be viewed as defining a 
trivial edge.

 We shall see that the following result also forces the existence of nil elements in completing the sets $\{\beta_i^*\}$.

 Let $\mathscr I:=J^+\cap J^-$ be a component of the double partition $\mathscr J^\pm$.  Let $c^+$ (resp. $c^-$) be the centre of $J^+$ (resp. $J^-$).  By symmetry we can assume without loss of generality that $c^+\leq c^-$.

 \begin {lemma} Let $f \in \mathscr I$ be a fully fixed point. Then $c^+<c^-$ and $f \in [c^+,\,c^-]$.  Moreover $\mathscr I$ has at most one fully fixed point.
 \end {lemma}

 \begin {proof}

 The first assertion follows from Proposition \ref {3.8.4}$(i)$
if $\mid \mathscr I\mid>1$. Moreover if $\mid \mathscr I\mid=1$ and if $c^+=c^-$ then only $J^+$ (resp. $J^-$) is reduced to one point $c=c^\pm$ which then is fixed under $\sigma^+$ (resp. $\sigma^-$). But this point $c$ cannot be fixed under $\sigma^-$ (resp. $\sigma^+$) by \ref{2.5} since the half-integer arc below (resp. above) $L$ joining $c\pm\frac{1}{2}$ is necessarily marked because the half-integer meander joining these two points is a loop passing only through $c\pm\frac{1}{2}$ and having a fictitious arc which is the one above (resp. below) $L$ joining these two points (see \ref{3.1}). Note that $\mid J^+\mid=\mid J^-\mid=1$ (and $\mathscr I\neq\emptyset$) cannot occur since $\pi^+\cup\pi^-=\pi$.
 

 
The second assertion follows from the Remark after Lemma \ref{3.8.1}.

 The restriction of $\sigma^+$ (resp. $\sigma^-$) to $\mathscr I$ has $\leq 2$ fixed points with equality only if $|J^+|$ (resp. $|J^-|$) is even. Thus $\mathscr I$ has at most two fully fixed points and this only when $c^+,\,c^-$ are half-integer.  Then by \ref {2.5} the two fixed points of $\sigma^+$ (resp. $\sigma^-$) lie on opposite sides of $c^+$ (resp. $c^-$). Thus the required conclusion results from the first part.
 \end {proof}

 \subsection{}\label{8.6}

 Now consider a straightened edge $E^*$ as a single arc joining the points $a^*,\,b^*$ of the horizontal line $L$.  Suppose $a^* < b^*$.  Then $\iota_{E^*}$ is the positive root $\varepsilon_{a^*}- \varepsilon_{b^*}$ and by Corollary \ref {8.4} an element of $-K$.  Hence by Lemma \ref {1.8}$(i)$,\,$(ii)$ it belongs to $-K^-$.  This means that the arc representing $E^*$ 
 attains a gap defined by $\hat{\pi}^-$ and we draw it above $L$.  Conversely if $a^* > b^*$, then this arc 
attains a gap defined by $\hat{\pi}^+$ and we draw it below $L$.

 In the parabolic case one has $\hat{\pi}^+=\emptyset$, so all such arcs lie above $L$. However this fails in the biparabolic case.  Moreover there is no reason why an arc corresponding to a straightened edge $E^*$ should not 
attain both a gap defined by $\hat{\pi}^+$ and a gap defined by $\hat{\pi}^-$.  However by 
Corollary \ref{8.4} there is an important exception, namely when $E=E^*$.  This will be used in the proof of Lemma \ref {9.3}.

 A non-trivial edge $E$ (coming from $S$) may be viewed as a sequence of arcs passing through a set $p_E$ of integer points of $L$. Observe that its straightening $E^*$ is also a sequence of arcs and passes through the same set $p_E$ of integer points of $L$ though in a different order (and by different arcs).  Let $p_S$ denote the (necessarily) disjoint union of the $p_E$ when $S$ is expressed as a disjoint union of non-trivial edges $E$. The following is immediate.

 \begin {lemma} $[1,\,n]\setminus p_S$ is just the set of fully fixed points.
 \end {lemma}

\subsection{}\label{8.7}

We may also deduce the last part of Lemma \ref {8.5} by a more elegant argument using Proposition \ref {3.2}.  This even gives a little more.

\begin {lemma}  Let $f_1,\,f_2,\,\ldots,\,f_k$ be distinct fully fixed points.  For all $i=1,\,2,\,\ldots,\, k-1$ let $\gamma_i$ be the arc joining $f_i$ and $f_{i+1}$ viewed as a root.  Then the $\gamma_i\,:\,i=1,\,2,\,\ldots,\,k-1$ are linearly independent and any sum of these elements lying in $\Delta$ lies in $K\cup -K$.
\end {lemma}

\begin {proof}  Let $S^\perp$ denote the orthogonal of $S$ in $\mathfrak h$.  By Proposition \ref {3.2} one has $S^\perp \cap \mathfrak h_{\Lambda} =\{0\}$ and $S^\perp+\mathfrak h_{\Lambda}=\mathfrak h$.
Since fully fixed points do not lie in $p_S$ we conclude that $\gamma_i \in S^\perp$, for all $i=1,\,2,\,\ldots,\, k-1$ under the identification of $\mathfrak h^*$ with $\mathfrak h$ through the Killing form.  Under 
the same identification, we have $M \subset \mathfrak h_{\Lambda}$, since $M$ lies in the semisimple part of the Levi factor of $\mathfrak q$.  Then the assertion follows through the above and Lemma \ref {1.8}$(iv)$.
\end {proof}

\section{Main Theorem.}\label{9}

\subsection{Ultimate Goal.}\label{9.1}

In \cite {J1} we described at least one and in general several adapted pairs for every truncated biparabolic subalgebra $\mathfrak q_{\Lambda}$ of a simple Lie algebra $\mathfrak g$ of type $A$.  We do not know if all such pairs are obtained (up to equivalence) or what are the equivalences between the given pairs  (except if $\ell(\mathfrak q_{\Lambda})=1$, (see \ref {00.2}) in which case $\mathscr N(\q_{\Lambda})$ is irreducible).

In \cite {FJ1} we showed that if $\ell (\mathfrak q_{\Lambda})=1$, then for the unique, up to equivalence, adapted pair $(h,\,\eta)$ one can express $\eta$ as the image of a \textit{regular nilpotent} element of $\mathfrak g^*$, though not in general uniquely.  This means in particular that $\eta$ can be associated to an element $w$ of the Weyl group $W$.  Of course $w$ will not be uniquely determined and in any case can be modified by any element of the Weyl group of the Levi factor of $\mathfrak q$.  Ideally $w$ will not even depend on the choice of the biparabolic subalgebra.

Fix a double partition $\mathscr J^\pm$ of $[1,\,n]$ and let $\mathfrak q = \mathfrak q_{\mathscr J^-,\,\mathscr J^+}$ be the corresponding biparabolic subalgebra of $\mathfrak g = \mathfrak {sl}_n(\bf k)$ defined as in \ref {1.8}.  Let $\mathfrak q_{\Lambda}$ be the (canonical) truncation of $\mathfrak q$ obtained by replacing the (standard) Cartan subalgebra $\mathfrak h$ of $\mathfrak q$ by its canonical truncation $\mathfrak h_{\Lambda}$ as defined in \ref {3.2}.

Let $\Xi$ be the set of adapted pairs for $\q_{\Lambda}$ constructed in \cite{J1} (recall also \ref{1.8} and in particular that such an adapted pair in $\Xi$ has its first element in $\h_{\Lambda}$).
Our goal is to prove the following

\begin {thm}  Let $(h,\,\eta) \in \Xi$.  Then $\eta$ is the restriction to $\q_{\Lambda}$ of a regular nilpotent element $y\in\mathfrak g^*$.

\end {thm}

\begin{proof}

For every non-trivial edge $E$ (of cardinality $e$) obtained from $S$, recall (\ref{3.3}) that $S_E=\{\epsilon_i\beta_i\}_{1\le i\le e-1}$ is the set of the roots corresponding to the non-trivial directed arcs forming $E$ and let $\Pi_E^*=\{\beta_i^*\}_{1\le i\le e-1}$ be the simple root system constructed in Section \ref{6}. Set $I_E=I_{1,\,e}=[1,\,e-1]$ and 
$$y_E=\sum_{i\in I_E}x_{\epsilon_i\beta_i}+\sum_{i\in I_E\mid\beta_i^*\neq\epsilon_i\beta_i}x_{\beta_i^*}.$$
Then by \ref{6.1}(b) the restriction of $y_E$ to $\q_{\Lambda}$ is just $$\eta_E=\sum_{i\in I_E}x_{\epsilon_i\beta_i}$$
so then $\eta$ is recovered by summing over all non-trivial edges (see also \ref{3.3}).
On the other hand $y_E$ by \ref{6.1}(a) is up to conjugation a Jordan block of size $e$. In order to obtain a regular nilpotent element of $\g^*$ we must add to the sum of the $y_E$, as $E$ runs over the non-trivial edges defined by $S$, a set of elements which join these Jordan blocks to form a Jordan block of size $n$ and whose weights lie in $K$ (so then vanish under restriction). This is achieved in subsections \ref{9.2} -- \ref{9.12}.

\end{proof}

\textbf{Remark}.  The hypothesis on $\eta$ means that it is given by \ref {1.8}$(*)$, where $S$ is constructed as in \cite [Sect. 5]{J1}. By \cite [Sect. 8]{J1} $\eta$ is regular in $\mathfrak q_{\Lambda}^*$.  In our present analysis we consider a somewhat more general choice of $S$ using anti-Toeplitz involutions (\ref {2.1}).  Sometimes the resulting $\eta$ can fail to be regular (for example if $\mathfrak q$ is the Borel subalgebra of $\mathfrak g$) and sometimes it is always regular (for example if $\mathfrak q =\mathfrak g$).  Unfortunately it is difficult to tell when elements of this more general class are regular and even more difficult to tell if they lead to new inequivalent adapted pairs.

\subsection{}\label{9.2}

Recall \ref {8.4} that a straightened edge value lies in $-K$.  It is either a positive (resp. negative) root and so lies in $-K^-$ (resp. $-K^+$).  Let $\mathscr F^+$ (resp. $\mathscr F^-$) denote the
 set of straightened edges whose value lies in $-K^-$ (resp. $-K^+$).  We may write
$\mathscr F^+ =\{E^*_1,\,\ldots,\,E^*_r\},\,\mathscr F^- =\{E^{\prime *}_1,\,\ldots,\,E^{\prime *}_s\}$ so that their starting points are strictly increasing.  We represent $E^*_i$ (resp. $E_j^{\prime *}$) by an arc $a_i^*\rightarrow b_i^*$ (resp. $b_j^{\prime *} \leftarrow a_j^{\prime *}$) drawn above (resp. below) $L$.

\begin {lemma}  For all $i=1,\,2,\,\ldots,\,r-1$ (resp. $j=1,\,2,\,\ldots,\,s-1$), $\varepsilon_{b_{i+1}^*}-\varepsilon_{a_i^*}$ (resp. $\varepsilon_{b_j^{\prime *}}-\varepsilon_{a_{j+1}^{\prime *}}$) is a root lying in $K$.
\end {lemma}

\begin {proof} That these are roots is obvious. Again $\varepsilon_{b_{i+1}^*}-\varepsilon_{a_{i+1}^*}$ is a negative root lying in $K$, hence in $K^-$. By construction $\varepsilon_{a_{i+1}^*}-\varepsilon_{a_i^*} \in \Delta^-$.  Since $(K^-+\Delta^-)\cap \Delta \subset K^-$, the first part follows. The proof of the second part is similar.
\end {proof}

\subsection{}\label{9.3}

It is clear that the above lemma (together with Proposition \ref {6.7})
establishes Theorem \ref {9.1} in the parabolic case (in which case $\mathscr F^- =\emptyset$) when there are no fully fixed points.  Here the additional roots we must add to the union of the sets $\{\beta^*_i\}$ described in Section \ref{6} giving the straightened edges are exactly the $\gamma_i:=\varepsilon_{b_{i+1}^*}-\varepsilon_{a_i^*}$.

In the biparabolic case we need a further result to establish our main theorem.

Recall (\ref{1.9}) that $u,\,v \in [1,\,n]$ lie in the same component of the double partition $\mathscr J^\pm$ if and only if $\varepsilon_u-\varepsilon_v \in M$.

By \ref {3.4.5}, outside the equicentral case, all the turning points in  a double component $\mathscr I$ are sinks (resp. sources) if $c^+<c^-$ (resp. $c^+>c^-$).
Then this fact 
combined with Prop. \ref{8.2}$(i)$ implies that two end-points (one being a starting point and the other a finishing point) of straightened edges (not necessarily equal) can only lie in the same double component  if it is equicentral.

Moreover, by \ref{3.4.5.3}$(vii)$, in an equicentral  double component, all the turning points are end-points (of the unstraightened edges) and their number is at most two. Thus an equicentral double component cannot contain more than two end-points of straightened edges.\par


\begin {lemma}  
Let $E^*$ and $E^{\prime*}$ be distinct straightened edges with respective starting point $a^*,\,a^{\prime*}$ and respective finishing point $b^*,\,b^{\prime*}$. It cannot happen that  both elements of $\{a^*,\,b^{\prime*}\}$ lie in the same component 
$\mathscr I$ of the double partition, and that simultaneously 
both elements of $\{a^{\prime*},\,b^*\}$ lie in the same component 
$\mathscr I'$ of the double partition.
\end {lemma}

\begin{figure}[!h]
\centering
\input{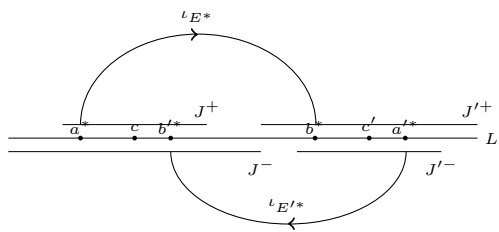}
\caption{\scriptsize{\it 
Illustration of the proof of Lemma \ref{9.3}.}
If the points $a^*$ and $b'^*$ (resp. $a'^*$ and $b^*$)
belong  to the
  same equicentral double component 
 then
  both double components
   are necessarily not linked.
The arc in $K^+$ (resp. in $K^-$) joining $b^{\prime*}$
 to $a^{\prime*}$
  (resp. $b^{*}$ to $a^{*}$)
  must attain a gap
      defined both by $\hat{\pi}^+$
       and by $\hat{\pi}^-$
       then cannot belong to $R^+_*$ (resp. $R^-_*$).}
\label{Figure 15}
\end{figure}

\begin {proof}  Suppose 
the contrary. The above observations imply that  
$\mathscr I$ and $\mathscr I'$ are equicentral and distinct (and then not linked).  Let 
$E$, resp. $E'$, be the edge corresponding to 
$E^*$, resp. $E^{\prime*}$. 

Yet by 
Lemma \ref {8.4} 
we conclude that 
$E=E^*$ and $E'=E^{\prime*}$.  However in this case both straightened edge values lie in $-R_*=R\cap -K$ by Corollary \ref{8.4}.  Again the hypothesis implies that 
$\iota_{E^*}\in \Delta^+$ and $\iota_{E^{\prime*}}\in \Delta^-$ or vice-versa.   Both cases are similar and we consider just the first. Then 
$\iota_{E^*}\in -R^-_*$ and so must 
attain a gap defined by $\hat{\pi}^-$, but cannot 
attain a gap defined by $\hat{\pi}^+$, whilst 
$\iota_{E^{\prime*}} \in -R^+_*$ and so must 
attain a gap defined by $\hat{\pi}^+$, but cannot 
attain a gap defined by $\hat{\pi}^-$.  Thus we cannot have both 
$\varepsilon_{a^*}-\varepsilon_{b^{\prime*}}\in M$ and $\varepsilon_{a^{\prime*}}-\varepsilon_{b^*}\in M$ by Lemma \ref{1.8}$(iii)$,\,$(ix)$ and $(x)$ and this contradiction proves the lemma.
(See Figure \ref{Figure 15}).
\end {proof}

\subsection{}\label{9.4}

In order to prove our main theorem in the biparabolic case, we must join (with notations of \ref{9.2}) the sequence of arcs

$$a_r^* \rightarrow b^*_r \rightarrow a_{r-1}^* \rightarrow \ldots \rightarrow a_1^*\rightarrow b^*_1,$$
written briefly as $a_r^* \rightarrow b_1^*$, to the sequence of arcs

$$b^{\prime *}_s\leftarrow a_s^{\prime *}\leftarrow b^{\prime *}_{s-1}\leftarrow a_{s-1}^{\prime *}\leftarrow\ldots\leftarrow a_2^{\prime *}\leftarrow b_1^{\prime *}\leftarrow a_1^{\prime *}$$
written briefly as 
$b_s^{\prime *}\leftarrow a_1^{\prime *}$.

In principle this is effected by adding either the arc $b^*_1 \rightarrow a_1^{\prime *}$ \textit{or} the arc $b_s^{\prime *} \rightarrow a_r^*$.  However we need know that at least one of the (non-zero) roots they represent lie in $K$.  It can happen that neither lie in $K$.  This question is studied in the next subsections.

\subsection{}\label{9.5}

 Recall the notations of \ref {9.2}.  

\begin {lemma}  Suppose that the elements of $\{b^*_1, \,a_1^{\prime *}\}$ (resp. $\{b_s^{\prime *},\,a_r^*\}$) do not lie in the same component of the double
partition.   Then $\varepsilon_{b^*_1}-\varepsilon_{a^{\prime *}_1}$ (resp.
 $\varepsilon_{b_s^{\prime *}}-\varepsilon_{a_r^*})$ lies in $K$.
\end {lemma}

\begin{figure}[!h]
\centering
\input{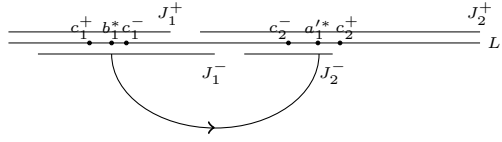}

\caption{\scriptsize{\it
Illustration of the proof of Lemma \ref{9.5}
when $J_1^+\neq J_2^+$ and $J_1^-\neq J_2^-$
 and $b_1^*<a^{\prime *}_1$. }
The arc joining $b_1^*$ to $a^{\prime *}_1$ is drawn
 from left to right and attains a gap defined by $\hat\pi^+$,
  hence is drawn below $L$ and corresponds to a root belonging to $K^+$. 
  Observe that in this case this arc also attains a gap defined by $\hat\pi^-$,
   hence could also be drawn above $L$ and corresponds to a root in $-K^-$.}
\label{Figure 16}
\end{figure}

\begin {proof}  Both cases are similar and we consider just the first.  The hypothesis implies that the root $\varepsilon_{b^*_1}-\varepsilon_{a^{\prime *}_1}$ does not lie in $M$ (by \ref{1.9}) and hence lies in $K\cup -K$.  Suppose $b^*_1$ lies in some component $J^+_1 \cap J^-_1$ and $a_1^{\prime *}$ in some component $J^+_2 \cap J^-_2$. 
Let $c_i^\pm$ be the centre of $J_i^\pm$\,:\, $i=1,\,2$.
Then by \ref {3.4.5} and Prop. \ref{8.2}$(i)$ taking account of sources and sinks one has $c_1^+\le c_1^-$ 
and $c_2^+\ge c_2^-$.  If both $J^+_1\neq J^+_2$ and $J^-_1\neq J^-_2$, then trivially the above root lies in $K\cap -K$ and we are done.  Thus we can assume that $J^+_1= J^+_2$ and $J^-_1\neq J^-_2$ or vice-versa.  Consider just the second case.  Then the above inequalities imply that $J^+_1$ lies strictly to the left of $J^+_2$ since in this case we must have $c_1^+\le c_1^-=c_2^-\le c_2^+$ with at most one equality on the right or on the left.  Consequently $\varepsilon_{b^*_1}-\varepsilon_{a^{\prime *}_1}$ corresponds to an arrow 
 from left to right 
attaining a gap defined by $\hat{\pi}^+$ 
hence drawn below $L$ so defining an element of $K^+$,
as required.  The remaining case is similar. (See Figure \ref{Figure 16}).
\end {proof}

\subsection{}\label{9.6}


It remains to consider the case when both elements of $\{b^*_1,\, a_1^{\prime *}\}$ lie in some equicentral double component $\mathscr I'_1$ 
 and both elements of $\{b_s^{\prime *},\,a_r^*\}$ lie in some equicentral double component $\mathscr I_r$ (different from $\mathscr I'_1$).  By Lemma \ref {9.3} this forces either $r>1$ or $s>1$.  We consider just the first case as the second is similar.  Notice we cannot delete the arc $a^*_r\rightarrow b^*_r$ because it represents several elements of the required simple root system $\Pi^*$.   Yet we may delete the arc $b^*_r \rightarrow a^*_{r-1}$ since it represents a single element of the required root system $\Pi^*$. 
 (This is why we need $r>1$.)  Thus we delete this arc and create new arcs $b^*_r \rightarrow a^{\prime *}_1,\, b_s^{\prime *} \rightarrow a^*_{r-1}$.

\begin {lemma}  The above two arcs both represent elements of $K$.
\end {lemma}

\begin{figure}[!h]
\centering
\input{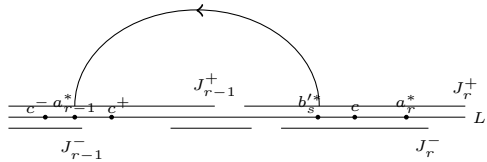}
\caption{\scriptsize{\it 
Illustration of the proof of Lemma \ref{9.6}
for the arc $ b_s^{\prime *} \rightarrow a^*_{r-1}$.}
Here the centre of $\mathscr I_r$ is denoted by $c$, whilst the centre of $J_{r-1}^+$, resp. $J_{r-1}^-$, is denoted by $c^+$, resp. $c^-$.
The arc joining $b^{\prime *}_s$ to $a_{r-1}^*$
  attains at least one gap defined by $\hat\pi^-$
   (and hence is drawn above $L$)
    so corresponds to a root in $K^-$.
     (An analogue of this Figure is also valid if $J_{r-1}^+=J_r^+$.)}
     \label{Figure 18}
     \end{figure}

\begin {proof}  
 By hypothesis $b_s^{\prime *}$ and $a^*_r$ lie in the same equicentral double component $\mathscr I_r:=J_r^+\cap J_r^-$. 
Then by the observations in the beginning of \ref{9.3}, $a^*_{r-1}$ must lie in a different double component $\mathscr I_{r-1}:=J^+_{r-1}\cap J^-_{r-1}$.  
We can now deduce the Lemma from Lemma \ref{9.5} above (with $b_s^{\prime *}$ and $a_{r-1}^*$ instead of $b_s^{\prime *}$ and $a_r^*$).
In more details since $a^*_{r-1} < a^*_r$,  $\mathscr I_{r-1}$ must lie to the left of $\mathscr I_r$.  Again since $a^*_{r-1}$ is a source it follows from \ref {3.4.5} that $J^-_{r-1}\neq J^-_r$. This implies that the arc $ b_s^{\prime *} \rightarrow a^*_{r-1}$ which passes from right to left attains a gap defined by $\hat{\pi}^-$ (and hence is drawn above $L$) and so represents an element of $K^-$, as required.
The proof for the arc $b^*_r \rightarrow a^{\prime *}_1$ belonging to $K$ is similar except that $b_r^*$ can lie to the left or to the right of $a^{\prime *}_1$. In the first case this arc attains a gap defined by $\hat\pi^+$ and then corresponds to a root in $K^+$. In the second case this arc attains a gap defined by $\hat\pi^-$ and corresponds to a root in $K^-$. (See Figure \ref{Figure 18}).
\end {proof}

\subsection{}\label{9.7}

As in \ref {9.3}, combining Lemmas \ref {9.2}, \ref {9.5} and \ref {9.6} proves Theorem \ref {9.1} in the biparabolic case when there are no fully fixed points.

\subsection{}\label{9.8}

Recall (see \ref{3.4.5.4}) the terminology of linked components of the double partition and that a component of the double partition contains at most one fully fixed point (Lemma \ref{8.5}).
Recall also (see \ref{3.4.5.4}) that if two fully fixed points $f_1$ and $f_2$ belong to components $\mathscr I_1$, resp. $\mathscr I_2$, which are not linked, then $\varepsilon_{f_1}-\varepsilon_{f_2}\in K\cap -K$.

The components of the double partition are ordered by the natural ordering in $[1,\,n]$. 
Starting at one component $\mathscr I_r$ and ending at a further component $\mathscr I_t$ to its right, we take $f_1 <f_2<\ldots <f_s$ to be a subset of fully fixed points lying in the components between (and including) $\mathscr I_r$ and $\mathscr I_t$.  

\begin {lemma} If $s\ge 2$ there exists a permutation $p$ of the set $[1,\,s]$ such that $\varepsilon_{f_{p(i)}}-\varepsilon_{f_{p(i+1)}}\in K$, for all  $i \in [1,\,s-1]$.
\end {lemma}

\begin {proof}  The proof is by induction on the number of fully fixed points starting from the left
the case $s=2$ being true since two fully fixed points $f_1$ and $f_2$ cannot belong to the same component of the double partition by Lemma \ref{8.5} and then $\varepsilon_{f_1}-\varepsilon_{f_2}\in K\cup -K$.  Thus assume the assertion holds for all 
$i\in[1,\,u-1]$ for some positive integer 
$u<s$ such that $u\ge 2$. Let $i_j$ be the image of $j$ under the above defined permutation of 
$[1,\,u]$.   By construction 
$f_{u+1}$ is strictly greater than both $f_{i_1}$ and 
$f_{i_u}$.  If 
$\varepsilon_{f_{i_u}}-\varepsilon_{f_{u+1}}\in K$, we are done.  Otherwise the double components $\mathscr I_2,\,\mathscr I_3$ containing the two fully fixed points $f_{i_u}$ and $f_{u+1}$ must be linked. Moreover we can write $\mathscr I_2 = J^+\cap J^-_2,\, \mathscr I_3 =J^+\cap J^-_3$ (resp. $\mathscr I_2=J_2^+\cap J^-$ and $\mathscr I_3=J_3^+\cap J^-$). Then $f_{i_1} \notin J^+$ (resp. $f_{i_1}\not\in J^-$) since $J^+$ (resp. $J^-$) has at most two $\sigma^+$ (resp. $\sigma^-$) fixed points.  This forces 
$f_{i_1}<f_{i_u}$ and then that the double components containing $f_{i_1}$ and 
$f_{u+1}$ are not linked.  Then 
$\varepsilon_{f_{u+1}}-\varepsilon_{f_{i_1}}\in K$ and so the induction proceeds.
\end {proof}

\textbf{Remark}.   Underlying this result is the fact that in the parabolic case there can be at most two fully fixed points 
since there are at most two $\sigma^+$ fixed points in $J^+=[1,\,n]$ in this case.

\textbf{Terminology}. For a permutation $(i_1,\ldots,\,i_s)$ as above, we shall call $f_{i_1}\rightarrow\ldots \rightarrow f_{i_s}$ the ordered chain $ C$  of fully fixed points defined by the set $\{f_i\}_{i=1}^s$ of fully fixed points. Let $a_C,\, b_C$ be the end-points of $C$ with the first we will call a source and the second a sink.

\subsection{}\label{9.9}

Recall the notation of \ref {9.4}.  By Lemmas \ref{9.5} and \ref {9.6} we may join the chains described in \ref {9.4} to form a single chain starting at a source $a^*$ and ending at a sink $b^*$ which links the straightened edges defined by $S$ through elements of $K$.

To prove our main theorem it remains to adjoin the set $F$ of fully fixed points to the ends of this chain by adding links which as roots lie in $K$. This will proceed as follows.

Decompose $F$ (not necessarily uniquely) as a disjoint union $F_{a^*}\sqcup F_{b^*}$, so that
$\varepsilon_f-\varepsilon_{a^*} \in K$ for all $f \in F_{a^*}$ and $\varepsilon_{b^*}-\varepsilon_f \in K$ for all $f \in F_{b^*}$.

Let $C_{a^*}$ (resp. $C_{b^*}$) be the ordered chain defined by $F_{a^*}$ (resp. $F_{b^*}$).  Obviously $b_{C_{a^*}} \in F_{a^*}$ and so $\varepsilon_{b_{C_{a^*}}}-\varepsilon_{a^*} \in K$.  Consequently the chain $C_{a^*}$ can be adjoined to $a^* \rightarrow b^*$ at its end-point $a^*$,
to form a chain $a_{C_{a^*}}\rightarrow b^*$ with all arrows aligned, incorporating $F_{a^*}$ in the chain $a^*\rightarrow b^*$, and this by adding only elements in $K$.  Similarly $\varepsilon_{b^*}-\varepsilon_{a_{C_{b^*}}} \in K$ and so the chain $C_{b^*}$ can be adjoined to
$a_{C_{a^*}}\rightarrow b^*$ at its end-point $b^*$, to obtain a chain $a_{C_{a^*}}\rightarrow b_{C_{b^*}}$ with all arrows aligned, incorporating $F$ in the chain $a^*\rightarrow b^*$, and this by adding only elements in $K$.

It is clear that this will complete the proof of Theorem \ref {9.1}.

To construct the required decomposition of $F$ we consider three possible cases treated below.

\subsection{}\label{9.10}

Suppose $a^*,\,b^*$ lie in the same component $\mathscr I$ of the double partition.

  Let $F^\ell$ (resp. $F^r$) be the set of fully fixed points in double components strictly to the left (resp. right) of $\mathscr I$.

 The hypothesis implies that $\mathscr I$ is an equicentral double component and hence by Lemma \ref {8.5} cannot contain any fully fixed point.  Hence $F=F^\ell \cup F^r$.

We can write $\mathscr I=J^+\cap J^-$ where the centres of $J^\pm$ coincide.  By our standing hypothesis that $\pi^+\cup \pi^- =\pi$ and $\pi^+\cap \pi^- \varsubsetneq\pi$, it follows that $J^\pm$ cannot coincide. Then either $J^-\varsupsetneq J^+$ or $J^-\varsubsetneq J^+$.

In the first case we take $F_{a^*} = F^\ell,\, F_{b^*}=F^r$. In the second case we make the reverse choice.  Then for example given $f \in F^\ell$, the arc $f\rightarrow a^*$ is a positive root and in the first case 
attains at least the gap defined by $J^-\setminus J^+$ and hence represents an element of $K^+$.  The remaining (three) required verifications are equally straightforward.

%
%
%
%

\subsection{}\label{9.11}

Suppose in this and the next subsection that $a^*,\,b^*$ lie in distinct components $\mathscr I_{a^*},\,\mathscr I_{b^*}$ of the double partition.
Write $\mathscr I_{a^*}=J_{a^*}^+\cap J_{a^*}^-$ and $\mathscr I_{b^*}=J_{b^*}^+\cap J_{b^*}^-$.

Suppose that $\mathscr I_{b^*}$ and $\mathscr I_{a^*}$ are not linked.

 There are 
two possibilities which we examine individually.  
 We give the detailed verification of our assertion only in the first case.

($i^\ell$) $\mathscr I_{b^*}$ lies to the left of $\mathscr I_{a^*}$. 

In this case we take $F_{a^*}$ to be the set of fully fixed points lying in a component $\mathscr I=J^+\cap J^-$ of the double partition such that $J^+$ lies to the left of $J^+_{b^*}$ (including $J^+_{b^*}$).  Then each of the roots 
$\varepsilon_f-\varepsilon_{a^*}$ : $f\in F_{a^*}$ is positive and 
attains the gap corresponding to $\hat{\pi}^+$ lying between $J^+$ and $J^+_{a^*}$, hence belongs to $K^+$.\par
Then we take $F_{b^*}$ to be the set of fully fixed points lying in a component $\mathscr I=J^+\cap J^-$ of the double partition such that $J^+$ lies strictly to the right of $J^+_{b^*}$.



($i^r$) $\mathscr I_{b^*}$ lies to the right of $\mathscr I_{a^*}$. 

In this case we take $F_{a^*}$ to be the set of fully fixed points lying in a component $\mathscr I=J^+\cap J^-$ of the double partition such that $J^-$ lies strictly to the right of $J^-_{a^*}$.\par
Then we take $F_{b^*}$ to be the set of fully fixed points lying in a component $\mathscr I=J^+\cap J^-$ of the double partition such that $J^-$ lies to the left of $J^-_{a^*}$ (including $J^-_{a^*}$).



\subsection{}\label{9.12}

Suppose that $\mathscr I_{b^*}$ and $\mathscr I_{a^*}$ are linked.

 There are two possibilities which we examine individually.  
 Verifications of the required properties are straightforward.

($i_+$).  Suppose that $J_{b^*}^-=J_{a^*}^-$.

Since $b^*$ is a sink and $a^*$ is a source this forces $\mathscr I_{b^*}$ to lie to the left of $\mathscr I_{a^*}$.  Then we take  $F_{b^*}$ to be the set of fully fixed points lying in a component $\mathscr I=J^+\cap J^-$ of the double partition such that $J^+$ lies strictly to the right of $J^+_{b^*}$
and $F_{a^*}$ to be the set of fully fixed points lying in a component $\mathscr I=J^+\cap J^-$ of the double partition such that $J^+$ lies to the left of $J^+_{b^*}$ (including $J^+_{b^*}$).

($i_-$).  Suppose that $J_{b^*}^+=J_{a^*}^+$.

Since $b^*$ is a sink and $a^*$ is a source this forces $\mathscr I_{b^*}$ to lie to the right of $\mathscr I_{a^*}$.  Then we take  $F_{a^*}$ to be the set of fully fixed points lying in a component $\mathscr I=J^+\cap J^-$ of the double partition such that $J^-$ lies strictly to the right of 
$J^-_{a^*}$ and $F_{b^*}$ to be the set of fully fixed points lying in a component $\mathscr I=J^+\cap J^-$ of the double partition such that $J^-$ lies to the left of $J^-_{a^*}$ (including $J^-_{a^*}$).

\subsection{}\label{9.13}

Combining the results in the last three subsections and the remarks in the end of \ref {9.9} completes the proof of Theorem \ref{9.1}.

\subsection{}\label{9.14}

It can happen that $S$ is the empty set.  Then $\eta$ (as defined by \ref {1.8}$(*)$) equals zero.  As $\eta$ is supposed to be regular in $\mathfrak q^*_{\Lambda}$ this can only occur if $\mathfrak q_{\Lambda}$ is commutative. One can easily check from our construction that this is indeed the case.  It follows that the elements of $\pi^+$ (resp. of $\pi^-$) are pairwise orthogonal (and so the resulting decomposition of $\pi$ is unique up to permutation).

Take for example $\pi^+=\{\alpha_1,\,\alpha_3\}, \,\pi^- =\{\alpha_2\}$.

 Then the set of fully fixed points here is $\{1,\,2,\,3,\,4\}$ and one possible choice of a simple root system $\Pi^*$ such that $\Pi^*\subset K$ is $\Pi^*=\{-\alpha_1-\alpha_2,\,\alpha_1+\alpha_2+\alpha_3,\,-\alpha_2-\alpha_3\}$. 
 
 The element of the Weyl group $W$ taking 
$\pi$ to $\Pi^*$ is $s_2 s_3 s_1$, which can be thought of as the square root of the unique longest element of $W$
(in the sense discussed in \cite [preamble and Sect. 2]{J2}).

In \cite [Sect. 3]{J2} a somewhat analogous result was obtained when $\mathfrak q$ is the Borel subalgebra of $\mathfrak g$ and moreover in some sense this result applied to all simple Lie algebras (more precisely except those of types $C$, $B_{2n}$ and $F_4$).

Of course these cases are rather simple and do not have the complication that the Levi factor admits a non-trivial Weyl subgroup by which the choice of $w \in W$ taking 
$\pi$ to $\Pi^*$ can be modified.

\section{The index one case.}

\subsection{}\label{10.1}
We apply our method to the index one case, that is when $\ell(\mathfrak q_{\Lambda})=1$ for the canonical truncation $\mathfrak q_{\Lambda}$ of a biparabolic subalgebra $\mathfrak q$ of a simple Lie algebra $\mathfrak g$. 
It is equivalent to the fact that $\g=\mathfrak {sl}_n(\bf k)$, $\pi^+=\pi$ and $\pi^-=\pi\setminus\{\alpha_p\}$ with $p$ and $n$ coprime by \cite{J0.5}.\par
Although this method is based on our previous paper (\cite{FJ1}) we obtain a considerable simplification of the original proof as well as a more precise result (Thm. \ref{10.3}).
\par
In this particular case, one can easily check that there is a single unmodified integer edge $E$ and also that
there is only one half-integer meander which is an edge with only one non-fictitious end.
Thus this end must be marked (\ref{3.1}) and this gives two modified integer meanders 
which may be also obtained by deleting from the whole unmodified integer meander $E$ the single arc corresponding to a simple root of $\mathfrak g$.
The latter was called in \cite{FJ1} the ``exceptional value'' and we will denote it by $\epsilon_s\beta_s$ ($s$ as simple root).
When $p=1$, one of the above (modified) edges is a trivial edge.\par
To be more precise recall that, by \cite{J0.5}, the second element $\eta$ of an adapted pair in the index one case is given by
$$\eta=\sum_{\gamma \in S} x_{\gamma}$$
with $I=[1,\,n-1]$, $S=\{\epsilon_i\beta_i\}_{i\in I\setminus\{s\}}=\mathscr K\cup(-\mathscr K')\setminus\{\epsilon_s\beta_s\}$, $\mathscr K$ (resp. $\mathscr K'$) being the Kostant cascade relative to $\pi$ (resp. to $\pi^-$) and $\epsilon_s\beta_s$ being the unique root in $\pi\cup(-\pi)$ which also belongs to $\mathscr K\cup(-\mathscr K')$.
Moreover by \cite{J0.5}, $V={\bf k}x_{\epsilon_s\beta_s}$ is an $\ad h$ stable complement of $(\ad \q_{\Lambda})(\eta)$ in $\q_{\Lambda}^*$ and then by \cite[Cor. 2.3]{JS} (see also \ref{00.2}) $\mathscr W=\eta+V$ is a Weierstrass section for $Y(\q_{\Lambda})$ in $\q_{\Lambda}^*$ in the sense of \ref{00.1}. Moreover by \cite[Thm. 8.2]{J0.7} $\mathscr W\subset (\q_{\Lambda}^*)_{reg}$.

\par
By what was explained above we obtain, by applying our present method, two modified integer meanders $E_1$ and $E_2$ and it remains to straighten them, if necessary.

Let us consider first the case when $p\neq 1$.
Straighten the non-trivial edges $E_1$ and $E_2$, 
we obtain $E_i^*=a_i^*\rightarrow b_i^*$ ($i=1,\,2$) labelled such that $a_1^*<a_2^*$. Then $a_2^*\rightarrow b_2^*\rightarrow a_1^*\rightarrow b_1^*$ gives a required 
simple root system $\Pi^*$ for $\mathfrak g$.\par
For $p=1$, only one modified integer edge say $E_1$ is non-trivial and there is one fully fixed point $f$ corresponding to $E_2$. Moreover $E_1$ is already straightened (since it has no internal turning point) that is $E_1=E_1^*$. Finally $f$ lies strictly to the right of $a=a^*=1$ thus the arc joining $f$ to $a$ attains a gap defined by $\hat\pi^-$, hence belongs to $K^-$.

\subsection{}\label{10.2}

There is a second method to obtain a simple root system $\Pi^*$ for $\mathfrak g$.  This consists of adjoining to $S$ the exceptional value $\epsilon_s\beta_s$ and applying the straightening procedure of Section \ref{6}.  By \cite [2.10]{FJ1}, $\epsilon_s\beta_s$ is never nil and is a boundary value to some unique turning point, which is internal if and only if $p\neq 1$.

 Suppose $p \neq 1$. Then the straightening procedure implies in particular that the short directed central line corresponding to $\epsilon_s\beta_s$ is replaced by a directed external arc.  Thus we obtain a new root system $\Pi^*$ such that $\epsilon_s\beta_s\in \mathbb N\Pi^*\setminus\Pi^*$.   As in \ref {9.1} (but with a slight difference) we set $I =[1,\,n-1],\, I^\prime =I\setminus \{s\}$ and
 $$y = \sum_{i \in I^\prime}x_{\epsilon_i\beta_i} + \sum_{i \in I \mid\beta^*_i \neq \epsilon_i\beta_i}x_{\beta^*_i}.$$
 
 As before by \ref {6.1}(b) the restriction of $y$ to $\mathfrak q_{\Lambda}$ is just $\eta$, whilst by \ref {6.1}(a), $y$ is conjugate to a Jordan block of size $n$ hence regular nilpotent. Furthermore $y+cx_{\epsilon_s\beta_s}: c\in \bf k$ is again regular nilpotent and has image by the restriction map $\rho\,:\,\g^*\rightarrow \q_{\Lambda}^*$, $\eta +cx_{\epsilon_s\beta_s}$.  Thus in this case the entire Weierstrass section $\mathscr W = \eta +{\bf k}x_{\epsilon_s\beta_s}$ is an image by $\rho$ of an irreducible subvariety of $\mathscr N(\mathfrak g)_{reg}$ identifying with $\bf k$.
 
 Suppose $p=1$.  When we adjoin $\epsilon_s\beta_s$ to $S$ the resulting edge is already straightened. Thus proceeding as above we find that $y+cx_{\epsilon_s\beta_s}=\sum_{i\in I'}x_{\epsilon_i\beta_i}+cx_{\epsilon_s\beta_s}\,: \,c\in \bf k$ is again regular nilpotent if and only if $c\neq 0$.  Consequently we obtain $\mathscr W \setminus \{\eta\}$ as an image by $\rho$ of an irreducible subvariety of $\mathscr N(\mathfrak g)_{reg}$ identifying with ${\bf k}\setminus \{0\}$.  In order to recover $\eta$ as an image by $\rho$ of a regular nilpotent element we ``artificially" replace $\epsilon_s\beta_s$ by $-\epsilon_s\beta_s$ and mark the latter.  This creates a turning point (at $\varphi(n-1)$) which by our general procedure causes the line representing $-\epsilon_s\beta_s$ to be replaced by an external arc representing the root $\beta_s^{\prime *}$ joining $\varphi(n)$ to $\varphi(1)$. Then as before we obtain a regular nilpotent element, namely $y'=\sum_{i\in I'}x_{\epsilon_i\beta_i}+x_{\beta_s^{\prime *}}$ mapping by $\rho$ to $\eta$.  
 
 We may summarize the above as follows
 
 \begin {lemma}  ($\mathfrak q_{\Lambda}$ a truncated biparabolic of index one of $\g=\s\l_n(\bf k)$).
 
 There is a one dimensional subvariety of $\mathscr N(\mathfrak g)_{reg}$ which maps bijectively through the restriction map $\rho:\g^*\rightarrow\q_{\Lambda}^*$ to the Weierstrass section $\mathscr W = \eta +{\bf k}x_{\epsilon_s\beta_s}$.  If $p\neq 1$ this subvariety is irreducible.  If $p=1$ it is the union of an irreducible subvariety of dimension $1$ and a point.
 \end {lemma}

 \subsection{}\label{10.3}
 
 The above result has the following interesting consequence.  First to preserve some generality let $\mathfrak a$ be an algebraic Lie subalgebra of a semisimple Lie algebra $\mathfrak g$.  Let $\overline{\textbf{A}}$ be the unique closed irreducible subgroup of the adjoint group $\textbf{G}$ of $\g$ with Lie algebra $\ad_{\g}(\mathfrak a)$. Observe that, if $\bf A$ denotes the adjoint group of $\a$, then 
any $\bf A$ orbit in $\a^*$ may be viewed as an $\overline{\bf A}$ orbit in $\a^*$ and vice-versa.
 
  The restriction map $\rho\,:\,\mathfrak g^* \rightarrow \mathfrak a^*$ is a morphism of $\overline{\textbf{A}}$ modules.  Consequently the image by $\rho$ of any $\overline{\textbf{A}}$ orbit in $\mathfrak g^*$ is an $\overline{\textbf{A}}$ orbit in $\mathfrak a^*$.   
 
 One may recall that $\mathscr N(\mathfrak g)_{reg}$ is a single $\textbf{G}$ orbit and one may ask what $\overline{\textbf{A}}$ orbits in $\mathfrak a^*$ occur in the image by $\rho$ of $\mathscr N(\mathfrak g)_{reg}$.  For example for $\a=\q_{\Lambda}$ of index one in $\g=\s\l_n(\bf k)$, with $p=1$ and $n>2$, the zero orbit in $\a^*$ is absent in the image by $\rho$ of $\mathscr N(\mathfrak g)_{reg}$, since the kernel of the restriction map $\rho$ contains no regular nilpotent element. 
 
 Yet one can still ask if $\mathfrak a^*_{reg}$ lies in the image by $\rho$ of $\mathscr N(\mathfrak g)_{reg}$, at least if $\mathfrak a$ is a proper truncated biparabolic subalgebra of $\mathfrak {sl}_n(\bf k)$.  Such a result which we have no idea how to prove would imply our main result (Thm. \ref{9.1}) though would not give a precise description of the regular nilpotent element $y$.  Nevertheless, for the index one case, Lemma \ref{10.2} gives the following

\begin{thm}
Let $\mathfrak q_{\Lambda}$ be a truncated biparabolic subalgebra of index one of a simple Lie algebra $\g$ over $\bf k$ and let $\overline{{\bf Q}_{\Lambda}}$ be the unique closed irreducible subgroup of the adjoint group $\bf G$ of $\g$ such that $Lie(\overline{{\bf Q}_{\Lambda}})=\ad_{\g}(\mathfrak q_{\Lambda})$. Then every regular $\overline{{\bf Q}_{\Lambda}}$ orbit in $\mathfrak q_{\Lambda}^*$ is the image, by restriction to $\mathfrak q_{\Lambda}$, of a $\overline{{\bf Q}_{\Lambda}}$ orbit in $\mathscr N(\g)_{reg}$.
\end{thm}

\begin{proof}
Since the index of $\mathfrak q_{\Lambda}$ is equal to one, the algebra $Y(\mathfrak q_{\Lambda})$ is a polynomial $\bf k$-algebra in one variable $f$ which is furthermore irreducible, by \cite{J0.5}.
Recall \ref{00.2}, especially the definition of the nullcone $\mathscr N(\mathfrak q_{\Lambda})$ of $\q_{\Lambda}$. Then $\mathscr N(\mathfrak q_{\Lambda})$ is the zero set of $f$ in $\mathfrak q_{\Lambda}^*$. Hence $\mathscr N(\mathfrak q_{\Lambda})$ is irreducible and by \cite[ Cor. 8.7]{J0.7}, every regular element in $\mathfrak q_{\Lambda}^*$ belongs to a $\overline{{\bf Q}_{\Lambda}}$ orbit of an element of the Weierstrass section $\mathscr W$. Finally by Lemma \ref{10.2}, each element of the Weierstrass section is the restriction to $\mathfrak q_{\Lambda}$ of an element in $\mathscr N(\g)_{reg}$. Since this restriction is a morphism of $\overline{{\bf Q}_{\Lambda}}$ modules, the conclusion of the theorem follows.

\end{proof}

\subsection{}\label{10.4}

Consider the projectivisation $\mathbb P(\mathscr N(\g))$ of $\mathscr N(\g)$ and $\mathbb P(\q_{\Lambda}^*)$ of $\q_{\Lambda}^*$.
The restriction map $\rho:\g^*\to\q_{\Lambda}^*$ induces a projective map
$\mathbb P\rho:\mathbb P(\mathscr N(\g))\setminus\mathbb P(\ker(\rho)\cap\mathscr N(\g))\to\mathbb P(\q_{\Lambda}^*)$.
Now the image by $\mathbb P\rho$ of a closed set is a closed set by \cite[Thm. 3, 5.2, Chap. I]{Scha}. Moreover  $\mathscr N(\g)$ is the closure of $\mathscr N(\g)_{reg}$ and $\mathscr N(\g)$ is conic (that is to say that, for every $c\in\bf k$ and $x\in\mathscr N(\g)$, $cx\in\mathscr N(\g)$).
Since the closure of $(\q_{\Lambda}^*)_{reg}$ is equal to $\q_{\Lambda}^*$, Theorem \ref{10.3} above gives the

 \begin{cor}
 ($\q_{\Lambda}$ a truncated biparabolic of index one of $\g$ simple.)
 
 The image of $\mathscr N(\g)$ under the restriction map $\rho\,:\,\g^*\rightarrow \q_{\Lambda}^*$ is $\q_{\Lambda}^*$.

 \end{cor} 
  
\textbf{Remark.} With respect to Theorem \ref{10.3}, the above result, even it were known for an arbitrary truncated biparabolic, is less interesting than an answer to the question raised in \ref{10.3}. Indeed it is already easy to check  that $\eta$ as given in \ref{1.8}$(*)$, with $S$ satisfying Proposition \ref{3.2}, is the image of some (not necessarily regular) element of $\mathscr N(\g)$ and for the obvious choice of the preimage of $\eta$, the conjugacy class of this preimage was calculated in \cite[Section 10]{J1}.
  
\bigskip

\begin{flushleft}
\appendix{\bf APPENDICES.}

\end{flushleft}

\section{Illustration of our method.}\label{11}

For $\g=\s\l_n(\bf k)$, and $1\le i\neq j\le n$, denote by $(ij)$ the root $\varepsilon_i-\varepsilon_j$
and by $x_{ij}$ the root vector $x_{\varepsilon_i-\varepsilon_j}$.

\begin{figure}[!h]
\centering
\input{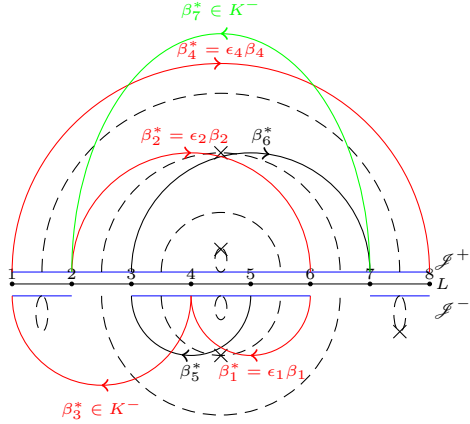}
\caption{\scriptsize{\it 
Case $\g=\mathfrak{sl}_8(\bf k)$,
$\pi^+=\pi$ and $\pi^-=\pi\setminus\{\alpha_2,\,\alpha_6\}$.}
In dashed lines are represented the half-integer meanders 
with a cross on each marked arc.
Hence we get two non-trivial modified integer meanders,
and no fully fixed point.
One edge is already straightened (drawn in black on the Figure),
and the other edge is, after straightening, drawn in red on the Figure.
In green is represented the arc joining both straightened edges.
The new system $\Pi^*$ of simple roots obtained then is
$\Pi^*=\{(53),\,(37),\,(72),\,(26),\,(64),\,(41),\,(18)\}$.
Moreover $\eta=x_{64}+x_{26}+x_{21}+x_{18}+x_{53}+x_{37}$
and $y=\eta+x_{41}+x_{72}$ is regular nilpotent in $\g^*$
and its restriction to $\q_{\Lambda}$ is equal to $\eta$.
The element $w$ of the Weyl group taking $\pi$ to $\Pi^*$
is $w=s_1s_2s_4s_6s_1s_3s_5s_2s_4s_6s_1s_3s_5s_6$.}
\label{Figure 24}
\end{figure}

\begin{figure}[!h]
\centering
\input{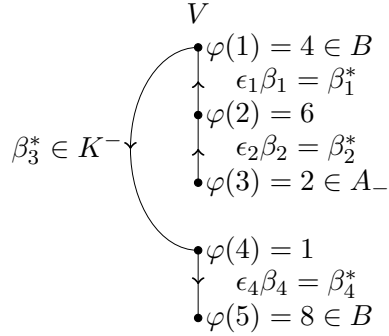}
\caption{\scriptsize{\it 
Straightening of the edge for the previous case.}
Here is represented the straightening
 of the edge $E=\{4,\,6,\,2,\,1,\,8\}$
 for the previous case.
Since $\varphi(3)$ is an internal turning point of type $A_-$,
the short directed arc representing the root $\epsilon_3\beta_3$ is deleted
and replaced by the directed external arc representing $\beta_3^*$.
Then $\Pi_E^*=\{\beta_2^*,\,\beta_1^*,\,\beta_3^*,\,\beta_4^*\}$
is a system a simple roots for $\s\l_5(\bf k)$ in which
$\epsilon_3\beta_3\in\mathbb N\Pi_E^*$ and $\beta_3^*\in K$.}
\label{Figure 25}
\end{figure}

\clearpage

\section{ Index of notations.}\label{12}

\begin{tabular}{ll}

\ref{00}&$\bf k$\\
\ref{00.0}& $\g,\,{\bf G},\,Y(\g),\,W,\,\q_{\Lambda},\,\mathscr N(\q_{\Lambda})$\\
\ref{00.1} & $\a,\,Y(\a)$\\
\ref{00.2} & $\a^f,\,\ell(\a),\,\a^*_{reg},\,{\bf A},\,(h,\,\eta)$\\
\ref{00.2} &$\mathscr N(\a),\,\mathscr N(\a)_{reg},\,,\,Sy(\a)$\\
\ref{00.3} & $\a_{\Lambda},\,\h_{\Lambda},\,\q,\,\q_{\Lambda},\,{\bf Q}_{\Lambda}$\\
\ref{00.5}& $\rho$\\
\ref{0.1}& $\mathscr J^\pm,\,S,\,\kappa^\pm,\,S^\pm,\,\sigma^\pm$\\
\ref{0.2}&$S_E,\,E,\,\eta_E,\,y_E,\,E^*$\\
\ref{0.3}&$I,\,K,\,\iota_I$\\
\ref{0.4} (\ref{4.3})&$\epsilon_I$\\
\ref{0.5}&$\Pi_E^*$\\
\ref{0.6} (\ref{8.2})&$\iota_{E^*}$\\
\ref{1.1}&$\varepsilon_i,\,[r,\,s],\,c_{[r,\,s]},\,L,\,(r,\,s),\,\alpha_i,\,\Delta^\pm$\\
\ref{1.1}&$\Delta,\,\pi,\,{\rm supp}(\,),\,x_{\gamma}$\\
\ref{1.3}&$e,\,\varphi$\\
\ref{1.4}&$T,\,T_0$\\
\end{tabular}
\begin{tabular}{ll}

\ref{1.5}&$I_{s,\,t},\,\beta_i,\,\iota_I,\,{\rm Supp}(\iota_I),\,\Phi(I)$\\
\ref{1.6}&$\kappa_J$\\
\ref{1.7}&$\mathscr J,\,\kappa_{\mathscr J},\,\mathscr K,\,\kappa_{\mathscr J^\pm}\, (\kappa^\pm),\pi_i^\pm,\pi^\pm$\\
\ref{1.8}&$\q_{\pi^+,\pi^-},\,\q_{\mathscr J^+,\mathscr J^-},\,\hat\pi^{\pm},\,S$\\
\ref{1.8}&$R,\,K,\,R_*,\,M,\,R_*^\pm,\,K^\pm$\\
\ref{2.1}&$\mathscr M,\,\mathscr K,\,,\overline{\pi},\,\mathscr M_{r,\,s},\,\kappa_{\mathscr M_{r,\,s}},\,\sigma_{\mathscr M_{r,\,s}},\,\sigma$\\
\ref{2.2.1}&$A_{r,\,s},\,A'_{r,\,s}$\\
\ref{2.4}&$\sigma^\pm$\\
\ref{3.1}&$\iota^\pm,\,\widetilde\pi$\\
\ref{3.3}&${\bf E}(S),\,S_E,\,\Pi,\,\epsilon_i$\\
\ref{3.6}&$\overline{\pi}^\pm$\\
\ref{3.8.4}&$\mathscr I,\,\mathscr E,\,F,\,M_f$\\
\ref{6.1}&$\Pi^*,\,\beta_i^*$\\
\ref{6.2}&$V$\\
\ref{6.4}&$A,\,B$\\
\ref{6.5}&$sg(t),\,A_{\pm},\,B_{\pm}$

\end{tabular}

\section{ Index of notions.}

\begin{tabular}{ll}
adapted pair & \ref{00.2}\\
 anti-Toeplitz involution & \ref{2.1}\\
arrows aligned & \ref{1.4}\\
assigned to : marked value (resp.  directed external arc) assigned &\ref{6.4}, resp. \ref{6.5}\\
batch of arcs & \ref{2.2.1}\\
 biparabolic subalgebra &\ref{1.8}\\
 boundary (point or value) &\ref{4.2}\\
canonical truncation &\ref{00.3}\\
component (of the double partition) or double component &\ref{1.9}\\
compound (interval or interval value) &\ref{1.5}\\
cross a point (an arc crosses a point ... ) &\ref{6.5}\\
directed arc &\ref{1.1}\\
directed external arc &\ref{6.2}\\
edge &\ref{1.3}\\
end-point (of an edge, resp. of an arc)&\ref{1.4}, resp. \ref{1.1}, \ref{6.2}\\
equicentral (double) component & \ref{3.4.5.4}\\
finishing point (of a directed arc, resp. an edge)&\ref{1.1}, resp. \ref{1.3}\\
\end{tabular}

\begin{tabular}{ll}

fully fixed point &\ref{8.5}\\
gap &\ref{4.1}\\
half-integer arc &\ref{3.1}\\
half-integer meander &\ref{3.1}\\
internal turning point &\ref{1.4}\\
interval, interval value &\ref{1.5}\\
 isolated value &\ref{6.3}\\
join a point (an arc joins a point) &\ref{1.1}\\
Kostant cascade &\ref{2.1}\\
linked components (of the double partition) &\ref{3.4.5.4}\\
loop &\ref{1.3}\\
marked value &\ref{6.4}\\
marking &\ref{2.1}\\
meander &\ref{1.3}\\
meet a point (an arc meets a point) &\ref{1.1}\\
meet simply a point &\ref{6.5}\\
modified (unmodified) integer arc &\ref{2.4}\\
modified (integer) meander &\ref{2.4}\\
modified involution &\ref{2.4}\\
nil (element-value-arc-root) &\ref{4.2}\\
nullcone&\ref{00.2}\\
odd compound interval (interval value) &\ref{4.3}\\
partition &\ref{1.7}\\
short directed arc &\ref{6.2}\\
signature (positive or negative) &\ref{6.5}\\
simple (interval or interval value) &\ref{1.5}\\
sink &\ref{1.4}\\
source &\ref{1.4}\\
starting point (of a directed arc, resp. of an edge) &\ref{1.1}, resp. \ref{1.3}\\
straightened edge &\ref{6.6}\\
straightened edge value&\ref{8.2}\\
trivial arc &\ref{1.1}\\
trivial component &\ref{3.7}\\
trivial (non-trivial) edge &\ref{1.3}\\
 trivially (non-trivially) modified integer arc &\ref{2.4}\\
turning point &\ref{1.4}\\
value of an edge &\ref{1.3}\\

\end{tabular}

\bigskip

\bibliographystyle{elsarticle-num}

\end{document}